\documentclass[10pt]{article}

\usepackage[english]{babel}

\usepackage[letterpaper,top=2cm,bottom=2cm,left=3cm,right=3cm,marginparwidth=1.75cm]{geometry}

\usepackage[colorlinks=true, allcolors=blue]{hyperref}

\usepackage{amsmath,amsthm,amsfonts,amssymb,amscd}
\usepackage{stmaryrd}
\usepackage{enumerate}
\usepackage{mathrsfs}
\usepackage{longtable,geometry}
\usepackage{mathtools}
\usepackage{yhmath}
\usepackage{pgf,tikz,pgfplots}
\pgfplotsset{compat=1.15}
\usetikzlibrary{arrows}
\usepackage{dsfont}

\usepackage{enumerate}
\usepackage{xcolor}
\usepackage{graphicx}
\usepackage{listings}
\usepackage{hyperref}
\usepackage{stmaryrd}
\usepackage{verbatim}
\usepackage{subcaption}

\usepackage{appendix}

\newcommand{\R}{\mathbb{R}}

\newcommand{\Haus}{\mathcal{H}^{2}}

\newcommand{\Ccal}{\mathcal{C}}

\DeclareMathOperator*{\argmin}{arg\,min}
\newcommand{\Hausone}{\mathcal{H}^{1}}
\newcommand{\bdy}{\Gamma}
\newcommand{\bpar}{\gamma}
\newcommand{\ra}{\rightarrow}
\newcommand{\eps}{\varepsilon}
\newcommand{\curve}{\mathrm{Hom}}
\newcommand{\mres}{\mathbin{\vrule height 1.6ex depth 0pt width 0.13ex\vrule height 0.13ex depth 0pt width 1.3ex}}
\newcommand{\surf}{S}
\newcommand{\AT}{\mathcal{AT}_{\eps}}

\title{Phase field approximation for Plateau's problem: a curve geodesic distance penalty approach}
\author{Matthieu Bonnivard \and \'Elie Bretin \and Antoine Lemenant \and Eve Machefert}

\newtheorem{thm}{Theorem}[section] 
 
\newtheorem{definition}{Definition}[section] 
\newtheorem{prop}{Proposition}[section]
\newtheorem{cor}{Corollary}[section]
\newtheorem{rmq}{Remark}[section]

\newtheorem{lemme}{Lemma}[section]

\begin{document} 
	\maketitle
	
	\begin{abstract}
		This work focuses on a phase field approximation of Plateau's problem. Inspired by Reifenberg's point of view, we introduce a model that combines the Ambrosio-Torterelli energy with a geodesic distance term which can be considered  as a  generalization of the approach developed in~\cite{lemSant,bonnivard2015approximation} to approximate solutions to the Steiner problem. First, we present a $\Gamma$-convergence analysis of this model in the simple case of a single curve located on the edge of a cylinder. In a  numerical section, we  detail the numerical optimisation schemes used to minimize this energy for numerous examples, for which good approximation of solutions to Plateau's problem  are found.
	\end{abstract}

	\vspace{1cm}
	\tableofcontents

	\section{Introduction}

	\indent In this paper we introduce a phase field approximation of  the famous Plateau's problem, which consists in finding a surface of minimal area spanning a given set of curves. The precise formulation of this topological constraint raises challenging mathematical questions that contributed to the development of the geometric measure theory, from which different major schools of thought have emerged.\\
	
	On the one hand, Federer and Fleming introduced objects known as integral currents to describe the boundary constraint, in the case of oriented surfaces~\cite{FF60}. They established the existence of an integral current with fixed boundary and minimal mass, which can be interpreted as a weighted area. On the other hand, Reifenberg developed  a different approach using  \v{C}ech homology to characterize the property of spanning a curve~\cite{reifenberg1960}. He obtained the existence of surfaces of minimal area in a class including non-orientable surfaces, as well as those studied by Federer and Fleming. At the same time De Giorgi developped the theory of sets of finite perimeter which provides  a powerful  approach in the co-dimension one case  (see~\cite{degiorgi1,degiorgi2} and the books~\cite{maggi_sets_2012,giusti}). More recently, some new approaches have been introduced as the one   of Harrison and Pugh~\cite{harrison2014soap,harrison2016existence} using a very natural  spanning criterium together with an existence result.    Another nice framework is  the very general   existence theory  in~\cite{de2017direct,maggi2023plateau}. Alternatively, the notion of sliding minimizers of~\cite{davidSliding} attempts to give yet a better description of what a soap film should be.\\

	Numerous numerical methods have been developed to approximate minimal surfaces. In references~\cite{MR1502829, MR0458923, MR942778, MR1613695, MR1613699}, the surface is represented parametrically and then discretized using a finite element approach.  Numerical techniques involving an implicit representation of minimal surfaces have been developed, for example, in references~\cite{MR1214016} and~\cite{MR2143330}, in which a level set method is employed. Concerning phase field methods, the approach proposed in~\cite{Wang:2021:CMS} offers  particularly suitable numerical approximations, but is restricted to oriented surfaces. It is also worth mentioning that a phase field approximation for Plateau's problem  was already proposed  by Chambolle, Merlet and Ferrari in~\cite{Merlet1, Merlet2}, or more precisely for the $\alpha$-Mass energy, which corresponds to Plateau's problem in the limit $\alpha\to 0$. In these works,  the authors start with another interpretation of the    Steiner approximation introduced  in~\cite{lemSant}, by looking at a dual formulation for the geodesic distance. This allows them to transform the topological constraint into a divergence constraint on the second variable, and gave rise to the paper~\cite{Merlet1} which focuses on the 2D-Steiner problem. It was then extended by the same authors for any dimensions and co-dimensions in~\cite{Merlet2}, but the topological constraint in that paper is again a divergence constraint on the second variable that can be seen as a boundary constraint on a $k$-current.   Let us also mention the paper~\cite{maggi2023hierarchy}, in which the authors study the $\Gamma$-limit of a phase field approximation but where the topology is prescribed as a constraint in the level set of the field.  The paper~\cite{bretin2024penalizedallencahnequationmean} also proposes a numerical investigation of Plateau's problem  based on Allen–Cahn's equation and a penalty term localized on the  skeleton of the interface in order to control its topology during the computation. \\

	In this paper,  we propose a new phase field model to approximate a certain instance of Plateau's problem, that more closely resembles Reifenberg's example than the formulation with currents proposed in~\cite{Merlet2}. 
	In particular, the present work also contributes to the challenging question of prescribing the topology of the level set of a function.  For that purpose we use a penalty term involving a certain geodesic distance, weighted by $u$, in the spirit of what was done for the Steiner problem in~\cite{lemSant,bonnivard2015approximation,bonnivard2018phase,BBL2020}. For the recall, the classical Steiner problem consists in finding the shortest path connecting a given finite set of points, and can be interpreted as a Plateau problem in lower dimension.  In this context, a geodesic distance has been used  to enforce the connectedness constraint on the level set of the phase field $u$. On the other hand the elliptic energy of Cahn-Hilliard type on  $u$ is the quantity that will encode the surface area. In this paper we aim to generalize this approach to more complicated topological constraints in higher dimensions.   Our approach is  different from~\cite{maggi2023hierarchy}, in the sense that we add a penalty term that prescribes the topology only at the limit, which is more relevant in view of performing numerical results. Actually, showing that this penalty term is indeed able to prescribe the topology is exactly the core of our results.
	\\
	
	{\bf A curve geodesic phase field model for Plateau's problem.} Let us now describe how to extend the mathematical formulation of the phase field approximation of Steiner's problem to Plateau's problem. In the rest of the paper, we will always consider $2$-dimensional surfaces in $\R^3$, since it corresponds to the physical observations of soap films. However, we point out that the model and the analysis performed in Section~\ref{Section:Analysis} can be generalized in codimension $1$ for any dimension.

	The approach developed  in~\cite{lemSant,bonnivard2015approximation,bonnivard2018phase, BBL2020} for Steiner's problem over a set $\Omega\subset \R^2$, consists in adding to the classical Ambrosio-Tortorelli length approximation term 
	\[
	\AT(u) :=\int_{\Omega} \left(\eps |\nabla u|^{2}  + \frac{1}{4\eps}(1-u)^{2}\right)dx,
	\]
	a penalization term involving a geodesic distance to enforce the topological constraint of Steiner's problem. The energy to minimize reads
	\begin{equation} \label{SteinerEnergy} 
		\AT(u)   + \frac{1}{c_{\eps}}\sum_{i=1}^{n} d_{u,\eps}(a_0, a_{i}),\end{equation}
	where $d_{u,\eps}(a_0, a_{i})$ is the geodesic distance between the points $a_0$ and $a_i$, with weight $u^2+\delta_\eps$, defined as follows:
	\[d_{u,\eps}(a_0, a_{i}) := \inf \left \{ \int_\gamma (u^2 +\delta_\varepsilon)  d\mathcal{H}^1, \gamma \text{ Lipschitz curve connecting } a_0 \text{ to } a_i \right  \}.\]
	In the previous definition, $(\delta_\eps)$ is a sequence of positive numbers converging to zero, whose role in the modeling is to control the length of the geodesics associated with the energy~\eqref{SteinerEnergy}, and gain compactness.

	The goal of this paper is to propose and justify a variational approximation of the form~\eqref{SteinerEnergy} for Plateau's problem in $\R^3$, where the Ambrosio-Tortorelli term $\AT(u)$ naturally approximates the area of the sought surface, and the penalized term enforces the topological constraint of spanning a set of curves.
	By analogy with Steiner's problem, where geodesic curves connect given points, optimal surfaces for Plateau's problem will connect fixed closed curves.

	More precisely,  we intend  a set  $E$ to be spanned by a boundary $\Gamma$, when  any pair of  closed curves $\gamma_1,\gamma_2$ contained  in the boundary $\Gamma$ can be joined continuously by an homotopy inside $E$. This is inspired by Reifenberg's topological assumption that any element of the first group of homology of $\Gamma$ should be trivial in the homology group of $E$. By interpreting these homotopies  connecting two given closed curves as a geodesic distance between the two curves, we obtain a very natural analogy between Steiner's approximation and Plateau's approximation. 
	
	To explain this fact more in detail, let us define as follows, the set of admissible homotopies
	\[
	\curve(\gamma_i,\gamma_j) : = \{ \ell \in \mathrm{Lip}([0,1] \times \mathbb{S}^{1}, \R^3) \text{ such that } \ell(0) = \gamma_i \text{ and } \ell(1) = \gamma_{j} \}.
	\]
	For the sake of clarity, in the sequel, we will use the term ``curve" only to refer to the given curves associated with the topological constraint of Plateau's problem. We shall employ the term ``homotopy" (or ``path") to identify an element of $\curve(\gamma_i,\gamma_j)$.

	For all paths $\ell \in \curve(\gamma_i,\gamma_j)$, we define the surface $\surf_\ell$ as the image
	\begin{equation}\label{Def:surface_image}
		\surf_{\ell} := \ell([0,1] \times \mathbb{S}^{1}).
	\end{equation}
	Notice that the surface $\surf_{\ell}$ is $\mathcal{H}^{2}$- rectifiable and has finite $\mathcal{H}^{2}$-measure.


	
	We are now in position to define the generalized geodesic distance between curves, associated with a given continuous function $u$.
	\begin{definition}\label{Def:geoDistance}
		Let $p \in [1,+\infty]$ and $(\delta_\eps)_{\eps>0}$ be a sequence of positive numbers that converges to zero. We define the $p$-geodesic distance between $\gamma_i$ and $\gamma_j$ as follows:
		\begin{align*}
			&d^{p}_{u,\eps}(\gamma_i,\gamma_j) := \inf\left \{  \int_{\surf_{\ell}}( |u|^{p} + \delta_\eps) d\Haus \ |\ \ell \in \curve(\gamma_i,\gamma_j)\right \}  \text{ if } p<\infty, \\
			&d^{\infty}_{u,\eps}(\gamma_i,\gamma_{j}) := \inf \left \{ \sup_{\surf_{\ell}} (|u|+ \delta_\eps)\ |\ \ell \in \curve(\gamma_i,\gamma_j)   \right \}.\\
		\end{align*}
	\end{definition}
	
	Subsequently, for the sake of simplicity, we will use the notation $d^{p}_{u}$ instead of $d^{p}_{u,\eps}$. We can now introduce the approximation functional for Plateau's problem using this approach, with the given closed curve $\gamma_0,\ldots, \gamma_d$ contained in an open set $\Ccal$, as follows:
	\begin{equation}
		\label{functionalGeneral}
		F^{p}_{\eps}(u) := \eps \int_{\mathcal{C}}|\nabla u|^{2} dx + \frac{1}{4\eps}\int_{\mathcal{C}}(1-u)^{2}dx + \frac{1}{c_{\eps}}\sum_{i=1}^{d}d^{p}_{u}(\gamma_0, \gamma_{i}) .
	\end{equation}
	In the above definition, $(\delta_\eps)$ and $(c_\eps)$ are a sequences of positive numbers that converge to zero, such that the ratio $\delta_\eps/c_\eps$ converges to zero.
	The phases $u$ are $H^1$ functions with values in $[0,1]$, that we also assume to be continuous. They will satisfied extra conditions, defined in the next Section, to ensure that the topological constraint is fulfilled. \\

	{\bf Main theoretical result.} Our main result states that from a sequence of minimizers of $F_\varepsilon^p$, we can construct a solution of Plateau's problem. More precisely, this solution is obtained as the limit of connected components of the complementary set of a level set $\{u_\varepsilon \leqslant t_\varepsilon\}$. For the sake of simplicity, we will prove this result in the particular case of a single curve, contained in the boundary of a cylinder. In this particular situation we obtain the following Theorem.

	\begin{thm}
		\label{th : csq Gamma cv}
		Let $p\in [1,+\infty]$, $d=1$ (the constant in \eqref{functionalGeneral}), and  $\gamma_1$ be a single   curve which is the graph of a Lipschitz function lying on the boundary of  a cylinder. We also let   $\gamma_0$ be a constant   curve inside the support of $\gamma_1$ (as precisely described in Section \ref{defprob}). Then, for any quasi-minimizing sequence $u_{\eps}$ for $F_\eps^{p}$ (in the sense of Definition~\ref{def quasi min}), and setting $g_\eps:=u_\eps-{u_\eps^2}/{2}$, there exists $s_\varepsilon=O(c_\varepsilon)$ and  $t_{\varepsilon} \in [s_{\varepsilon},1/2-\varepsilon]$ such that the connected component of  the level set $\{g_\varepsilon >t_\varepsilon\}$ containing the upper part of the cylinder,  converges in $L^1$ (up to a subsequence)  to a solution of Plateau's problem (as defined in \eqref{def:pblimit}).
	\end{thm}
	
	The proof of Theorem \ref{th : csq Gamma cv} follows from a $\Gamma$-convergence type analysis, combining a limsup and a liminf inequality. The proof of the limsup inequality is rather standard, and relies on constructing the appropriate recovery sequence using the 1D-profile. The more delicate part is for the liminf inequality, for which the argument used in~\cite{bonnivard2015approximation} does not work. Indeed, for the Steiner problem in~\cite{bonnivard2015approximation} the liminf was achieved by use of an integral geometric functional involving the length of the projection over lines in every direction, that was able to bound the total length of the union of the geodesic curves. This argument was purely 1D and here in higher dimensions,  the analogue integral geometric functional cannot be used in the same manner. We therefore have developed a completely different argument, based on the co-area formula and on the selection of certain good level sets of $u_\eps$. The main issue is to show that those level sets satisfy the desired topological constraint,  and that the convergence is robust enough to preserve the topology at the limit. \\

\begin{figure}[!htbp]
		\includegraphics[width=.32\textwidth]{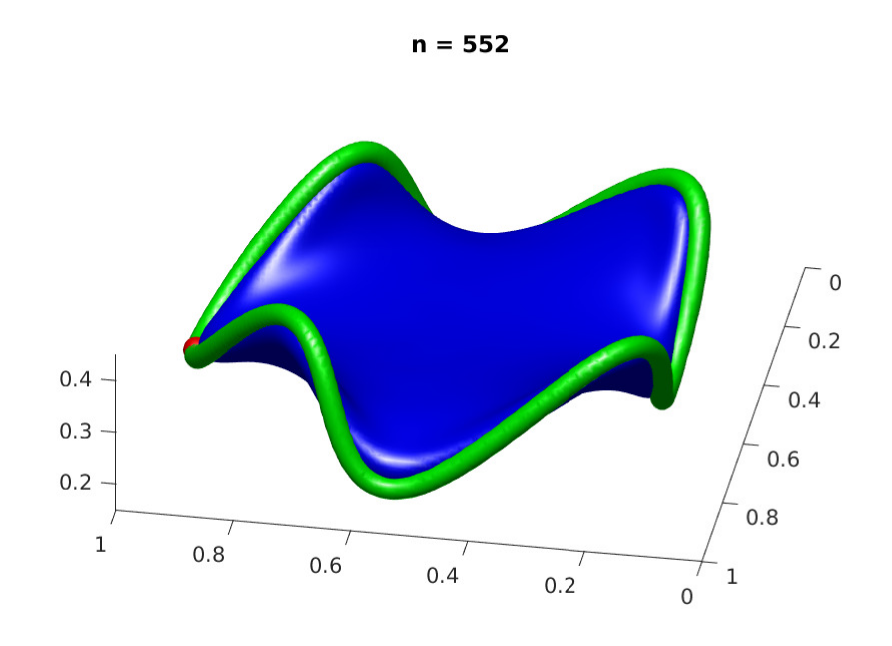}
		\includegraphics[width=.32\textwidth]{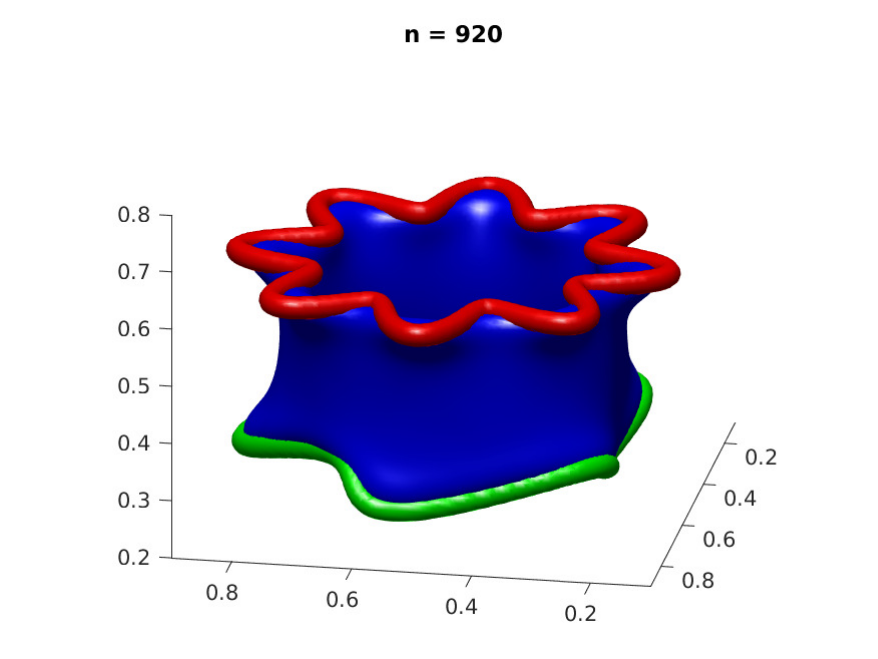}
		\includegraphics[width=.32\textwidth]{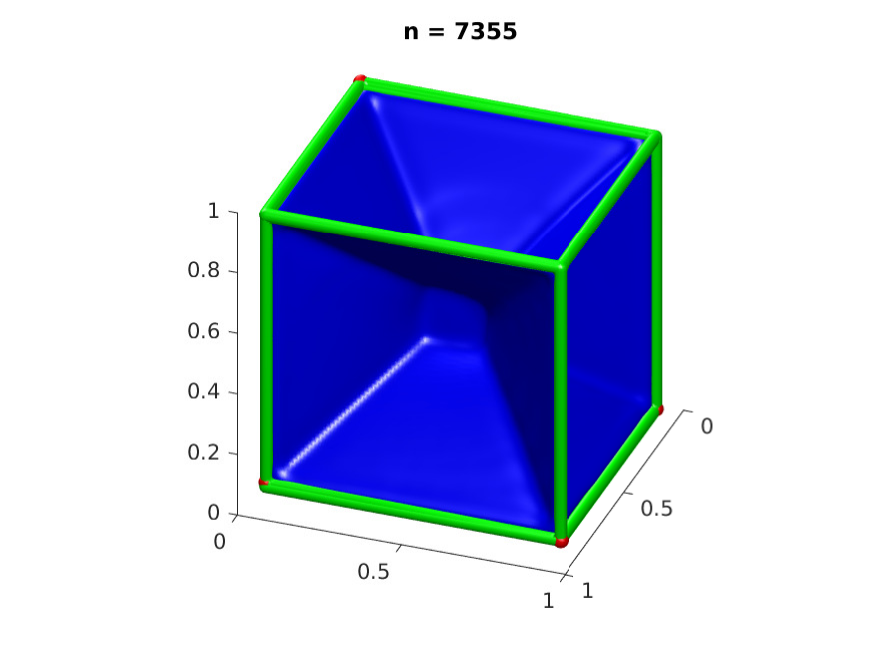}
		
		\caption{Examples of numerical approximations of the solution of the Plateau problem   }
		\label{fig:plateau_cube}
	\end{figure}

	{\bf Numerical experiments.}  The objective of the numerical part is twofold. First, it aims to demonstrate the effectiveness of our approach for computing  approximations of solutions to Plateau's problem in a framework that goes well beyond the scope of $\Gamma$-convergence proofs, taking into account multiple curves, singularities, and non-oriented solutions.  For example, Figure~\ref{fig:plateau_cube} shows three numerical approximations of solutions to Plateau's problem. These include the well-known cube example, which has a singular solution. The second objective is to propose an efficient scheme to minimize the functional $F^2_{\eps}$ by considering a discretization of it's $L^2$-gradient flow following the relaxation proposed in~\cite{BBL2020}. We also consider a variant consisting of replacing the Ambrosio-Tortorelli length approximation with the Willmore-Cahn-Hilliard energy recently introduced in~\cite{bretin2024} in order to improve the regularity of the phase field function $u$ and accelerate the optimisation process. The minimization of the geodesic term is also a delicate point that requires some approximation to use fast-marching algorithms  as in the case of the Steiner problem. \\
	
	{\bf Outline of the paper.} The rest of the paper is organized as follows. In Section 2, we justify precisely the convergence result outlined above in the simplified setting of a single curve prescribed in a cylinder. In Section 3, we detail a numerical approach to justify the relevance of this approximation and provide a proof of concept of the method in a wide class of examples.

	\medskip
	
	{\bf Acknowledgments.} The authors would like to thank Laurent Mazet for fruitfull discussions  about the proof of Proposition \ref{th : 2.1} and for giving us nice references on minimal surfaces. We also warmly thank Camille Labourie for providing us with the reference  about Borsuk Theorem and for many ``coffee-discussions" about this work.  This work was partially supported by the ANR Project ``STOIQUES'' and by the Institut Universitaire de France (IUF).


	\section{Gamma-convergence of the approximate functional}~\label{Section:Analysis}
	
	To carry out the convergence analysis of the phase field approximation that we propose, we consider the simplified setting where the boundary of the minimal surface is composed of only one closed curve, contained in the side edge of a cylinder. As will be detailed in the next subsection, these geometric assumptions will allow us to exploit the theory of sets with finite perimeter to get existence of a solution to the limit problem. They will also imply separation properties that will be at the core of the level set method developed in the proof of the $\liminf$ result.
	

	\subsection{Definition of Plateau's problem}
	
	\label{defprob}

	We choose to work in a cylinder to guarantee the separation property described in detail in Subsection~\ref{Sect:separation}. Although the results remain valid for any cylinder with convex basis, we assume this basis to be a disk for the sake of simplicity. This open cylinder is denoted by $\Ccal_0$, and is assumed to be of radius $1$, centered at the origin and of height $2h>0$. 
	
	The prescribed closed curve $\Gamma$   is assumed to be contained in the lateral edge of $\Ccal_0$. Indeed, we assume that $\Gamma$ is a Lipschitz graph and denote by $\tilde{\gamma}$ the Lipschitz function such that $\Gamma = \{(x',\tilde{\gamma}(x')),\ x' \in \mathbb{S}^{1}\}$. This setting will allow us to impose the boundary constraint in Plateau's problem using a radial extension of $\Gamma$ outside $\Ccal_0$, from which we will construct two open sets $\mathcal D^\pm$ defined as follows.
	
\begin{figure}
		\centering
		\subfloat[\centering Dilated cylinder]{\includegraphics[scale = 0.36]{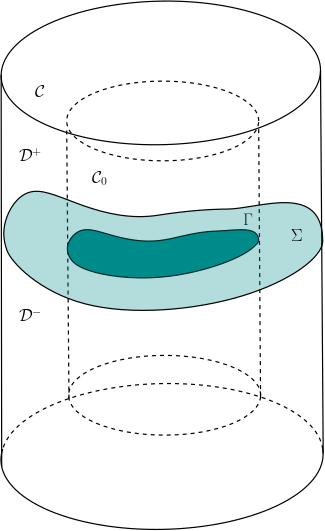} \label{fig:2_cylindres_emboites}} \quad \quad \quad \quad \quad
		\subfloat[\centering Section of the nested cylinders]{\includegraphics[scale = 0.3]{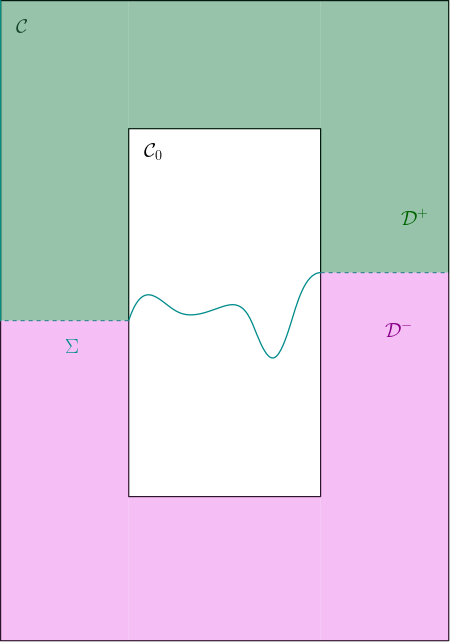} \label{fig:3_cylindres_emboites}} 
		\caption{Geometrical setting of Plateau's problem. The boundary constraint is defined through a radial extension of a prescribed curve $\Gamma$ contained in the lateral edge of a cylinder $\Ccal_0$.}
		\label{fig:cylindes_emboites}
	\end{figure}

	We denote by $\mathcal{C}$ the dilatation of cylinder $\mathcal{C}_0$ by a factor $\lambda > 1$, in other words $\mathcal{C}=\lambda \mathcal{C}_0$, as shown in Figure \ref{fig:2_cylindres_emboites}, so that $\overline{\mathcal{C}_0} \subset \mathcal{C}$. Analogously, we denote by 
	$\Sigma$ the radial extension of $\Gamma$ in $\Ccal$, which consists in the surface described in cylindrical coordinates by
	\[\Sigma := \{(r,\omega, \tilde{\gamma}(1,\omega)),\ r\in ]1,\lambda[ \text{ and } \omega \in \mathbb{S}^1\}\subset \Ccal.\]
	This surface $\Sigma$ separates $\Ccal\setminus \overline{\mathcal{C}_0}$ in two open connected components $\mathcal{D}^{+}$ and $\mathcal{D}^{-}$, located respectively above and below, and defined by
	\begin{align*}
		&\mathcal{D}^{+} := \{(r,\omega,z) \in \Ccal\setminus \overline{\mathcal{C}_0}, z >\tilde{\gamma}(\omega) \},\\
		&\mathcal{D}^{-} := \{(r,\omega,z) \in \Ccal\setminus \overline{\mathcal{C}_0}, z < \tilde{\gamma}(\omega) \}.
	\end{align*}
	We represent the geometrical setting in Figure~\ref{fig:cylindes_emboites}.


	Since the prescribed boundary is composed of only one curve, denoted by $\gamma := (Id_{\mathbb{S}^1},\tilde{\gamma})$, in order to use the geodesic distance between curves as in Definition~\ref{Def:geoDistance}, we introduce a constant curve $\gamma_0$ associated to an arbitrary fixed point $x_{0} \in \Gamma$. Upon translating $\Gamma$ vertically, we can assume that $x_0$ is at height $0$.



	Having introduced the objects required to define the boundary constraint, we are now in position to define the competitors for the problem.
	
	\begin{definition}\label{Def:competitor}
		We say that a Borel set $\Omega\subset \Ccal$ is a \emph{competitor} if its characteristic function $\chi_{\Omega}$ has bounded variation in $\Ccal$, $\mathcal{D}^{+} \subset \Omega$ and $\mathcal{D}^{-} \subset \Omega^{c}$. 
	\end{definition}

	\begin{definition}
		The Plateau problem that we consider is defined as
		\begin{equation}
			\label{def:pblimit}
			\inf \{ P(\Omega ,\mathcal{C}) ,\ \Omega \text{ is a competitor}\}.
		\end{equation} 
	\end{definition}
	The existence of a minimizer follows from the   theory of sets with finite perimeter, and the direct method of calculus of variations. Moreover, in this particular simplified case we also know from the standard theory of minimal surfaces  that the minimizer is actually a graph inside the cylinder $\Ccal_0$, that coincides with $\Gamma$ on the boundary of $\Ccal_0$. In the following statement we gather the known properties of the minimizer that will be used later.
	
	\begin{prop}
		\label{th : 2.1}
		Plateau's problem \eqref{def:pblimit} admits a solution $\Omega$ which is equivalent to an open set such that $\Haus$-a.e. $\partial^*\Omega=\partial \Omega$ (where $\partial^*\Omega$ denote the essential boundary of $\Omega$). Moreover $\partial \Omega\cap \Ccal_0$ is a smooth graph which is locally analytic in $\Ccal_0$, and $\partial \Omega \cap \partial \Ccal_0 = \Gamma$.   In particular, $\partial \Omega \cap \Ccal_0$ is of the form $S_{\ell}$ for some $\ell \in {\rm Hom}(\gamma, \gamma_0)$, and $S_\ell$ is Ahlfors regular ``up to the boundary $\Gamma$'' in the sense that  there exist some constants $C_{1},C_{2} > 0$ such that for all $x \in \partial\Omega \cap \Ccal$ and $r>0$ such that $B(x,r)\subset \Ccal$,
		\begin{equation}
			\label{th : 2.2}
			C_{1} \leqslant \frac{P(\Omega,B(x,r))}{r^{2}} \leqslant C_{2}.
		\end{equation}
	\end{prop}
	
	\begin{proof} The existence of a set of finite perimeter $\Omega$ that is a solution for the Problem \eqref{def:pblimit} follows from the  direct method of calculus of variations (see for instance~\cite[Proposition 12.29]{maggi_sets_2012} or~\cite{Machefert2025}). The fact that this minimizer   is equivalent to an open set such that $\partial^*\Omega=\partial \Omega$ $\Haus$-a.e., together with the Ahlfors-regularity (up to the boundary) is one of the main results of~\cite{Machefert2025}. It is also known from the regularity theory of perimeter minimizers, that $\partial \Omega$ is locally an analytic surface inside $\Ccal_0$ (see~\cite[Theorem 27.3]{maggi_sets_2012}). It remains to prove that $\partial \Omega$ is a graph that coincides with the curve $\Gamma$ on the boundary of $\Ccal_0$, in other words that  $\partial \Omega \cap \partial \Ccal_0 = \Gamma$. This will directly imply that $\partial \Omega \cap \Ccal_0= S_{\ell}$ for some $\ell \in {\rm Hom}(\gamma, \gamma_0)$. This claim actually follows from the theory of minimal surfaces. Indeed, it was proved by Rado~\cite{Rado}  (see also~\cite[Theorem 16 page 94]{lawson}) that when a Jordan curve $\Gamma$ admits a one-to-one orthogonal projection onto a convex Jordan curve in the plane (which is clearly the case here for our choice of $\Gamma$, that is a graph above a 2D-circle), then  there is a unique solution to the minimal surface equation satisfying the Dirichlet condition $\Gamma$. Moreover, this solution can be expressed as a graph. Now we claim that our BV-solution $\Omega$ must coincide with this minimal graph solution. There are at least two different ways to prove it: the first one is to use the graph solution as a barrier, to touch $\partial \Omega$ from above and below and then invoke the unique continuation property for the minimal surface equation.  The second way to see it is to use  that  a minimal graph admits a  calibration, obtained by extending the normal to the surface. 
	\end{proof}


	\subsection{Statement of the main results}
	
	We recall the Definition of the approximation functional in the case of one closed curve. 
	
	\begin{definition}
		Let $p \in [1,+\infty] $, $(c_{\eps})$ and $(\delta_{\eps})$ sequences of positive numbers converging to zero such that $\delta_\eps / c_\eps \to 0$, and $ u \in H^{1}(\Ccal)\cap C(\overline{\Ccal})$ such that $0\leqslant u \leqslant 1$ and $u = 1$ on $\partial{\Ccal}$. We define the energy functional 
		\begin{equation}
			\label{functional}
			F^{p}_{\eps}(u) := \eps \int_{\mathcal{C}}|\nabla u|^{2} dx + \frac{1}{4\eps}\int_{\mathcal{C}}(1-u)^{2}dx + \frac{1}{c_{\eps}}d^{p}_{u}(\gamma, \gamma_0) .
		\end{equation}
	\end{definition}

	The $\Gamma$-convergence type result that we prove in this paper is divided, as usual, into a $\Gamma-\limsup$ and a $\Gamma-\liminf$ inequality.

	\begin{rmq}
		In the following two statements we will use the notation $P(\Omega, F)$ for the relative perimeter  of $\Omega$ in the closed set $F$. This means $\mathcal{H}^{N-1}(\partial^* \Omega\cap F)$, which is well defined for any set $\Omega$ that has finite perimeter in $\R^N$.
	\end{rmq}

	\begin{thm}[$\Gamma-\limsup$]\label{Thm:limsup}
		Let $p \in [1,+\infty] $ and $\Omega \subset \Ccal$ be a competitor for the Plateau problem~\eqref{def:pblimit}, in the sense of Definition~\ref{Def:competitor}. 
		Then, there exists a sequence $(u_{\eps}) \in H^{1}(\Ccal)\cap C(\overline{\Ccal})$ such that $0\leqslant u_{\eps} \leqslant 1$, $u_{\eps} = 1 $ on $\partial{\Ccal}$ and
		\begin{equation}\label{Ineq:limsup}
			\limsup_{\eps \to 0} F^{p}_{\eps}(u_{\eps}) \leqslant P(\Omega,\overline{\mathcal{C}_0}).\end{equation}
	\end{thm}

	\begin{thm}[$\Gamma-\liminf$]
		\label{Thm:liminf}
		Let $p \in [1,+\infty]$ and $\Omega$ be a solution to the Plateau problem~\eqref{def:pblimit}.
		For any sequence $u_{\eps} \in H^{1}(\Ccal)\cap C(\overline{\Ccal})$ such that $0\leqslant u_{\eps} \leqslant 1$ and $u_{\eps} = 1$ on $\partial \Ccal$,
		\begin{equation}\label{Ineq:liminf}
			\liminf_{\eps \to 0} F^{p}_{\eps}(u_\eps) \geqslant P(\Omega, \overline{\mathcal{C}_0}).
		\end{equation}
	\end{thm}
	
	As already mentioned in the Introduction, by combining  these Theorems we obtain Theorem \ref{th : csq Gamma cv} which is a classical consequence of a $\Gamma$-convergence type result, and justifies the relevance of this approach. Here below is the definition of a quasi-minimizing sequence that was used in the statement of Theorem \ref{th : csq Gamma cv}. 
	
	\begin{definition}[Quasi-minimizing sequence for $F_\eps^p$]
		\label{def quasi min}
		Let $p\in [1,\infty]$. We say that a sequence $u_{\eps} \in H^{1}(\Ccal)\cap C(\overline{\Ccal})$ such that $0\leqslant u_{\eps} \leqslant 1$ and $u_{\eps} = 1$ on $\partial\Ccal$ is a quasi-minimal sequence for $F_\eps^p$ if 
		\begin{equation}
			\label{quasi minimal seq}
			F_\eps^{p}(u_\eps) - \inf_{u} F_\eps^{p}(u) \xrightarrow[\eps \to 0]{} 0.
		\end{equation}
	\end{definition}


	\subsection{Technical tools}

	\subsubsection{Separation}\label{Sect:separation}
	
	In this section, we discuss the property of separation satisfied by the surfaces $\surf_{\ell}$, where we consider more precisely the set $\curve(\gamma,\gamma_0)$ of Lipschitz curves in $\overline \Ccal_0$ connecting $\gamma$ to $\gamma_{0}$,
	\[
	\curve(\gamma,\gamma_0) : = \{ \ell \in \mathrm{Lip}([0,1] \times \mathbb{S}^{1}, \overline{\mathcal{C}_0}) \text{ such that } \ell(0) = \gamma \text{ and } \ell(1) = \gamma_{0} \}.
	\]
	This topological feature will be crucial for the proof of the $\Gamma-\liminf$ inequality.

	\begin{prop}
		\label{Thm:separation}
		\label{E_phi_separe}
		Let $\ell \in \curve(\gamma,\gamma_0)$. The set $\surf_{\ell}$ defined by~\eqref{Def:surface_image} separates the cylinder $\overline{\mathcal{C}_0}$, in the sense that $\overline{\mathcal{C}_0} \setminus \surf_{\ell}$ has at least two connected components, one containing the north pole $N=(0,0,h)$, the other one containing the south pole $S=(0,0,-h)$.
	\end{prop}
	
	In order to show this result we will use the following Theorem (see for instance~\cite[Chapter XVII]{dugundji1966topology}). 
	
	\begin{lemme}[Borsuk theorem] \label{Borsuk}
		Let $A \subset \R^{n}$ be compact, and let $p,q$ be in distinct connected components
		of $\R^{n}\setminus A$. If A is deformed over $\R^{n}$ into a set $B$, and if $A$ never crosses either $p$ or $q$ during this deformation, then $p,q$ are still in distinct connected components of $\R^{n}\setminus B$.
	\end{lemme}
	
	\begin{rmq}
		The hypothesis of deformation, in the above Lemma, can be written as the existence of a continuous function $\phi : A \times [0,1] \to \R^{n}$ such that, $\forall x \in A$, $\phi(x,0) = x$, $\phi(A\times \{1\}) \subset B$ and $p,q \notin \phi(A \times [0,1]).$
	\end{rmq}
	
	\begin{proof}[Proof of Proposition~\ref{Thm:separation}] Upon considering a higher cylinder, we can assume that there exists a cylinder $A_1$, with the same radius as $\mathcal{C}_0$ and smaller height, containing $S_{\ell}$ in its adherence. For simplicity, we also suppose that the height of the cylinder $A_{1}$ is 2 and that it is centered, \emph{i.e}, $z \in [-1,1]$. This way $\overline{A_{1}}$ contains the unitary closed disk with height 0, denoted $\mathbb{D} := \{(x,y,0) \in \R^{3}, \sqrt{x^{2}+y^{2}} \leqslant 1 \}$.
		
		
		
		We consider $A := \partial \mathcal{C}_0 \cup \mathbb{D}$ (which is compact) and $B := \partial \mathcal{C}_0 \cup \surf_{\ell}$.
		
		
		
		It is clear that the disk separates the cylinder, thus we can choose $p,q \in \mathcal{C}_0\setminus \overline{A_{1}}$ belonging to distinct connected components of $\R^{n}\setminus A$.
		We now check the assumptions of the Borsuk theorem. We will justify that we can deform $A$ into $B$. To this aim, we introduce the following function on $A \subset \R^{2}\times \R$:
		\[\phi_{1}(r,\omega,z) := \begin{cases}
			\ell(1-r,e^{i\omega}) &\text{ if } z=0\\
			(r,\omega, (1-|z|)\tilde{\gamma}(\omega) + z|z|) &\text{ if } |z| \leqslant 1\\
			(r,\omega,z) &\text{ if } |z|>1
		\end{cases}.\]
		This function is continuous on $A$ since for $z=0$ and $r=1$  we have $\ell(0,\omega) = \gamma(e^{i\omega}) = (1,\omega, \tilde{\gamma}(\omega))$ from the assumption that $\Gamma$ is a Lipschitz graph over the unitary circle. 
		
		Then, we define $\phi(\cdot,t) := (1-t)id + t \phi_{1}$ on $A$.
		It is clear that $\phi$ is continuous on $A \times [0,1]$ and that for all $x\in A$, we have $\phi(x,0) = x$  and $\phi(x,1) = \phi_{1}(x) \in B$. Finally, we must make sure that the deformation does not cross either $p$ or $q$. One can easily check that $\phi(\overline{A_{1}}\times [0,1]) \subset \overline{A_{1}}$ and $\phi(\partial \mathcal{C}_0 \setminus \overline{A_{1}} \times [0,1])  \subset \partial \mathcal{C}_0$, which guarantees that $p,q \notin \phi(A \times [0,1])$.

		Thus, we can apply Borsuk's Theorem and conclude that $p$ and $q$ are points in the cylinder $\Ccal_0$ and in distinct connected components of $B^{c} = (\partial \mathcal{C}_0 \cup S_\ell)^{c}$. Hence $\surf_{\ell}$ separates $\overline{\mathcal{C}_0}$.
	\end{proof}

	\begin{cor}\label{Coro:separation}
		Let $\ell \in \curve(\gamma,\gamma_0)$. The set $\surf_{\ell}\cup \Sigma$ separates the cylinder $\mathcal{C}$. 
	\end{cor}
	\begin{proof}
		This is an immediate consequence of the fact that the Lipschitz graph $\Sigma$ separates $\Ccal \setminus \overline{\Ccal_0}$ and Proposition~\ref{Thm:separation}.	
	\end{proof}


	\subsubsection{Covering a set by an open set with controlled perimeter}
	
	The following elementary lemma will be needed later for the proof of the liminf inequality.
	
	\begin{lemme}\label{OpenCovering}
		Let $K\subset \R^N$ be any set such that $0< \mathcal{H}^{N-1}(K) <+\infty$. Then for every $\varepsilon>0$ there exists an open set $A_\eps \subset \R^N$ such that:
		\begin{enumerate}
			\item $K\subset A_\eps$
			\item $A_\eps \subset \{x \in \R^N \; : \; dist(x,K)\leq \varepsilon\}$
			\item $P(A_\eps)\leq C \mathcal{H}^{N-1}(K)$,
		\end{enumerate}
		where $C>0$ is a universal constant.
	\end{lemme}
	
	\begin{proof} Let $K \subset \R^N$ and $\varepsilon>0$ be given. By definition of the Hausdorff measure, there exists a countable family of closed sets $A_i$ such that ${\rm diam}(A_i)\leq \varepsilon$ and  
		$$\sum_{i\in I} {\rm diam}(A_i)^{N-1}\leq C \mathcal{H}^{N-1}(K).$$
		For each $i\in I$ we let $B_i$ be an open ball of radius ${\rm diam}(A_i)$ such that $A_i\subset B_i$. Then the balls $\{B_i\}$ form a new covering of $K$.  We define $A_\eps:= \bigcup_{i\in I} B_i$. Then,
		$$P(A_\eps)\leq \sum_{i\in I} P(B_i) =  \sum_{i\in I} C_N {\rm diam}(B_i)^{N-1} \leq C  \mathcal{H}^{N-1}(K),$$
		thus $A_\eps$ fulfils all the requirements of the statement  and the lemma follows.
	\end{proof}



	\subsubsection{An average formula for finite measures}
	
	In this section we give  a standard average formula for finite measures that will be used later in the proof of the liminf inequality.
	
	\begin{lemme}
		Let $A\subset \R^n$ be a non empty set. Let $\mu$ be a finite measure on $A$ and $f:A\rightarrow \R$ be a $\mu$-measurable function such that $\int_A f\, d\mu \geqslant 0$.
		Then, there exists $t_{0} \in A$ such that 
		\begin{equation}
			\label{average formula}
			\int_{A} f d\mu \geqslant \mu(A) f(t_{0}).
		\end{equation}
	\end{lemme}

	\begin{proof}
		
		If $\mu(A)=0$, the statement~\eqref{average formula} is trivial, hence we may assume that $\mu(A)>0$ and set $m= \frac{\int_{A}f d\mu}{\mu(A)}$. If $m=0$, then either $f = 0$ $\mu$-a.e. on $A$, or there exists $t_0\in A$ such that $f(t_0)\leqslant 0$. In both cases, there exists $t_0\in A$ such that~\eqref{average formula} holds. Otherwise, if $m>0$, we recall the following inequality, valid for any $t>0$:
		\[\mu(\{f>t\}) = \int_{\{f>t\}}d\mu < \frac{1}{t}\int_{\{f>t\}}f d\mu \leqslant \frac{1}{t}\int_{A}f d\mu.\]
		Taking $t = m$, we get
		\[\mu(\{f>m\}) < \frac{1}{m}\int_{A}f d\mu = \mu(A).\]
		Thus, $\mu(\{f\leqslant m\})>0$, so
		there must exist $t_{0} \in A$ such that $f(t_{0}) \leqslant m$. This proves~\eqref{average formula}.
	\end{proof}


	\subsection{Proof of Theorem~\ref{Thm:limsup}: the limsup inequality}
	
	This section is devoted to the proof of the limsup inequality.

	\begin{proof}[Proof of Theorem \ref{Thm:limsup}]
		
		Let $\Omega \subset \Ccal$ be a competitor and  $\Omega_{0}$ be a minimizer for problem \eqref{def:pblimit}. It is enough to show that there exists a sequence $(u_{\eps}) \in H^{1}(\Ccal)\cap C(\overline{\Ccal})$ such that $0\leqslant u_{\eps} \leqslant 1$, $u_{\eps} = 1 $ on $\partial \Ccal$, for which the $\limsup$ inequality holds with $\Omega_{0}$, \emph{i.e.}
		\begin{equation}\label{limsupOmega0}
			\limsup_{\eps \to 0} F^{p}_{\eps}(u_{\eps}) \leqslant P(\Omega_0,\overline{\mathcal{C}_0}).
		\end{equation}
		Indeed, $P(\Omega_0,\overline{\mathcal{C}_0})\leqslant P(\Omega,\overline{\mathcal{C}_0})$ by definition of~\eqref{def:pblimit} (notice that in $\Ccal \setminus \overline{\Ccal_0}$ the essential boundaries both coincide with $\Sigma$), therefore inequality~\eqref{limsupOmega0} implies~\eqref{Ineq:limsup}.

		Let $\Omega_0$ be a solution of problem \eqref{def:pblimit}. To prove~\eqref{limsupOmega0}, we introduce the set $K := \partial^{*}\Omega_{0}\cap \overline{\mathcal{C}_0}$, where we recall that $\partial^{*}\Omega_{0}$ is the essential boundary of $\Omega_{0}$. Thus, $P(\Omega_{0},\mathcal{C}) = \Haus(\partial^{*} \Omega_{0} \cap \mathcal{C}) = \Haus(\partial^{*} \Omega_{0} \cap \overline{\mathcal{C}_0}) +\Haus(\partial^{*} \Omega_{0} \cap (\mathcal{C} \setminus \overline{\mathcal{C}_0})) = \Haus(K) + \Haus(\Sigma)$, since $\partial^{*} \Omega_{0} \cap (\mathcal{C} \setminus \overline{\mathcal{C}_0}) = \Sigma$, by definition of the boundary constraint.

		Then, we follow the standard construction for the $\Gamma-\limsup$ inequality, which is based on the well-known optimal profile of minimizers for a Modica-Mortola type energy. We recall from Proposition~\ref{th : 2.1} that, since $\Omega_{0}$ is a minimizer, its essential boundary is Ahlfors regular. Moreover, we also know from Proposition~\ref{th : 2.1} that $K\cap \partial \Ccal_0=\emptyset$. As a result, for all $x \in K$, and all $r >0$ such that $B(x,r) \subset {\Ccal}$,
		\begin{equation}
			\label{th 4.1 - AR}
			\Haus(K\cap B(x,r)) = P(\Omega_{0}, B(x,r))\geqslant C r^2.
		\end{equation}
		Applying~\cite[Theorem 2.104]{ambrosio2000oxford}, we deduce that the Minkowski content and the Hausdorff measure of $K$ coincide:
		\begin{equation}
			\label{Minkowski content}
			\lim_{r\to 0} \frac{\mathcal{L}^{3}(K_{r})}{2r} = \Haus(K),
			\quad\text{where } K_{r} := \{x \in \R^{3}\, |\, d(x,K) \leqslant r \}.
		\end{equation}
		

		Let $(k_{\eps})_{\eps>0}$ be a sequence of positive numbers, converging to $0$, that will be specified later on. 
		We define $a_{\eps} := -2\eps\ln{(\eps)}$, $b_{\eps} := \eps^{2}$ and $\lambda_{\eps} := \frac{1-k_{\eps}}{1-\eps}$ and we consider the function 
		\[u_{\eps} := \begin{cases}
			k_{\eps} & \text{ on } K_{b_{\eps}},\\
			k_{\eps} + \lambda_{\eps}\left(1- \exp{(\frac{b_{\eps}-d(x,K)}{2\eps})}\right) & \text{ on } K_{b_{\eps}+a_{\eps}} \setminus K_{b_{\eps}},\\
			1 & \text{ on } \overline{\Ccal}\setminus K_{a_{\eps}+b_{\eps}}.
			
		\end{cases}\]
		We also take $\eps$ small enough so that $K_{a_{\eps}+b_{\eps}} \subset \mathcal{C}$. We notice that with this definition, $u_{\eps}$ is continuous and equal to 1 on $\partial \Ccal$. We can also easily show that $u_{\eps}$ is Lipschitz on $\Ccal$, therefore $u_{\eps} \in H^{1}(\Ccal)\cap C(\overline{\Ccal})$.
		
		We claim that the sequence of functions $(u_\eps)_{\eps>0}$ satisfy the inequality
		\begin{equation}\label{Claim:prooflimsup}
			\limsup_{\eps \to 0} \left ( \int_{\mathcal{C}}\eps|\nabla u_{\eps}|^{2} + \frac{(1-u_{\eps})^{2}}{4\eps} \right) \leqslant \Haus(K).
		\end{equation}
		On $\mathcal{C} \setminus K_{b_{\eps}+a_{\eps}}$ it is clear that the integral is exactly $0$, and the contribution of the integral on $K_{b_{\eps}}$ is of order $\eps$ since  $\mathcal{L}^{3}(K_{b_{\eps}}) \leqslant C \eps^{2}$. The main contribution of the integral is thus attained on $K_{b_{\eps}+a_{\eps}} \setminus K_{b_{\eps}}$. Setting $\tau(x) := d(x,K)$, we may observe that, in this region, $u_{\eps}$ is a function of $\tau$, since $u_{\eps}(x) = \lambda_{\eps}f_{\eps}(\tau(x))+k_{\eps}$ with $f_{\eps}(t) := (1-\exp(\frac{b_{\eps}-t}{2\eps}))$. Moreover, this last function is the solution of the ordinary differential equation 
		\[\begin{cases}
			& f_{\eps}' = \frac{1-f_{\eps}}{2\eps},\\
			& f_{\eps}(b_{\eps}) = 0.
		\end{cases}\]
		Applying the coarea formula to the Lipschitz function $\tau$ (see for instance~\cite[Theorem 3.11]{EvansGariepy}), and noticing that 
		\[1 - k_{\eps} -\lambda_{\eps}f_{\eps}(t) =  \lambda_\eps(1-\eps - f_{\eps}(t)),\]
		we get 
		\begin{align*}
			A_{\eps} &:= \int_{K_{a_{\eps}+b_{\eps}}\setminus K_{b_{\eps}}} \eps |\nabla u_{\eps}|^{2} + \frac{(1-u_{\eps})^{2}}{4\eps}  \\
			&= \int_{b_{\eps}}^{b_{\eps}+a_{\eps}} \left(\int_{\tau(y) = t} \eps |\nabla u_{\eps}(y)|^{2} + \frac{(1-u_{\eps}(y))^{2}}{4\eps} d\Haus(y)\right) dt\\
			&= \int_{b_{\eps}}^{b_{\eps}+a_{\eps}}\left (\eps \lambda_{\eps}^{2}|f_{\eps}'(t)|^{2} + \frac{(1-k_{\eps}-\lambda_{\eps}f_{\eps}(t))^{2}}{4\eps} \right )\Haus(\{\tau = t\}) dt\\
			&= \lambda_{\eps}^{2}\int_{b_{\eps}}^{b_{\eps}+a_{\eps}} \left (\frac{(1-f_{\eps}(t))^{2}}{4 \eps} + \frac{(1-\eps-f_{\eps}(t))^{2}}{4\eps}\right )\Haus(\{\tau = t\}) dt\\
			&\leqslant  \frac{\lambda_{\eps}^{2}}{2\eps} \int_{b_{\eps}}^{b_{\eps}+a_{\eps}} (1-f_\eps(t))^2\Haus(\{\tau = t\}) dt\\
			&=  \frac{\lambda_{\eps}^{2}}{2\eps} \int_{b_{\eps}}^{b_{\eps}+a_{\eps}} \exp\left(\frac{b_{\eps}-t}{\eps}\right)\Haus(\{\tau = t\}) dt.
		\end{align*}
		Once again, the coarea formula yields 
		\[g(t) := \mathcal{L}^{3}(K_{t}) = \int_{0}^{t}\Haus(\{\tau = t\}) dt,\]
		hence $g$ is differentiable a.e.\! and $g'(t) = \Haus(\{\tau = t\})$. Integrating by parts, we deduce the upper bound
		\[A_{\eps} \leqslant \frac{\lambda_{\eps}^{2}\eps}{2}g(a_{\eps}+b_{\eps}) - \frac{\lambda_{\eps}^{2}}{2\eps}g(b_{\eps}) + \frac{\lambda_{\eps}^{2}}{2\eps^{2}}\int_{b_{\eps}}^{a_{\eps}+b_{\eps}}\exp \left ( \frac{b_{\eps}-t}{\eps} \right) \mathcal{L}^{3}(K_{t})dt.\]
		Coming back to~\eqref{Minkowski content}, we know that for any $\eta > 0 $ and for $t$ small enough, $\mathcal{L}^{3}(K_{t}) \leqslant 2t(\Haus(K) + \eta)$. As a result, using that $\lim_{\eps\ra 0}\lambda_{\eps} = 1$ and  $\lim_{t\ra 0}g(t) = 0$, we can estimate $A_\eps$ as follows:
		\begin{align*}
			\limsup_{\eps\ra 0} A_{\eps} &\leqslant \limsup_{\eps\ra 0} \left(\frac{\lambda_{\eps}^{2}\eps}{2}g(a_{\eps}+b_{\eps}) + \frac{\lambda_{\eps}^{2}}{2\eps^{2}}\int_{b_{\eps}}^{a_{\eps}+b_{\eps}}\exp \left ( \frac{b_{\eps}-t}{\eps} \right) \mathcal{L}^{3}(K_{t})dt \right) \\
			&\leqslant \limsup_{\eps\ra 0} \left( \frac{1}{2\eps^{2}} \int_{b_{\eps}}^{a_{\eps}+b_{\eps}}\exp \left ( \frac{b_{\eps}-t}{\eps} \right) \mathcal{L}^{3}(K_{t})dt \right)  \\
			&\leqslant (\Haus(K) + \eta)\, \limsup_{\eps\ra 0} \left( \frac{1}{2\eps^{2}} \int_{b_{\eps}}^{a_{\eps}+b_{\eps}}\exp \left ( \frac{b_{\eps}-t}{\eps} \right) 2t dt \right)\\
			&= (\Haus(K) + \eta)\, \limsup_{\eps\ra 0} \left (\int_{0}^{a_{\eps}/\eps}\exp( -s)(\frac{b_{\eps}}{\eps}+s)  ds\right),
		\end{align*}
		where we have used the change of variables $t=b_\eps+\eps s$ in the last integral.
		Since $\frac{a_{\eps}}{\eps} = -2\ln (\eps) \to + \infty$ and $\frac{b_{\eps}}{\eps} = \eps \to 0$, we get 
		\[\limsup_{\eps\ra 0} A_{\eps} \leqslant (\Haus(K) + \eta) \int_{0}^{+\infty} \exp(-s)s ds = \Haus(K) + \eta,\]
		and letting $\eta \to 0$, 
		\[\limsup_{\eps\ra 0} A_{\eps} \leqslant \Haus(K).\]
		The claim~\eqref{Claim:prooflimsup} follows from the previous observations.

		Now, thanks to Proposition~\ref{th : 2.1} there exists $\ell^*\in \curve(\gamma,\gamma_0)$ such that $K = S_{\ell^*}$. By Definition~\ref{Def:geoDistance}, if $p=\infty$,
		\begin{align*}
			d^{\infty}_{u_{\eps}}(\gamma, \gamma_{0})
			&\leqslant \sup_{S_{\ell^*}}(|u_{\eps}| +\delta_\eps)
			= \sup_{K} (|u_{\eps}|)+\delta_\eps
			=   k_{\eps}+\delta_\eps,
		\end{align*}
		and if $p < + \infty$,
		\begin{align*}
			d^{p}_{u_{\eps}}(\gamma, \gamma_{0})
			&\leqslant  \int_{S_{\ell^*}} (|u_{\eps}|^{p} +\delta_\eps )d\Haus 
			=  \int_{K} |u_{\eps}|^{p} d\Haus + \delta_\eps \Haus(K)
			=   \Haus(K) (k_{\eps}+ \delta_\eps).
		\end{align*}
		Since, $\frac{\delta_\eps}{c_\eps} \to 0$, setting $k_{\eps} = c_{\eps}^{2}$ ensures that, for every $p\in [1,+\infty]$,
		\begin{equation}\label{ProofLimSupIneq2}
			\lim_{\eps\ra 0} \frac{1}{c_{\eps}}d^{p}_{\varphi_{\eps}}(\gamma, \gamma_{0}) = 0.
		\end{equation}
		This choice of $k_{\eps}$ also guarantees that for $\eps$ small enough, $u_{\eps}$ takes values in $[0,1]$. Combining~\eqref{Claim:prooflimsup} and~\eqref{ProofLimSupIneq2}, we obtain inequality~\eqref{limsupOmega0} and conclude the proof.
	\end{proof}

	\subsection{Proof of Theorem~\ref{Thm:liminf}: the liminf inequality }
	
	This section is devoted to the proof of the liminf inequality.

	\subsubsection{Definition of the relaxed energies $E^{p}_\eps$}

	The proof will also rely on relaxed energies associated with the functional $F_\eps^p$, depending not only on a phase field $u$, but also on a path $\ell\in \curve(\gamma,\gamma_0)$.
	
	\begin{definition}
		Let $p\in [1,+\infty]$, $u \in H^{1}(\Ccal)\cap C(\overline{\Ccal})$ such that $0\leqslant u \leqslant 1$ and $u = 1$ on $\partial{\Ccal}$ and $\ell \in \curve(\gamma,\gamma_0)$. We define 
		\begin{align*}
			&E^{p}_{\varepsilon}(u,\ell) := \varepsilon\int_{\mathcal{C}} |\nabla u|^{2} dx + \frac{1}{4\varepsilon}\int_{\mathcal{C}}(1-u)^{2} dx + \frac{1}{c_{\varepsilon}}\int_{\surf_{\ell}} (|u|^{p} + \delta_\eps) d\Haus  \quad \text{if } p<+\infty,\\
			&E^{\infty}_{\varepsilon}(u,\ell) := \varepsilon\int_{\mathcal{C}} |\nabla u|^{2} dx + \frac{1}{4\varepsilon}\int_{\mathcal{C}}(1-u)^{2} dx + \frac{1}{c_{\varepsilon}}\sup_{\surf_{\ell}} (|u| + \delta_\eps).
		\end{align*}
	\end{definition}

	\begin{rmq}
		From Definition~\ref{Def:geoDistance} of the geodesic distance between closed curves, we know that for all $u\in H^{1}(\Ccal)\cap C(\overline{\Ccal})$ such that $0\leqslant u \leqslant 1$ and $u = 1$ on $\partial{\Ccal}$,  
		\[F^{p}_{\varepsilon}(u) = \inf\left \{ E^{p}_{\varepsilon}(u,\ell)\ ,\ell \in \curve(\gamma,\gamma_0)  \right \}.\]
	\end{rmq}

	\subsubsection{Case $p=\infty$.}\label{Proof:liminfInfini}

	In this section we prove Theorem~\ref{Thm:liminf} in the particular case $p=\infty$.
	
	\begin{proof}[Proof of Theorem~\ref{Thm:liminf} for $p=\infty$] Let $\Omega$ be a solution to Plateau's problem~\eqref{def:pblimit}, and 
		$u_{\eps} \in H^{1}(\Ccal)\cap C(\overline{\Ccal})$ such that $0\leqslant u_{\eps} \leqslant 1$ and $u_{\eps} = 1$ on $\partial {\Ccal}$.  
		\medskip
		
		\emph{Step 1.} We start by using the relaxed energy $E_\eps^{\infty}$ to avoid taking an infimum over $\curve(\gamma,\gamma_0)$ but instead approximating the geodesic distance $d^{\infty}_{u_\eps}(\gamma,\gamma_0)$ by $\sup_{\surf_{\ell_{\varepsilon}}} (|u_{\varepsilon}|+ \delta_\eps)$, for a well-chosen $\ell_\eps$. For instance, by Definition~\ref{Def:geoDistance}, we can find $\ell_{\varepsilon} \in \curve(\gamma,\gamma_0)$ such that \[\sup_{\surf_{\ell_{\varepsilon}}} (|u_{\varepsilon}|) \leqslant \sup_{\surf_{\ell_{\varepsilon}}} (|u_{\varepsilon}|+\delta_\eps) \leqslant d^{\infty}_{u_{\varepsilon}}(\gamma, \gamma_{0}) + \varepsilon c_{\varepsilon},\]
		which implies  $F^{\infty}_{\varepsilon}(u_{\varepsilon}) + \varepsilon \geqslant E^{\infty}_{\varepsilon}(u_{\varepsilon},\ell_{\varepsilon})$.
		Since 
		\[\liminf_{\eps\to 0} F^{\infty}_{\varepsilon}(u_{\varepsilon}) = \liminf_{\eps\to 0} (F^{\infty}_{\varepsilon}(u_{\varepsilon}) + \varepsilon) \geqslant \liminf_{\eps\to 0} E^{\infty}_{\varepsilon}(u_{\varepsilon},\ell_{\varepsilon}),\]
		it is enough to show the following inequality:
		\begin{equation}\label{Ineq:liminfRelax}
			\liminf_{\eps\to 0} E^{\infty}_{\varepsilon}(u_{\varepsilon},\ell_{\varepsilon}) \geqslant P(\Omega,\overline{\mathcal{C}_0}),
		\end{equation}
		to deduce~\eqref{Ineq:liminf}.

		\medskip
		
		\emph{Step 2.} We identify a sublevel-set of $u_\eps$ that contains the surface $\surf_{\ell_{\varepsilon}}$, which possesses the crucial property of separation (see Proposition~\ref{Thm:separation}).
		Without loss of generality, we may assume that $\liminf_{\eps\ra 0} E^{\infty}_{\varepsilon}(u_{\varepsilon},\ell_{\varepsilon}) $ is finite, and achieved by a subsequence, still denoted $\varepsilon$. Consequently, the energy $E^{\infty}_{\varepsilon}(u_{\varepsilon},\ell_{\varepsilon})$ is uniformly bounded, so there exists a constant $C>0$ such that 
		\begin{equation*}
			\sup_{\surf_{\ell_{\varepsilon}}} (|u_{\varepsilon}|) < C c_{\varepsilon}.
		\end{equation*}
		In particular,
		\begin{equation}\label{InclusionSurfLevelSetu}
			\surf_{\ell_{\varepsilon}} \subset \{u_{\varepsilon} < Cc_{\varepsilon}\}.
		\end{equation}

		\medskip

		\emph{Step 3.} Using the above inclusion~\eqref{InclusionSurfLevelSetu} and a classical estimate of the Ambrosio-Tortorelli term in $E_\eps^{\infty}$, we now exhibit a sublevel set $\{g_\eps\leq t_\eps\}$ of the function $g_\eps:=u_\eps-\frac{u_\eps^2}{2}$, whose perimeter in $\Ccal$ allows us to estimate from below $\liminf_{\eps\ra 0} E^{\infty}_{\varepsilon}(u_{\varepsilon},\ell_{\varepsilon})$ by
		\begin{equation}\label{Proofliminf:goalStep3}
			\liminf_{\eps\to 0} E^{\infty}_{\varepsilon}(u_{\varepsilon},\ell_{\varepsilon}) \geqslant \frac{1}{2} \liminf_{\eps\to 0} P(\{g_{\varepsilon} \leqslant t_{\varepsilon} \}, \mathcal{C}).
		\end{equation}
		
		Inequality~\eqref{Proofliminf:goalStep3} can be proved as follows.
		Setting $s_{\varepsilon} = Cc_{\varepsilon} - \frac{C^{2}c_{\varepsilon}^{2}}{2}$, the inclusion~\eqref{InclusionSurfLevelSetu} yields by monotonicity of the mapping $t\mapsto t-\frac{t^2}{2}$ on $[0,1]$, the inclusion
		\begin{equation}\label{InclusionSurfLevelSetg}
			\surf_{\ell_{\varepsilon}} \subset  \{g_{\varepsilon} < s_\eps\}.
		\end{equation}
		Notice that $g_\eps$ takes values in $[0,1/2]$ a.e.\! in $\Ccal$, and that $s_\eps$ is positive for $\eps$ small enough. The definition of $g_\eps$ comes from the classical estimate
		\begin{align*}
			E^{\infty}_{\varepsilon}(u_{\varepsilon},\ell_{\varepsilon})
			&\geqslant \varepsilon \int_{\mathcal{C}}|\nabla u_{\varepsilon}|^{2}dx + \frac{1}{4\varepsilon}\int_{\mathcal{C}}(1-u_{\varepsilon})^{2}dx\\
			&\geqslant \int_{\mathcal{C}} |\nabla u_{\varepsilon}||1-u_{\varepsilon}|dx\\
			&= \int_{\mathcal{C}} |\nabla g_{\varepsilon}|dx.
		\end{align*}
		Since $g_\eps$ is an antiderivative of $\nabla u_{\varepsilon}(1-u_{\varepsilon})$, $u_{\varepsilon} \in H^{1}(\Ccal)$ and $\Ccal$ is bounded,  $g_{\varepsilon}$ is in  $W^{1,1}(\Ccal)$. In particular, it has bounded variations in $\Ccal$.
		As a result, we can apply the coarea formula for $BV$ fonctions, to obtain 
		\[\int_{\mathcal{C}} |\nabla g_{\varepsilon}|dx = \int_{\R} P(\{g_{\varepsilon} > t \}, \mathcal{C}) dt = \int_{0}^{1/2} P(\{g_{\varepsilon} > t \}, \mathcal{C}) dt \geqslant \int_{s_{\varepsilon}}^{1/2 - \varepsilon} P(\{g_{\varepsilon} > t \}, \mathcal{C}) dt . \]

		Then, we apply the average formula \eqref{average formula} with $\mu = \mathcal{L}$, $A = [s_{\varepsilon},1/2-\varepsilon]$ and $f(t) =P(\{g_{\varepsilon} > t \}, \mathcal{C})$. This gives us the existence of $t_{\varepsilon} \in [s_{\varepsilon},1/2-\varepsilon]$ that satisfies  
		\[ \int_{s_{\varepsilon}}^{1/2} P(\{g_{\varepsilon} > t \}, \mathcal{C}) dt \geqslant (\frac{1}{2}-\varepsilon -s_{\varepsilon})P(\{g_{\varepsilon} > t_{\varepsilon}\},\mathcal{C}).\]
		Hence, we get 
		\[\liminf_{\eps\to 0} E^{\infty}_{\varepsilon}(u_{\varepsilon},\ell_{\varepsilon}) \geqslant \liminf_{\eps\to 0} (\frac{1}{2}-\varepsilon-s_{\varepsilon})P(\{g_{\varepsilon} > t_{\varepsilon} \}, \mathcal{C}) = \frac{1}{2} \liminf_{\eps\to 0} P(\{g_{\varepsilon} \leqslant t_{\varepsilon} \}, \mathcal{C})  .\]
		This proves~\eqref{Proofliminf:goalStep3}.

		\medskip
		
		\emph{Step 4.} To conclude the proof of inequality~\eqref{Ineq:liminfRelax}, we will construct two competitors $\Omega_\eps^1$ and $\Omega_\eps^2$ such that the following inequality holds:
		\begin{equation}
			\label{eq competitors}
			P(\{g_\eps \leq t_\eps \},\mathcal{C}) \geqslant P(\Omega_\eps^1,\mathcal{C}) + P(\Omega_\eps^2,\mathcal{C}) - 2\Haus(\Sigma).
		\end{equation}
		
		The construction of these two competitors constitutes the main novelty of our approach. Its main ingredient is the topological property of separation satisfied by the surface $\surf_{\ell_\eps}$ (see Theorem~\ref{Thm:separation}).

		We denote 
		\begin{equation}
			\label{A_eps}
			A_{\varepsilon} := \{g_{\varepsilon} \leqslant t_{\varepsilon}\}
		\end{equation}
		and set $\Tilde{A_{\varepsilon}} := A_{\varepsilon}\cap \overline{\mathcal{C}_0}$. Since intersecting a set of finite perimeter with a convex set reduces its perimeter (see for instance~\cite[exercise 15.14]{maggi_sets_2012}, and~\cite{Machefert2025} for a proof), 
		\[P(\Tilde{A}_{\varepsilon},\mathcal{C}) = P(A_{\varepsilon}\cap \overline{\mathcal{C}_0},\mathcal{C})
		\leqslant P(A_{\varepsilon},\mathcal{C}).\]
		Then we recall that $S_{\ell_\eps} \subset \{g_{\varepsilon} < s_{\varepsilon} \}\cap \overline{\mathcal{C}_0} \subset A_{\varepsilon}\cap \overline{\mathcal{C}_0} = \Tilde{A}_{\varepsilon}$ because $s_{\varepsilon} \leqslant t_{\varepsilon}$. By Corollary~\ref{Coro:separation}, $\Tilde{A}_{\varepsilon}\cup \Sigma$ separates $\mathcal{C}$. This implies that $\mathcal{C} \backslash (\Tilde{A}_{\varepsilon}\cup \Sigma)$ contains at least two connected components. We denote by $\Omega_{\varepsilon}^{1}$ the connected component containing the north pole $N = (0,0,h)$ (where $2h$ is the height of the centered cylinder $\mathcal{C}$) and $\Omega_{\varepsilon}^{2}$ the one containing the south pole $S = (0,0,-h)$, as represented on Figure~\ref{fig:liminf}. By construction, $\Omega_\varepsilon^1$ and $(\Omega_\varepsilon^2)^c$ are two competitors.  Notice that, since $g_\varepsilon$ is continuous, then $\Omega_\varepsilon^1$ and $\Omega_\varepsilon^2$ are   open sets.

		Now to see that inequality \eqref{eq competitors} holds true, we will identify the essential boundary of $\Omega_\varepsilon^i$ in the three disjoint regions $\Ccal_0$, $\Ccal\setminus \overline{\Ccal_0}$, and $\partial \Ccal_0$.

		We start with the region $\Ccal_0$ and   notice that, 
		$$\partial \Omega_\varepsilon^1 \cap \partial \Omega_\varepsilon^2 \cap \Ccal_0=\emptyset.$$
		Indeed,  we know that $S_{\ell_\eps}$ separates $\overline{\Ccal_0}$ thus  $\overline{\Omega^1_\eps\cap \Ccal_0}$ is contained in the (different) connected component of $\overline{\Ccal_0}\setminus S_{\ell_\eps}$ containing $N$,  and $\overline{\Omega^2_\eps \cap \Ccal_0}$ is contained in the connected component of $\overline{\Ccal_0}\setminus S_{\ell_\eps}$ containing $S$ (notice that $S_{\ell_\varepsilon}$ cannot touch $\overline{\Omega^i_\eps\cap \Ccal_0}$ because  $g_\varepsilon$ is lower than $s_\eps$ on it which is strictly less than $t_\eps$).   This yields a contradiction in the case when there would exist $x\in \partial \Omega_\varepsilon^1 \cap \partial \Omega_\varepsilon^2 \cap \Ccal_0$, then it would belong to two different connected components of $\overline{\Ccal_0}\setminus S_{\ell_\eps}$. Since the essential boundary is always contained in the topological boundary, we deduce that 
		\begin{equation}
			\partial^* \Omega_\varepsilon^1 \cap \partial^* \Omega_\varepsilon^1 \cap \Ccal_0=\emptyset. \label{intersect}
		\end{equation}
		
		Now in $\Ccal \setminus \overline{\Ccal_0}$ we clearly have, by construction, that $\partial^*\Omega_\varepsilon^i$ coïncides with $\Sigma$. It remains to identify the essential boundary of $\Omega_\varepsilon^i$ on $\partial \Ccal_0$. But here again, we must have that 
		
		\begin{equation}
			\partial \Omega_\varepsilon^1 \cap \partial \Omega_\varepsilon^2 \cap \partial \Ccal_0 \setminus \Gamma=\emptyset. \label{amontrerencore}
		\end{equation}
		
		Indeed, the curve  $\Gamma$ separates $\partial \Ccal_0$ into an upper part $Z^+\subset \partial \Ccal_0\setminus \Gamma$ and a lower part $Z^- \subset \partial \Ccal_0\setminus \Gamma$. Let us show that  $\partial \Omega_\varepsilon^1 \cap \partial \Ccal_0 \setminus \Gamma \subset Z^+$ and $\partial \Omega_\varepsilon^2 \cap \partial \Ccal_0 \setminus \Gamma \subset Z^-$. If not, then there would be for instance a point $x \in Z^-$ and a sequence $x_n\in \Omega_\varepsilon^1$ such that $x_n\to x$. Since $\Omega_\varepsilon^1$ is disjoint from $\mathcal{D}^-$, this means that $x_n \in \overline{\Ccal_0}$ for all $n$. But this is not possible  since $S_{\ell_\varepsilon}$ separates $Z^+$ from $Z^-$ in $\overline{\Ccal_0}$, and that $\overline{\Omega_\varepsilon^1\cap \Ccal_0}$ belongs to the same connected component of $\overline{\Ccal_0}\setminus S_{\ell_\varepsilon}$ as $Z^+$. Since $x\in  \overline{\Omega_\varepsilon^1\cap \Ccal_0}$, it cannot be in $Z^-$. This proves that $\partial \Omega_\varepsilon^1 \cap \partial \Ccal_0 \setminus \Gamma \subset Z^+$.  The same holds for $\partial \Omega_\varepsilon^2$ and  \eqref{amontrerencore} follows.

		From \eqref{intersect}, \eqref{amontrerencore}, using also the fact that $\Gamma\cap \partial \Ccal_0$ has zero $\mathcal{H}^2$-measure and the fact that $\partial \Omega_\eps^i$ coincides with $\Sigma$ on $\Ccal \setminus \overline{\Ccal_0}$, we deduce that

		\begin{equation}
			\label{eq competitors2}
			P(\{g_\eps \leq t_\eps \},\mathcal{C}) \geqslant P(\Omega_\eps^1\cup \Omega_\eps^2,\mathcal{C})=  P(\Omega_\eps^1,\mathcal{C}) + P(\Omega_\eps^2,\mathcal{C}) - 2\Haus(\Sigma), \\
		\end{equation}
		and so follows \eqref{eq competitors}.
		
		\begin{figure}
			\centering
			\includegraphics[width = 0.4\textwidth]{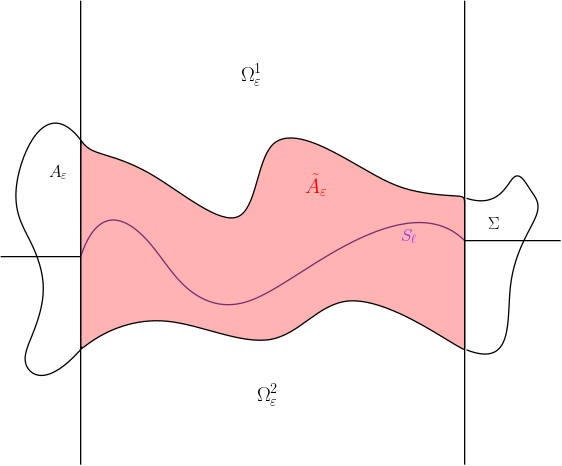}
			\caption{Construction of the competitors $\Omega_\varepsilon^1$ and $\Omega_\varepsilon^2$.}
			\label{fig:liminf}
		\end{figure}

		\medskip
		
		\emph{Step 5.}
		Since $\liminf_{\eps \to 0} E_\eps^\infty$ is finite, \eqref{eq competitors} yields that the perimeters of $\Omega_\varepsilon^1$ and $\Omega_\varepsilon^2$ are uniformly bounded. Thus,   there exist subsequences (not relabeled) of $(\Omega_\varepsilon^1)$ and $(\Omega_\varepsilon^2)$ and two sets of finite perimeter in $\Ccal$, $\Omega^1$ and $\Omega^2$, such that
		\[\lim_{\eps \to 0} \Omega_\varepsilon^1 = \Omega^1 \text{ and } \lim_{\eps \to 0} \Omega_\varepsilon^2 = \Omega^2 \text{ in } L^1(\Ccal). \]
		This naturally implies that $\Omega^1$ and $(\Omega^2)^c$ are competitors. And, the lower-semi-continuity of the perimeter, combined with \eqref{Proofliminf:goalStep3} and \eqref{eq competitors}, yields 
		\begin{align}
			\liminf_{\eps \to 0} E_\eps^\infty(u_\eps,\ell_\eps) &\geqslant \frac{1}{2} (P(\Omega^1, \Ccal) + P(\Omega^2, \Ccal)) - \Haus(\Sigma) \notag \\
			&=  \frac{1}{2} (P(\Omega^1, \overline{\Ccal_0}) + P(\Omega^2, \overline{\Ccal_0}))\\
			&=  \frac{1}{2} (P(\Omega^1, \overline{\Ccal_0}) + P((\Omega^2)^c, \overline{\Ccal_0})). \label{preuve liminf competiteurs limites}
		\end{align}

		By definition, $\Omega^1$ and $\Omega^2$ have larger perimeter in $\overline{\Ccal_0}$ than the minimizer $\Omega$, so inequality~\eqref{preuve liminf competiteurs limites}   clearly implies~\eqref{Ineq:liminfRelax}.
	\end{proof}

	\subsubsection{Case $p<\infty$.}

	\begin{proof}[Proof of Theorem~\ref{Thm:liminf} for $p<\infty$]
		\medskip
		
		\emph{Step 1.} As in the $p=\infty$ case, we introduce $\ell_{\varepsilon} \in \curve(\gamma,\gamma_0)$ such that
		\[\int_{\surf_{\ell_{\varepsilon}}} |u_{\varepsilon}|^p d\Haus \leqslant\int_{\surf_{\ell_{\varepsilon}}} (|u_{\varepsilon}|^p + \delta_\eps)d\Haus \leqslant d^{p}_{u_{\varepsilon}}(\gamma, \gamma_{0}) + \varepsilon c_{\varepsilon}.\]

		Hence, we need to prove the inequality
		\begin{equation}\label{Ineq:liminfRelaxpfinite}
			\liminf_{\eps\to 0} E^{p}_{\varepsilon}(u_{\varepsilon},\ell_{\varepsilon}) \geqslant P(\Omega,\overline{\mathcal{C}_0}).
		\end{equation}

		\medskip
		
		\emph{Step 2.} As in the case $p=\infty$, we may assume that $E^p_\eps(u_\eps,\surf_{\ell_\eps})\leqslant C$, in particular
		\[\int_{\surf_{\ell_{\varepsilon}}} |u_{\varepsilon}|^pd\Haus \leqslant C c_{\varepsilon}.\]
		Contrary to the $p=\infty$ case, we cannot conclude from the above inequality that $\surf_{\ell_\eps}$ is contained in the sublevel-set $\{u_\eps<Cc_\eps\}$. However, Tchebychev inequality yields that for all $\alpha > 0$, \[\Haus(\surf_{\ell_{\varepsilon}}\cap \{u_{\varepsilon}\geqslant \alpha\}) \leqslant \left(\frac{Cc_{\varepsilon}}{\alpha}\right)^{p}.\]
		Thus, taking $\alpha_{\varepsilon} := \sqrt{c_{\varepsilon}}>0$ leads to 
		\begin{equation}
			\label{limit K_eps}
			\Haus(\surf_{\ell_{\varepsilon}}\cap \{u_{\varepsilon}\geqslant \alpha_{\varepsilon}\}) \longrightarrow 0.
		\end{equation}
		Now we introduce 
		$
		K_{\varepsilon} := \surf_{\ell_{\varepsilon}}\cap \{u_{\varepsilon}> \alpha_{\varepsilon}\}$,
		and we  apply Lemma \ref{OpenCovering} to find an open set $U_\varepsilon \subset \Ccal$, such that $K_\varepsilon \subset U_\varepsilon$ and
		\begin{equation}
			P(U_\varepsilon)\leq C \mathcal{H}^2(K_\varepsilon),   \label{amontrerrr}
		\end{equation}
		for a universal constant $C>0$.  By construction we have 
		\begin{equation}\label{surfinclusionfinite}
			\surf_{\ell_{\varepsilon}} \subset \{u_{\varepsilon} \leq  \alpha_{\varepsilon}\} \cup U_{\eps}.
		\end{equation}

		The main idea of the proof relies on the fact that the set $\{u_{\varepsilon} \leq  \alpha_{\varepsilon}\} \cup U_\varepsilon$ still possesses the separation property whilst adding the set $U_\varepsilon$ to the sublevel set of $u_\eps$ does not contribute to the perimeter in the limit.

		\medskip

		\emph{Step 3.} As previously, we introduce the function $g_\eps:=u_\eps-\frac{u_\eps^2}{2}$ and using the monotonicity of $t\mapsto t-t^2/2$, the inclusion~\eqref{surfinclusionfinite}, and setting $s_\eps=\alpha_{\varepsilon}-\frac{\alpha_\eps^2}{2}$, we get
		\[\surf_{\ell_{\varepsilon}} \subset \{g_{\varepsilon} < s_\eps\} \cup U_\eps.\]
		The same computations as in case $p=\infty$ yield the existence of $t_\eps\in [s_\eps,1/2-\eps]$ such that
		\[E^p_{\eps}(u_\eps,\ell_\eps)\geqslant \int_{s_{\varepsilon}}^{1/2} P(\{g_{\varepsilon} > t \}, \mathcal{C}) dt \geqslant (\frac{1}{2}-\varepsilon -s_{\varepsilon})P(\{g_{\varepsilon} > t_{\varepsilon}\},\mathcal{C}).\]
		We claim moreover  that 
		\begin{equation}\label{Claim:liminfK}
			\liminf_{\eps\to 0} P(\{g_{\varepsilon} \leqslant t_{\varepsilon} \}, \mathcal{C}) = 
			\liminf_{\eps\to 0} P(\{g_{\varepsilon} \leqslant t_{\varepsilon} \}\cup U_\eps, \mathcal{C}).
		\end{equation}
		This follows from the following estimates,
		\[
		P(\{g_{\varepsilon} \leqslant t_{\varepsilon} \}\cup U_\eps, \mathcal{C}) \leqslant P(\{g_{\varepsilon} \leqslant t_{\varepsilon} \}, \mathcal{C}) + P(U_\eps, \mathcal{C}),
		\]
		and  
		\[
		\liminf_{\eps\to 0} P(U_\eps, \mathcal{C})=0
		\]
		because of \eqref{limit K_eps} and \eqref{amontrerrr}.

		Then, following exactly the same reasoning as in \emph{Step 4} from Section~\ref{Proof:liminfInfini}, replacing the sublevel-set $\{g_\eps\leqslant t_\eps\}$ by $\{g_\eps\leqslant t_\eps\}\cup U_\eps$, we obtain the desired $\liminf$.
	\end{proof}

	\subsection{Proof of Theorem~\ref{th : csq Gamma cv}}
	
	We finally give a proof for  the main result stated in the Introduction.

	\begin{proof}[Proof Theorem~\ref{th : csq Gamma cv}] We start with the case $p=+\infty$. We recall the conclusion \eqref{preuve liminf competiteurs limites} from the proof of Theorem~\ref{Thm:liminf}, where we have constructed two competitors $\Omega^1$ and $(\Omega^2)^c$ such that 
		
		\[\liminf_{\eps \to 0} F_\eps^p(u_\eps) \geqslant\liminf_{\eps \to 0} E_\eps^p(u_\eps,l_\eps) \geqslant \frac{1}{2} (P(\Omega^1, \overline{\Ccal_0}) + P((\Omega^2)^c, \overline{\Ccal_0})).\]
		
		We denote by $\tilde{u_\eps}$ the sequence obtained in Theorem~\ref{Thm:limsup}. Since $u_\eps$ is assumed to be a quasi-minimizing sequence for $F_\eps^p$, it follows that 
		
		\[F_\eps^p(u_\eps) \leqslant F_\eps^p(\tilde{u_\eps}) + o(1).\]
		
		Thus passing to the limit and applying Theorem~\ref{Thm:limsup} yield 
		
		\begin{align*}
			\frac{1}{2} (P(\Omega^1, \overline{\Ccal_0}) + P((\Omega^2)^c, \overline{\Ccal_0}))  & \leqslant \liminf_{\eps \to 0} F_\eps^p(u_\eps)
			\leqslant \liminf_{\eps \to 0} F_\eps^p(\tilde{u_\eps})\\
			&\leqslant \limsup_{\eps \to 0} F_\eps^p(\tilde{u_\eps})
			\leqslant P(\Omega_0, \overline{\Ccal_0}),
		\end{align*}
		
		where $\Omega_0$ is a minimizer of Plateau's problem \eqref{def:pblimit}. Hence, 
		\begin{equation}
			\label{eq limit perimeter}
			\frac{1}{2} (P(\Omega^1, \overline{\Ccal_0}) + P((\Omega^2)^c, \overline{\Ccal_0})) \leqslant P(\Omega_0, \overline{\Ccal_0}).
		\end{equation}

		This impose that at least $P(\Omega^1, \overline{\Ccal_0})$ or $P((\Omega^2)^c, \overline{\Ccal_0})$ is smaller than $P(\Omega_0, \overline{\Ccal_0})$. Let assume that $P(\Omega^1, \overline{\Ccal_0})$ is smaller that $ P(\Omega_0, \overline{\Ccal_0})$. Then, since $\Omega^1$ is a competitor and $\Omega_0$ is a minimizer, it follows that $\Omega^1$ is also a minimizer of Plateau's problem, \emph{i.e.}, $P(\Omega^1, \overline{\Ccal_0}) = P(\Omega_0, \overline{\Ccal_0})$. Thus, \eqref{eq limit perimeter} yields that $P((\Omega^2)^c,\overline{\Ccal_0}) \leqslant P(\Omega_0,\overline{\Ccal_0})$, which implies that $(\Omega^2)^c$ is also a minimizer. 
		
		Hence, both $\Omega^1$ and $\Omega^2$ are minimizers of Plateau's problem~\eqref{def:pblimit}, and this achieves the proof of Theorem~\ref{th : csq Gamma cv} in the case when $p=+\infty$.

		Now if $p<+\infty$, the only difference is that $\Omega_1$ is no more the limit of the connected component $\Omega_\varepsilon$ of the level set  $\{g_\eps > t_\eps\}$, but we had to remove the set $U_\varepsilon$ for topological reasons. But   since $P(U_\varepsilon)\to 0$, we deduce from the isoperimetric inequality (see~\cite[Proposition 12.35]{maggi_sets_2012}) that $|U_\varepsilon|\to 0$. This proves that the $L^1$-limit  of $\Omega_\varepsilon$   is the same as $\Omega_\varepsilon \setminus U_\varepsilon$, and achieves the proof of of Theorem~\ref{th : csq Gamma cv}.
	\end{proof}

	

	

	\section{Numerical experiments}\label{Section:numerics}
	
	This section details the numerical discretization schemes that enabled us to perform the experiments presented in the Introduction.
	We consider Plateau's problem where the given boundary $\Gamma$ is a finite union of $d$ closed curves:
	\[ 
	\bdy = \cup_{i=1}^d \bdy^{i}.
	\]
	We also assume that each component $\bdy^i$ is parametrized by a Lipschitz function $\bpar^i: \mathbb{S}^{1}\ra \R^3$.  The case analysed in Section~\ref{Section:Analysis} corresponds to $d=1$  where the geodesic distance term connects the curve $\bpar:=\bpar^1$ to a fixed point $x_2\in \bdy^1$, represented by a constant curve $\bpar^2$.  
	
	To extend this approach to cases with multiple components ($d\ge 2$), we can add a finite collection of distinct points $x^j \in \bdy^{j-d}$, indexed by $d+1\le j\le N=2d$ and  identified to a constant curve $\bpar^j$. In this general setting, the first step towards the generalization of the approximation functional $F_\eps$ defined in case $d=1$ in Section~\ref{Section:Analysis}, will be to set
	\begin{equation}\label{Def:NumericalFunctional}
		F_\eps(u) = 	\AT(u) + \frac{1}{c_\eps}\sum_{(i,j) \in I_{\gamma}}d^2_u(\bpar^i,\bpar^j)
	\end{equation}
	where $I_{\gamma}$ contains the couples $(i,j)$ for which
	we propose to connect the curves $\bpar^i$, $\bpar^j$.
	Here $d_u$ is still  the geodesic distance  introduced in Definition~\ref{Def:geoDistance} (with $p=2$) and defined by
	\[
	d^2_{u}(\bpar^i,\bpar^j) := \inf\left \{ \int_{\surf_{\ell}} (u^2 + \delta_{\eps} )d\Haus |\  \ell  \in \curve(\bpar^i,\bpar^j) \right \},\\
	\]
	with
	\[
	\curve(\bpar^i,\bpar^j) : = \{ \ell \in \mathrm{Lip}([0,1] \times \mathbb{S}^{1}, \overline{\mathcal{C}_0}) \text{ such that } \ell(0) = \bpar^i \text{ and } \ell(1) = \bpar^j \},
	\]
	and $\surf_{\ell} = \ell([0,1]\times \mathbb{S}^{1})$.

\begin{figure}
		\begin{minipage}{0.5\textwidth}
			\centering
			\includegraphics[scale=0.4]{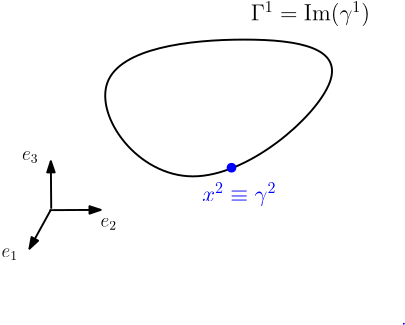}
		\end{minipage}
		\begin{minipage}{0.5\textwidth}
			\centering
			\includegraphics[scale=0.4]{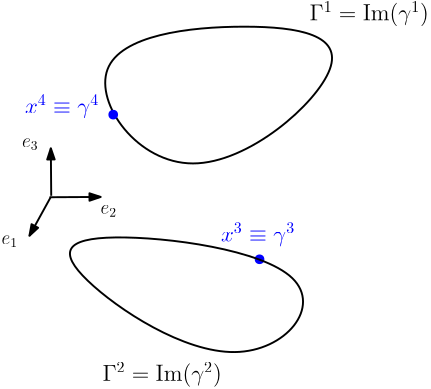}
		\end{minipage}
		\caption{Left: case of one curve $\Gamma^1$, parametrized by $\gamma^1$ and connected to a point $x^2\in \Gamma^1$, represented by the constant curve $\gamma^2$. Right: case of two curves $\Gamma^1$ and $\Gamma^2$, resp.\! parametrized by $\gamma^1$ and $\gamma^2$. Several possibilities occur: $\Gamma^1$ can be connected to a point $x^3\in \Gamma^2$, or $\Gamma^2$ to a point $x^4\in \Gamma^1$, or the curves $\Gamma^1$ and $\Gamma^2$ can be connected globally through the use of the geodesic distance $d_u(\bpar^1,\bpar^2)$.}
		\label{fig_plateau_type}
	\end{figure}

	As shown in Figure~\ref{fig_plateau_type}, let us give some examples to illustrate the possible choices of energies $F_\eps$ of the form~\eqref{Def:NumericalFunctional}.

	\paragraph{Case $d=1$:}  $\bdy^1$ is the image of one closed curve  $\bpar^1:\mathbb{S}^1\ra \R^3$. In this situation, we use the exact same framework as the one studied in Section~\ref{Section:Analysis}: for a given point $x^2\in \bdy^1$, identified with a constant curve $\bpar^2$, the energy $F_\eps$ reads
	\[
	F_\eps(u) =	\AT(u)  + \frac{1}{c_\eps}d^2_u(\bpar^1,\bpar^2).
	\]
	
	\paragraph{Case $d=2$:}  $\bdy$ is the union of two distinct curves $\bdy^1$ and $\bdy^2$, which are respectively parametrized by $\bpar^1$ and $\bpar^2$. In such configuration, one may consider several strategies in order to preserve the topology of the approximate optimal surface.
	\begin{enumerate}
		\item[-] Connect $\bdy^1$ to $\bdy^2$ globally, which conduces to
		\begin{equation}\label{EnergyFepsConfig1}
			F_\eps(u) =	\AT(u)  + \frac{1}{c_\eps}d^2_u(\bpar^1,\bpar^2).
		\end{equation}
		\item[-] Connect $\bdy^1$ to a point $x^3 \in \Gamma^2$, and $\bdy^2$ to a point $x^4 \in \Gamma^1$ as well. In this second configuration, the energy reads
		\begin{equation}\label{EnergyFepsConfig2}
			F_\eps(u) = 	\AT(u)  + \frac{1}{c_\eps}\Big[d^2_u(\bpar^1,\bpar^3)+d^2_u(\bpar^2,\bpar^4)\Big].
		\end{equation}
		
	\end{enumerate}

	Since Plateau's problem can be interpreted as an extension of Steiner's problem in three dimensions and the phase field models are ultimately very similar, we begin by recalling  in Section~\ref{Sect:SteinerNumerics} the numerical strategies developed in~\cite{BBL2020} to obtain efficient approximations to Steiner's problem. In particular, we review the following concepts:
	\begin{itemize}
		\item[-] relaxing the energy $F_\eps$ avoids the need to compute the derivative of the geodesic term with respect to $u$;
		\item[-] regularizing the geodesic term  to improve the smoothness of phase field function $u$;
		\item[-] using a variant of the Ambrosio-Tortorelli term to improve the regularity of the phase field solution $u$, and facilitate its minimization;
		\item[-] detailing the schemes used to solve the phase field PDE.
	\end{itemize}
	
	In section~\ref{Sect:PlateauNumerics}, we apply the same strategy to Plateau's problem. The novelty mainly lies in our treatment of the geodesic term, for which we propose two approaches. The first one, which is limited to circles, demonstrates how to reduce the problem of computing an optimal geodesic in infinite dimensions to a reduced space. For this smaller space, we can use fast marching algorithms in 2D. However, this approach is too restrictive, so we then explain how to compute a non-optimal geodesic easily allowing us to deal with a large number of cases.   Despite the non-optimality of the calculated geodesics, our model's ability to find approximations of the minimal surface in many configurations is illustrated by numerical experiments, extending far beyond the mathematical study of Plateau's problem of a curve in a cylinder.


	\subsection{Phase-field models and discretization of Steiner's problem}\label{Sect:SteinerNumerics}
	
	In this section, we focus on Steiner's problem in dimension two to observe the influence of the phase field model on numerical solutions. We first recall the numerical discretization strategy proposed in the article~\cite{BBL2020}, following the approximation strategy introduced in~\cite{bonnivard2015approximation} and refined in~\cite{bonnivard2018phase}. We then compare this method with a new phase field model for which we replace the Ambrosio-Tortorelli term $\AT$ with a higher-order version whose profile is smoother and better located around the diffuse interface.

	\medskip
	
	Consider a bounded and convex open set $\Omega \subset \R^2$. The Steiner problem consists in finding, for a given collection of points $a_0,\ldots, a_N \in \Omega$, a compact connected set $K \subset \Omega$ containing all the $a_i$’s and having minimal length. The idea is to obtain a minimizer of energy
	$$ F_{\varepsilon}(u) =  \AT(u) + \frac{1}{\lambda_{\varepsilon}} \sum_{i=1}^{N} {\bf D}(u^2 + \delta_{\varepsilon}; a_0, a_i), $$
	where ${\bf D}(w; a, b)$ is now defined by
	$$   {\bf D}(w; a, b) := \inf_{\Gamma \in G_{a,b}(\Omega)} \int_{\Gamma } w\,  d \Hausone.$$
	Here, $G_{a,b}(\Omega)$ is the set of Lipschitz curves in $\Omega$ connecting $a$ and $b$.

	\subsubsection{Reminder of the discretization scheme from~\cite{BBL2020}}
	From a numerical point of view,  as the minimization of $F_{\varepsilon}$ requires the computation of the gradient of the geodesic terms with respect to $u$ and raises some numerical difficulties, a relaxation approach was proposed in~\cite{BBL2020}. It consists in
	introducing an extra variable ${\bf \gamma} := ( \gamma_i )_{1 \leq i \leq N}$, where each $\gamma_i : [0,1] \to \Omega$  is a Lipschitz
	curve joining the base point $a_0$ to the point $a_i$ in $\Omega$, and then consider the relaxed energy $E_{\varepsilon}$ defined by
	$$ E_{\varepsilon}(u, {\bf \gamma}) =  \AT(u)  + R_{\varepsilon}(u,\gamma),\quad  \text{where } R_{\varepsilon}(u,\gamma) = \frac{1}{\lambda_{\varepsilon}} \sum_{i=1}^{N} \int_{\gamma_i([0,1])} (\delta_{\varepsilon} + u^2) d\Hausone.$$
	In particular, this approach can be viewed as a relaxation since
	$$F_{\varepsilon}(u) = \inf_{{\bf \gamma \in \mathcal{P}_{{\bf a}}}} \left\{  E_{\varepsilon}(u,\gamma) \right\},$$
	where $\mathcal{P}_{\bf a}$ is the set of $N$-uplets of Lipschitz curves defined by 
	$$  \mathcal{P}_{\bf a} :=  \{ \gamma =  (\gamma_i)_{i=1:N} : \gamma_i \in \text{Lip}([0,1];\Omega), \gamma_i(0) = a_0 \text{ and } \gamma_i(1) = a_i \}.$$
	An  advantage is that its  $L^2$-gradient flow
	$$
	\begin{cases}
		u_t &= - \nabla_u E_{\varepsilon}(u,\gamma), \\
		\gamma &= \argmin_{\gamma \in \mathcal{P}_{\bf a}} \{  E_{\varepsilon}(u,\gamma) \},
	\end{cases}
	$$
	can be approximated at time $ t = n \delta_t$ by the solution $(u^{n}, \gamma^n)$ of a time-decoupled scheme which alternates
	\begin{itemize}
		\item[-] a geodesic computation:
		$$ \gamma^{n} = \argmin_{\gamma \in \mathcal{P}_{\bf a}} E_{\varepsilon}(u^n,\gamma) =  \argmin_{\gamma \in \mathcal{P}_{\bf a}} R_{\varepsilon}(u^n,\gamma);$$
		\item[-] a phase field evolution: $u^{n+1} \simeq v(\delta)$ where  $v$ is the solution of the PDE:
		$$  v_t = - \nabla_u E_{\varepsilon}(v,\gamma^n) = \varepsilon \Delta v - \frac{1}{\varepsilon} V'(v)  - \frac{2}{\lambda_{\varepsilon} }  \left[  \sum_{i=1}^{N}  \Hausone|_{\gamma^n_i([0,1])} \right] v,$$
		with $ v(0) = u^{n}$ and where the potential $V$ is defined $\forall s \in \R$ by $V(s) = \frac{1}{4}(1-s)^2$.
	\end{itemize}

	\paragraph{Regularization of the geodesic term.} The method presented above makes it possible to avoid differentiating the geodesic term with respect to the phase field function u.
	However, another difficulty arises from the lack of regularity of the solution $v$, which is a consequence of the singularity of the geodesic contribution. $ \sum_{i=1}^{N}  \Hausone|_{\gamma^n_i([0,1])}$.
	In particular, without additional regularization, the numerical approximation of the solution is likely to be highly dependent on the spatial discretization used, with non-negligible anisotropic effects. An initial solution proposed in the work~\cite{BBL2020} consisted of filtering this geodesic term using a convolution kernel
	$$
	v_t = \varepsilon \Delta v - \frac{1}{\varepsilon} V'(v)  - \frac{2}{\lambda_{\varepsilon} } \omega_{\varepsilon,\alpha}[\gamma] v ,\quad  \text{where }
	\omega_{\varepsilon,r}[\gamma] = \sum_{i=1}^{N} \left(\rho_{\varepsilon^{r}}* \Hausone \mres {\gamma_i([0,1])}\right), 
	$$
	and $\rho_{\varepsilon^{r}}$ if a kernel of size $\varepsilon^{r}$, i.e. $\rho_{\varepsilon} = \frac{1}{\varepsilon^2} \rho(\cdot/\varepsilon)$.
	Notice that this new term is still variational and corresponds to the geodesic term
	$$ R_{\varepsilon,r}(u,\gamma) =  \frac{1}{\lambda_{\varepsilon}}  \int_{\Omega}  \omega_{\varepsilon,r}[\gamma] (\delta_{\varepsilon} + u^2) dx = \frac{1}{\lambda_{\varepsilon}}  \sum_{i=1}^{N} \int_{\gamma_i([0,1])} \rho_{\varepsilon^{r}}*(\delta_{\varepsilon} + u^2) d\Hausone. $$

\begin{figure}[!htbp]
		\centering
		\includegraphics[width=.24\textwidth]{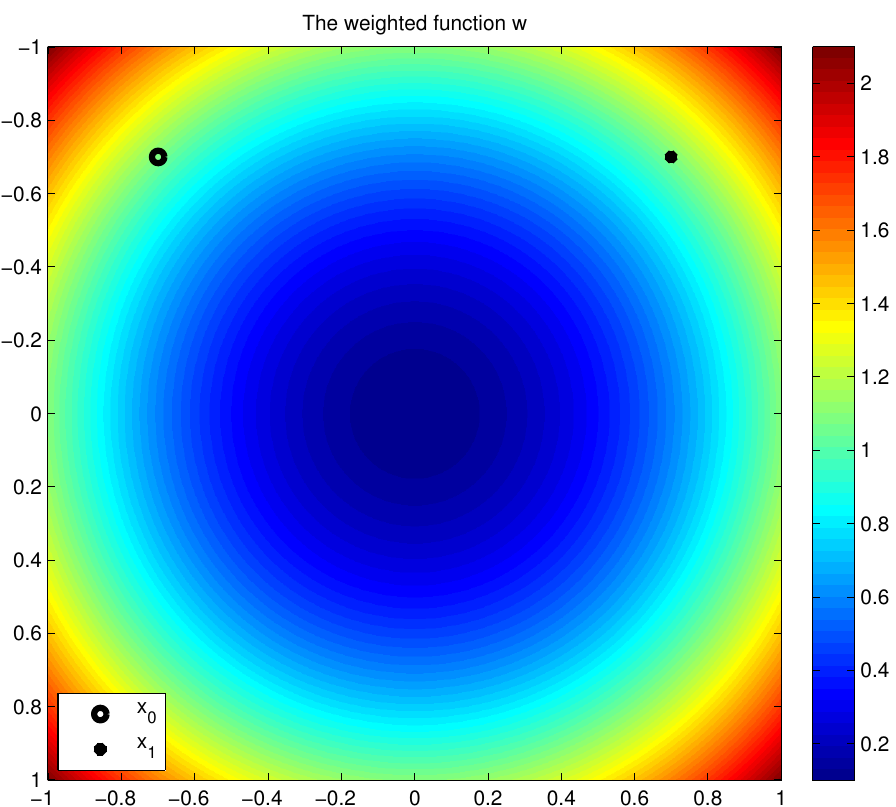}
		\includegraphics[width=.21\textwidth]{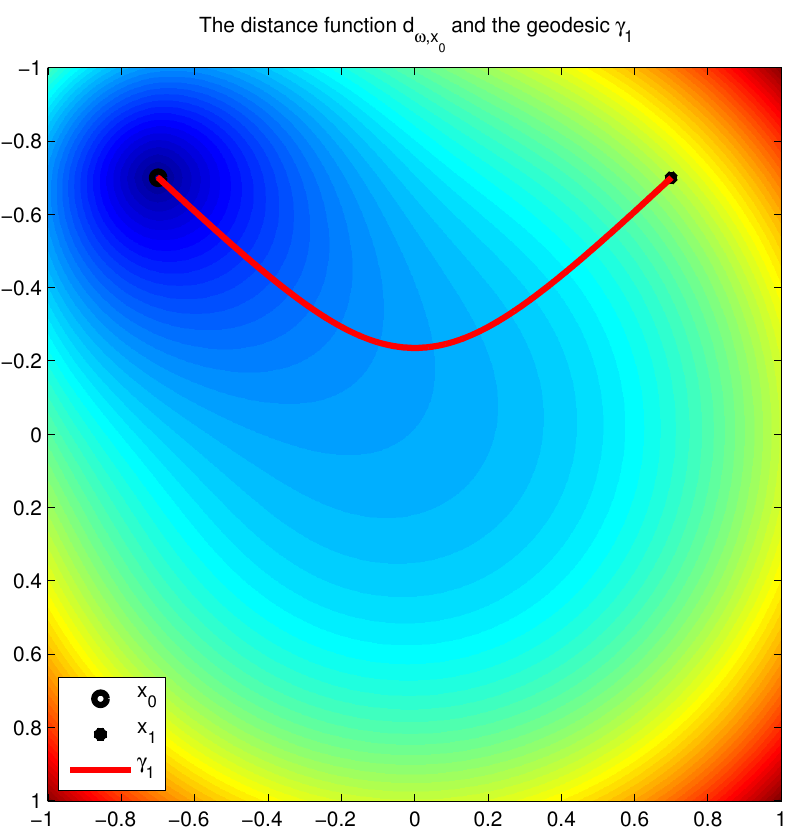}
		\caption{Numerical example of geodesic computation. Left: the weight function
			$w = \| x \|^2$ and the endpoints $a_0 = [-0.7, 0.7]$ and $a_1 = [0.7, -0.7]$.
			Right: the distance function $d_{\omega,a_0}$ computed using a fast Marching algorithm, and the geodesic $\gamma_1$ between $a_0$ and $a_1$.}	\label{fig:Geodesic}
	\end{figure}

	\paragraph{Computation of  $\gamma^n$.} The computation of each geodesic $\gamma^n = ( \gamma^n_i)_{1 = 1:N}$ defined as   
	$$\gamma^{n} =  \argmin_{\gamma \in \mathcal{P}_{\bf a}} R_{\varepsilon}(u^n,\gamma),$$
	can then be carried out in two steps:
	\begin{itemize}
		\item[1)]   Use the Fast Marching Method~\cite{fast_marching_sethian,1999lsmf_book} to compute the weighted distance function $x \mapsto d_{a_0,\omega}(x)$ corresponding to the distance function from $a_0$ to $x$,
		associated to the weight $w =  \rho_{\varepsilon^{r}}*(\delta_{\varepsilon + u^2})$. Notice that it can be defined as a viscosity solution of the non linear Eikonal equation
		$$ | \nabla  d_{a_0,\omega}(x) | = \omega(x), \text{ in } \Omega, \text{ with }   d_{a_0,\omega}(a_0) = 0. $$
		In practice, we use the Toolbox Fast Marching proposed by G. Peyre in the Matlab environment {\it http://www.mathworks.com/matlabcentral/fileexchange/.}
		\item[2)]  Compute each geodesic $\gamma_i:[0,1] \mapsto \Omega$ satisfying $\gamma^i(0)= a_0$ and $\gamma^i(1)= a_i$ by considering a discrete version of the flow
		$$ (\tilde{\gamma}^i)'(s) = - \nabla d_{a_0,\omega}(\tilde{\gamma}^i), \text{ with }  \tilde{\gamma}^i(0) = a_i.$$
		Here, we can use the  {\it Matlab} function {\it compute$\_$geodesic} from the Fast Marching Toolbox to obtain an approximation of each geodesic.
		Figure~\ref{fig:Geodesic} shows an example of distance function $d_{a_0,\omega}$ and the associated estimate geodesic $\gamma_1$ connecting the points
		$(a_0,a_1)$.
		
	\end{itemize}
	
	\paragraph{Computation of  $u^{n+1}$.} As far as the numerical resolution of PDEs is concerned,  we assume here that the computation box $\Omega$ is a square $\Omega= [0,1]^2$ and that all PDEs are considered 
	with additional periodic boundary conditions.  In order to get a high accuracy approximation in space, we also use a semi-implicit convex-concave Fourier-spectral method~\cite{EDTK_Eyre,ShenWWW12} where for instance $u^{n+1}$ can be defined as follows:
	
	$$  \frac{u^{n+1}- u^{n}}{\delta_t} =  \varepsilon \Delta u^{n+1} - \frac{1}{2 \varepsilon}(u^{n+1} - 1) - \alpha u^{n+1} - \left( \frac{2}{\lambda_{\varepsilon}} \omega_{\varepsilon,r}[\gamma^n] - \alpha   \right) u^{n}.$$
	As explained in~\cite{BBL2020}, $\alpha$ can be viewed as a stabilization parameter which ensures the decrease of $E_{\varepsilon}(\cdot, \gamma^{n})$, i.e   $E_{\varepsilon}(u^{n+1}, \gamma^{n}) \leq E_{\varepsilon}(u^{n}, \gamma^{n})$ as soon as
	$$ \alpha \geq \sup_{x \in \Omega} \left\{ \frac{2}{\lambda_{\varepsilon}}   \omega_{\varepsilon,r}[\gamma^{n}]  \right\}. $$
	In practice, we will consider a stabilization parameter $\alpha$ which evolves along the iterations and  take $\alpha^{n} = \sup_{x \in \Omega} \left\{ \frac{2}{\lambda_{\varepsilon}}   \omega_{\varepsilon,r}[\gamma^{n}]  \right\}$. Finally, $u^{n+1}$ is computed using the following formula
	$$ u^{n+1} = \left( I_d + \delta_t (\varepsilon \Delta + \frac{1}{2 \varepsilon} + \alpha^n)  \right)^{-1} \left( u^{n} + \delta_t \left( \frac{1}{2 \varepsilon} - \left( \frac{2}{\lambda_{\varepsilon}}  \omega_{\varepsilon,r}[\gamma^n] - \alpha^n \right) u^n  \right)  \right),$$
	where the operator $\left( I_d + \delta_t (\varepsilon \Delta + \frac{1}{2 \varepsilon} + \alpha^n)  \right)^{-1}$ can be computed easily in Fourier space using Fast Fourier Transform.\\

	\paragraph{Numerical experiments.}
	
\begin{figure}[!htbp]
		\centering
		\includegraphics[width=.21\textwidth]{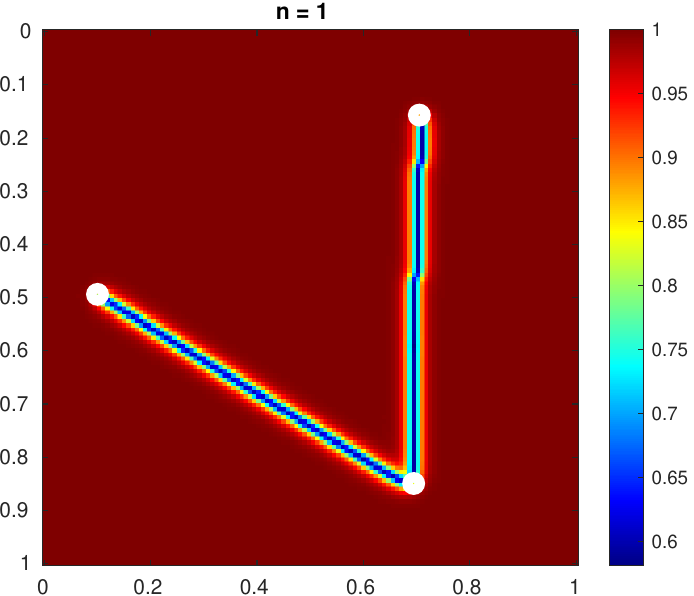}
		\includegraphics[width=.21\textwidth]{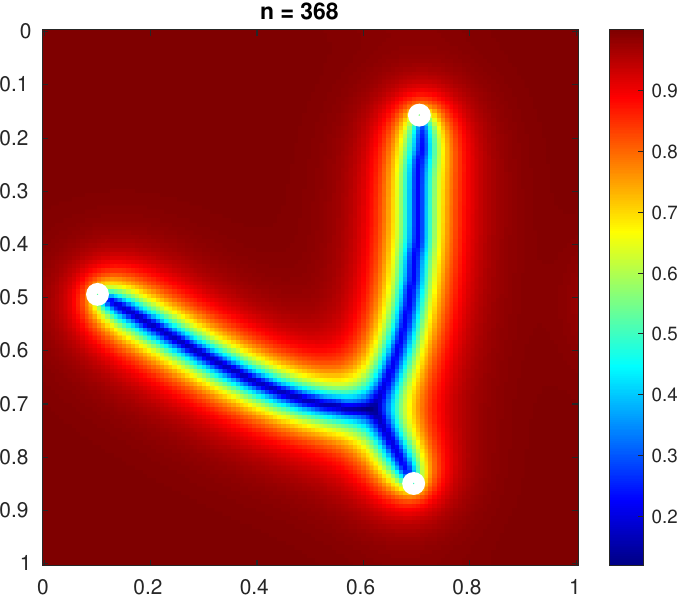}
		\includegraphics[width=.21\textwidth]{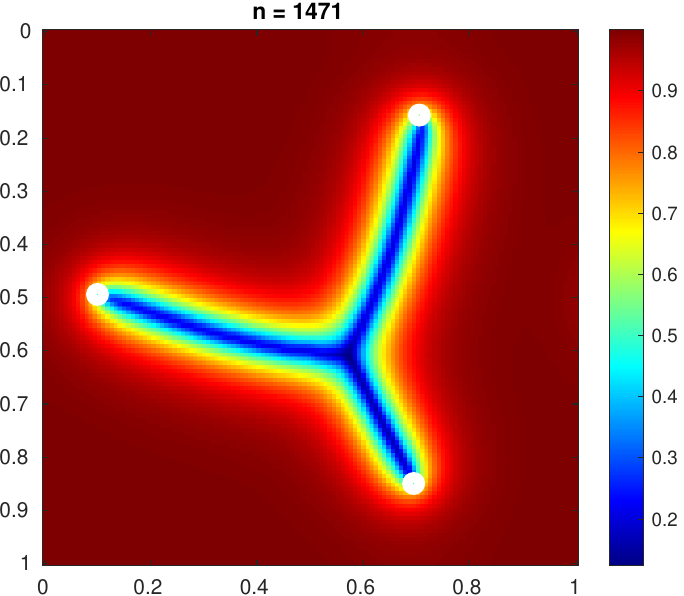}
		\includegraphics[width=.21\textwidth]{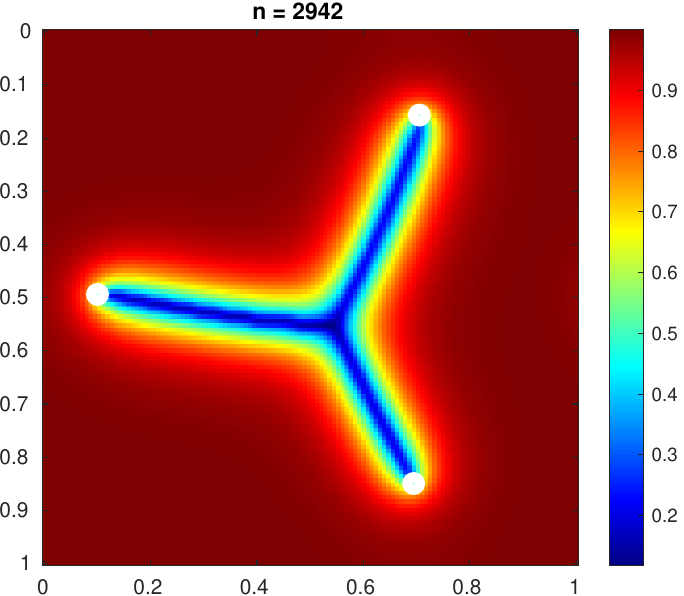} \\
		
		\includegraphics[width=.21\textwidth]{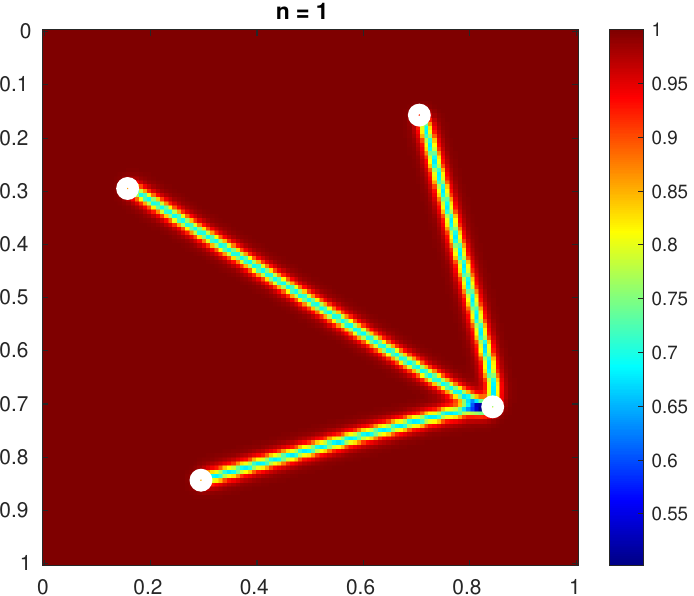}
		\includegraphics[width=.21\textwidth]{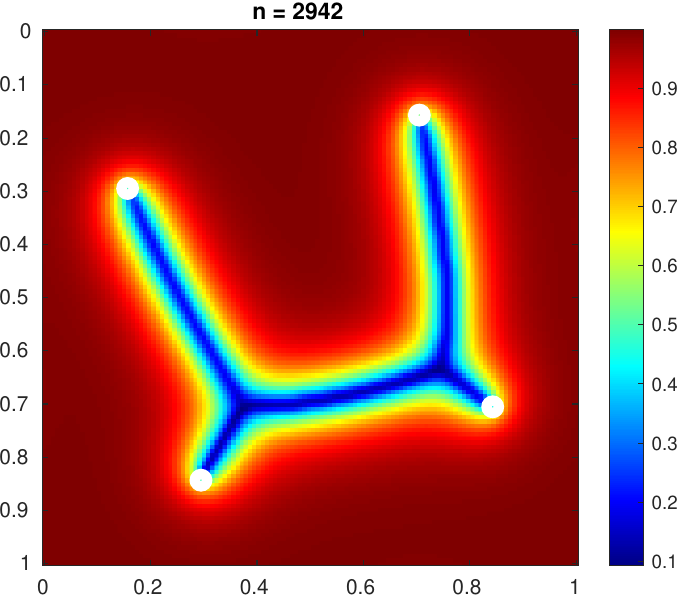}
		\includegraphics[width=.21\textwidth]{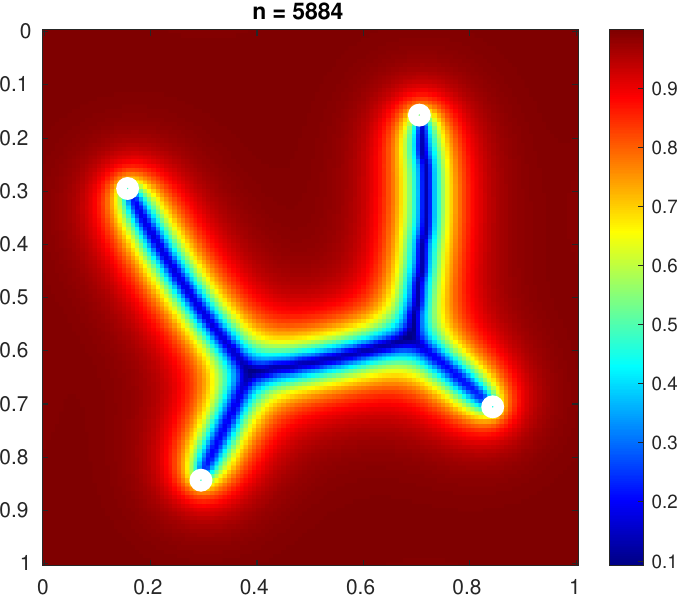}
		\includegraphics[width=.21\textwidth]{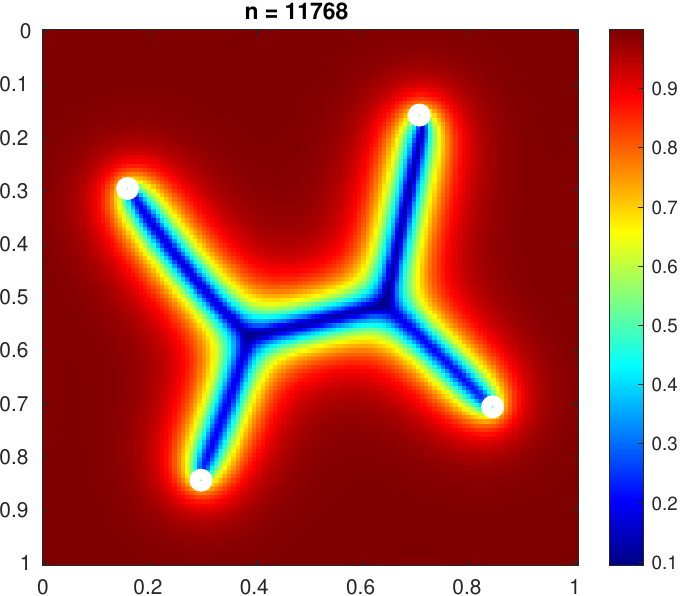} \\
		
		\caption{ Numerical experiments to approximate solutions of Steiner's problem.
			First line with $3$ points, second line with $4$ points. In each line, we plot
			the phase field function $u$ at different time $t_n$ along the iterations.}
		\label{fig:steiner_ambrosio tortorelli}
	\end{figure}

	Figure~\ref{fig:steiner_ambrosio tortorelli}  shows two examples of numerical approximations of solutions to Steiner's problem in the case of $3$ and $4$ points uniformly distributed over a circle. These simulations were carried out with a spatial discretization step of $\delta_x = 1/P$ with $P= 2^7$.
	We have also chosen an approximation parameter $\varepsilon = 4 \delta_x $. Without going into detail about the choice of other simulation parameters, these two experiments show that the numerical approach does indeed provide a relatively good approximation to the solution of Steiner's problem. On the other hand, it was difficult to use a smaller $\varepsilon$ parameter without resorting to a multi-resolution approach, where the $\varepsilon$ parameter could decrease during the iterations. Another remark concerns the phase field solution $u$, which appears rather singular and poorly localised around the interface $K$, with an expected profile of the form
	$$ u_{\varepsilon}(x) \simeq  1-  \exp \left( \frac{- dist(x,K)}{2\varepsilon} \right),$$
	which admits  singularity on  $K$.
	
	\subsubsection{ A new Willmore Cahn-Hilliard phase field model  }

	Although the previous numerical experiments have all shown the ability of our phase field model to approximate the solutions to Steiner's problem, the minimization of the $E_{\varepsilon}$ function itself raises a number of practical difficulties.
	\begin{itemize}
		\item Firstly, the solution $u$  of the phase field model is not very regular, which leads to a strong influence of the discretization parameter field on the solution itself. Therefore, we now propose to modify the Ambrosio-Tortorelli length approximation term to improve the regularity of  the associated profile.
		
		\item One of the major difficulties in the optimisation of $E_{\varepsilon}$ is that it contains $2$ terms whose action is ambivalent, which leads to ill-conditioned problems. It is indeed the alternating action of computing the geodesic and minimizing the phase field model that causes the interface to evolve towards Steiner's solution, but each of the two energies $\AT$ and $R_{\varepsilon}$ does not independently allow the interface to be moved in order to reduce its length. Here, we propose to use a Willmore Cahn-Hilliard energy~\cite{bretin2024}, which approximates the perimeter of the set $K$ and whose $L^2$-gradient flow  allows us to approximate mean curvature flow of interface $K$.
	\end{itemize}

	The aim of this section is to show how replacing the  term $\AT$ by this new phase field model leads to a much more numerically efficient model, without fundamentally changing the spirit of our phase field model. In particular, the $\Gamma$-convergence result  should always be true and obtained in a similar way.
	
	\paragraph{Phase field approximation of mean curvature flow in non oriented interface.}
	
\begin{figure}[!htbp]
		\centering
		\includegraphics[width=.21\textwidth]{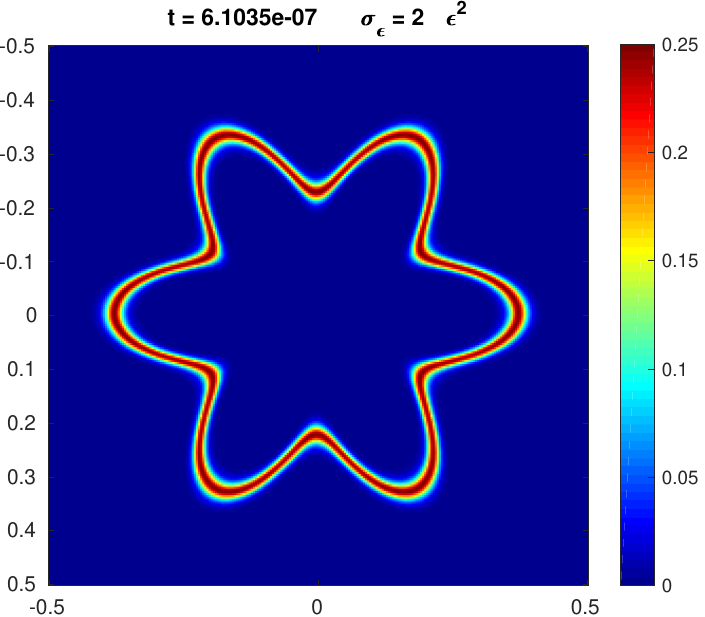}
		\includegraphics[width=.21\textwidth]{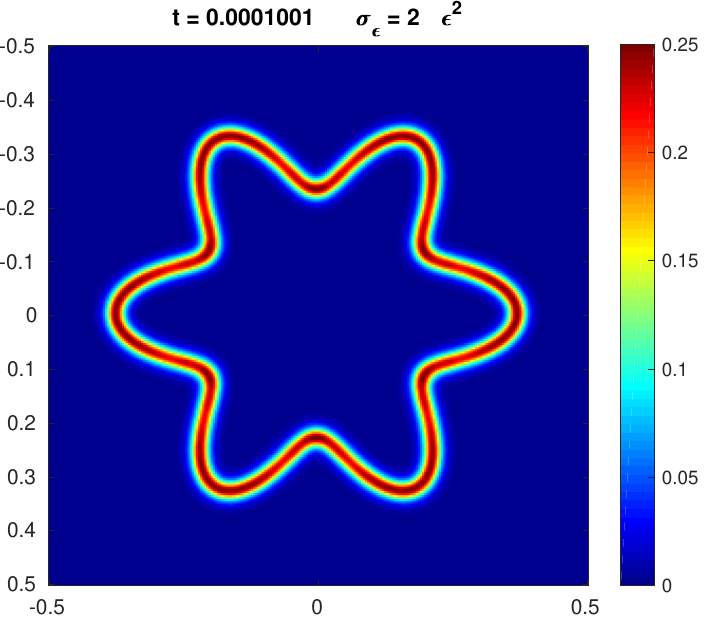}
		\includegraphics[width=.21\textwidth]{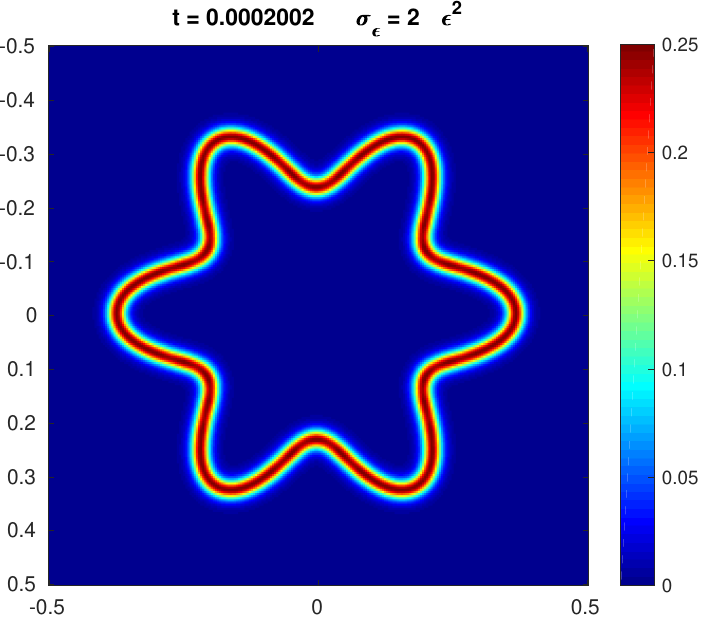}
		\includegraphics[width=.21\textwidth]{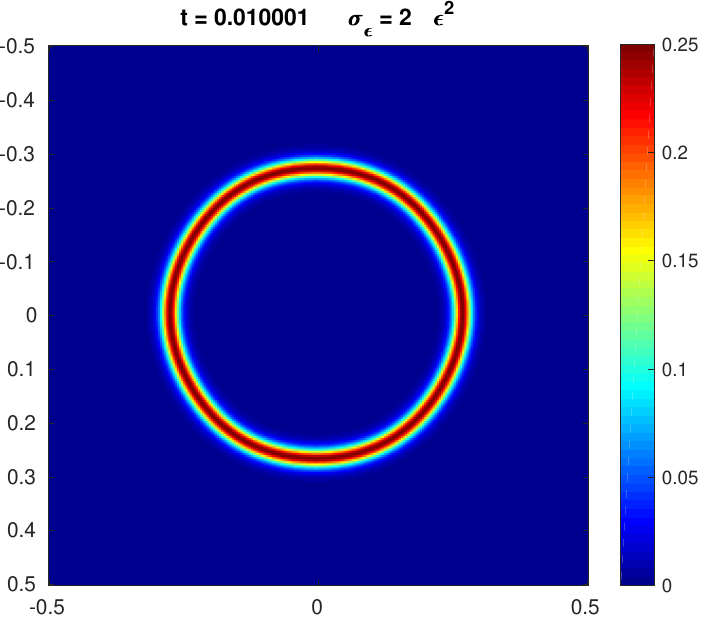}
		
		\caption{ Example of numerical approximation of  mean curvature flow using the $\varepsilon$-gradient flow of
			$P_{\varepsilon}$. The phase field function $u$ is plotted along the iterations. }
		\label{fig:non_oreinted_mean_curvature}
	\end{figure}

	In~\cite{bretin2024}, the authors propose a $\Gamma$-convergence analysis of the modified Cahn-Hillard functional
	$$ P_{\varepsilon}(u) = \int_{\Omega} \left( \frac{\varepsilon}{2} |\nabla u|^2 + \frac{1}{\varepsilon} W(u)  \right) dx + \frac{\sigma_{\varepsilon}}{2 \varepsilon} \int_{\Omega} \left( \varepsilon \Delta u  - \frac{1}{\varepsilon} W'(u) \right)^2 dx,$$
	where the real novelty concerns the non smooth potential $W$ which is defined as
	$$ W(s) = \begin{cases}
		s^2 (1/2 - 2s)  &\text{ if } s \leq \frac{1}{4}, \\
		+ \infty  &\text{ otherwise}.
	\end{cases}
	$$
	Here, $\sigma_{\varepsilon}>0$ satisfies $\varepsilon^2/\sigma_{\varepsilon} \to 0$ to obtain a phase field profile sufficiently stable.
	Notice that the potential $W$  admits a smooth well in $s=0$ and an obstacle at $s = 1/4$. More precisely, it has been constructed
	so that the derivative $y = -q'$ of the standard phase field profile $q(s) = \frac{1}{2}(1 - \tanh(s/2)))$ satisfies the equation
	$$ |y'(s)| = \sqrt{2 W(q)}, \text{ and }  y''(s) =  W'(y(s)).$$
	An optimal minimizing sequence $u_{\varepsilon}$ of $P_{\varepsilon}$ is then expected of the form
	$$ u_{\varepsilon}(x) = -q'\left( \frac{dist(x,K)}{\varepsilon} \right) = \frac{1}{4}\left(1 - \tanh\left(\frac{dist(x,K)}{2 \varepsilon}\right)^2 \right),$$
	as soon as the contribution of the second order term of $P_{\varepsilon}$, which can be viewed as an approximation of the Willmore energy,
	is sufficiently important to stabilize the profile $y = -q'$. In particular, the authors show in~\cite{bretin2024} that the $\varepsilon$-gradient flow of $P_{\varepsilon}$ which reads
	$$u_{t} =   \Delta u - \frac{1}{\varepsilon^2} W'(u) + \sigma_{\varepsilon} \left( - \Delta (\Delta u - \frac{1}{\varepsilon^2} W'(u))  + \frac{W''(u)}{\varepsilon^2} ( \Delta u - \frac{1}{\varepsilon^2} W'(u)) \right),$$
	can be used to approximated  the non oriented mean curvature flow. Without going into detail about the simulation parameters and the choice of numerical discretization, Figure~\ref{fig:non_oreinted_mean_curvature} shows a numerical approximation example of a mean curvature flow based on the modified Cahn Hilliard energy $P_{\varepsilon}$. Note that the profile is stable and the minimization of $P_{\varepsilon}$ allows the interface $K$ to evolve, which was clearly not the case by considering the $L^2$-flow of Ambrosio-Tortorelli energy.

	\paragraph{Application to the case of Steiner's problem and discretization.}
	
	We now propose to replace $\AT$ by  $P_{\varepsilon}$. One just has to be careful about the profile, which does not have the same properties. In the case of $\AT$, the profile vanishes  on the interface $K$ and takes the value $s=1$ far from it, whereas for $P_{\varepsilon}$, the  profile $u$ vanishes far from $K$ and takes the value $s = \frac{1}{4}$ on it. We thus consider the relaxed version $\tilde{E}$ defined by
	$$ \tilde{E}_{\varepsilon,r}(u, {\bf \gamma}) =  P_{\varepsilon}(u)  + R_{\varepsilon,s}(1-4u,\gamma).$$
	In particular, its $L^2$-gradient flow reads
	$$
	\begin{cases}
		u_t &= - \nabla_u \tilde{E}_{\varepsilon}(u,\gamma), \\
		\gamma &= \argmin_{\gamma \in \mathcal{P}_{\bf a}} \{  \tilde{E}_{\varepsilon}(u,\gamma) \}.
	\end{cases}
	$$ 
	The solution $(u^{n}, \gamma^n)$ at time $ t = n \delta_t$
	can be still approximated  using  a time-decoupled scheme, which alternates
	\begin{itemize}
		\item[-] a geodesic computation:
		$$ \gamma^{n} = \argmin_{\gamma \in \mathcal{P}_{\bf a}} \tilde{E}_{\varepsilon}(u^n,\gamma) =  \argmin_{\gamma \in \mathcal{P}_{\bf a}} R_{r,\varepsilon}(1-4u^n,\gamma),$$
		\item[-] a phase field evolution: $u^{n+1} \simeq v(\delta)$ where  $v$ is the solution of the PDE:
		\begin{align*} 
			& v_t = \varepsilon \left( \Delta u - \frac{1}{\varepsilon^2} W'(u) + \sigma_{\varepsilon} \left( - \Delta (\Delta u - \frac{1}{\varepsilon^2} W'(u))  + \frac{W''(u)}{\varepsilon^2} ( \Delta u - \frac{1}{\varepsilon^2} W'(u))) \right) \right) \\
			& \qquad \qquad - \frac{2}{\lambda_{\varepsilon} } \omega_{\varepsilon,r}[\gamma^n] (4 v - 1)   ,
		\end{align*}
		with $ v(0) = u^{n}$.
	\end{itemize}
	
	We compute $u^{n+1}$ by a semi-implicit convex-concave Fourier-spectral method which reads
	$u^{n+1} = L[g(u^n)]$, where the explicit term $g$ is defined as
	\begin{align*}
		g(u) =  u + \delta_t \varepsilon & \left( - \frac{F'(u)}{\varepsilon^2}  +  \sigma_\varepsilon \left(- \Delta  [F'(u)/\varepsilon^2] + \frac{F''(u)}{\varepsilon^2} \left[ \Delta u - F'(u)/\varepsilon^2 \right] \right)  + \alpha u - \beta \Delta u  \right.\\
		& \quad \left.- \frac{2}{\varepsilon \lambda_{\varepsilon} } \omega_{\varepsilon,r}[\gamma^n] (4u  - 1)  \right),
	\end{align*}
	and where the homogeneous linear operator $L$ satisfies
	$L = \left( I_d +  \varepsilon \delta_t (- \Delta + \sigma_\varepsilon \Delta^2 + \alpha I_d - \beta \Delta  ) ) \right)^{-1}$.
	Here,  $\alpha$ and $\beta$ can be viewed as stabilization parameters which  need to be chosen sufficiently large to guarantee  the decrease of $\tilde{E}_{\varepsilon}(\cdot, \gamma^{n})$.
	
	\paragraph{Numerical experiments.}
	
\begin{figure}[!htbp]
		\centering
		\includegraphics[width=.21\textwidth]{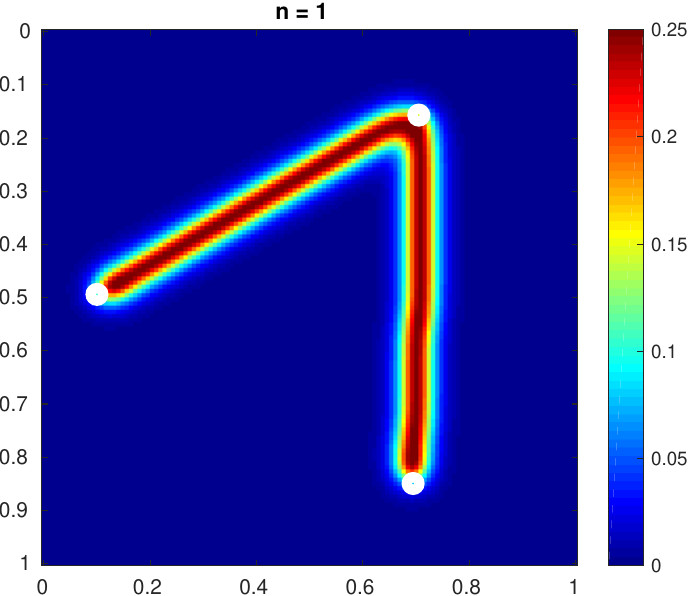}
		\includegraphics[width=.21\textwidth]{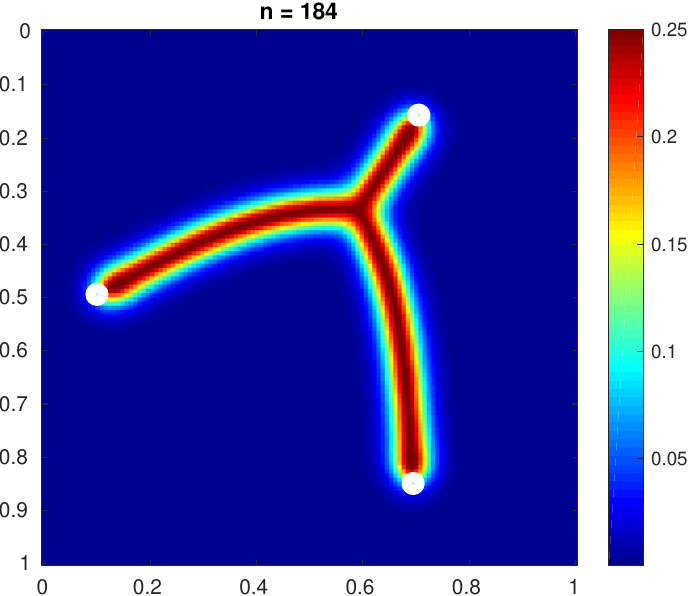}
		\includegraphics[width=.21\textwidth]{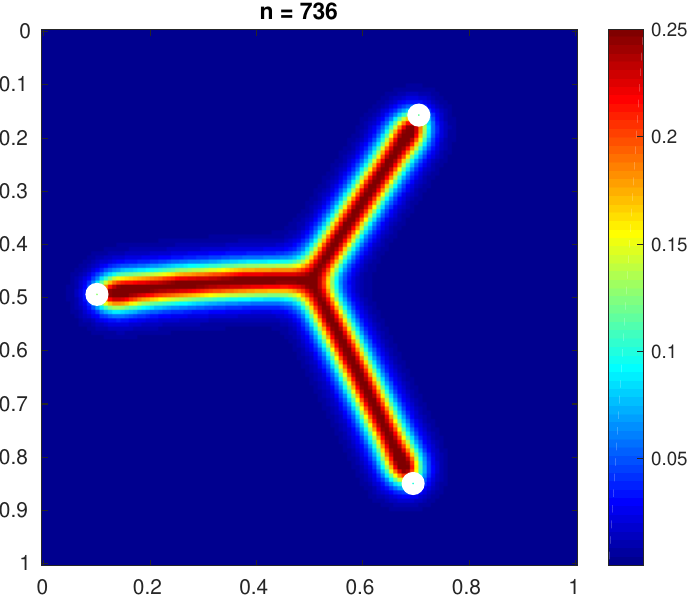}
		\includegraphics[width=.21\textwidth]{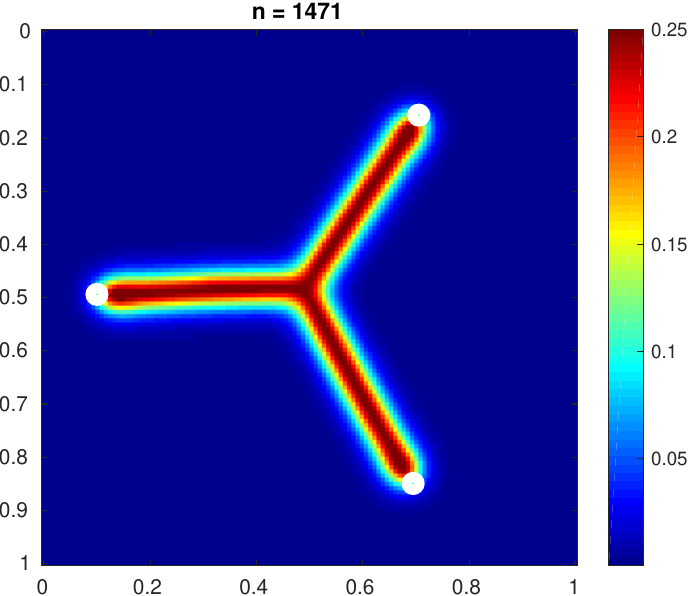} \\
		\includegraphics[width=.21\textwidth]{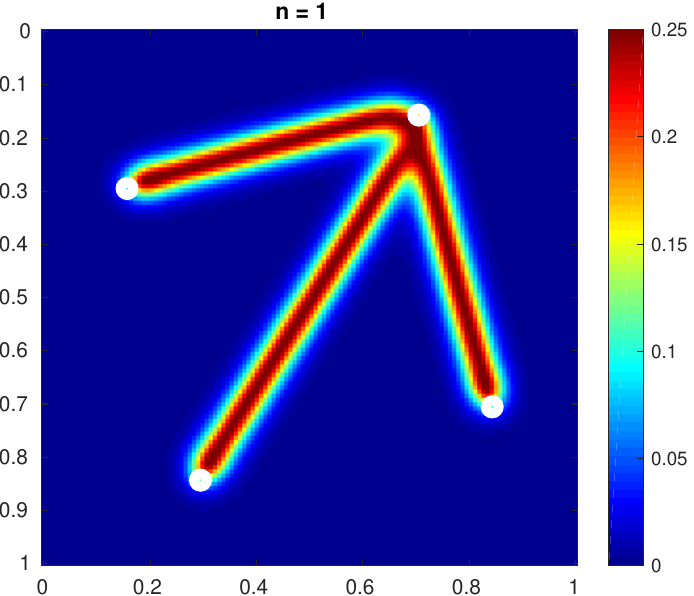}
		\includegraphics[width=.21\textwidth]{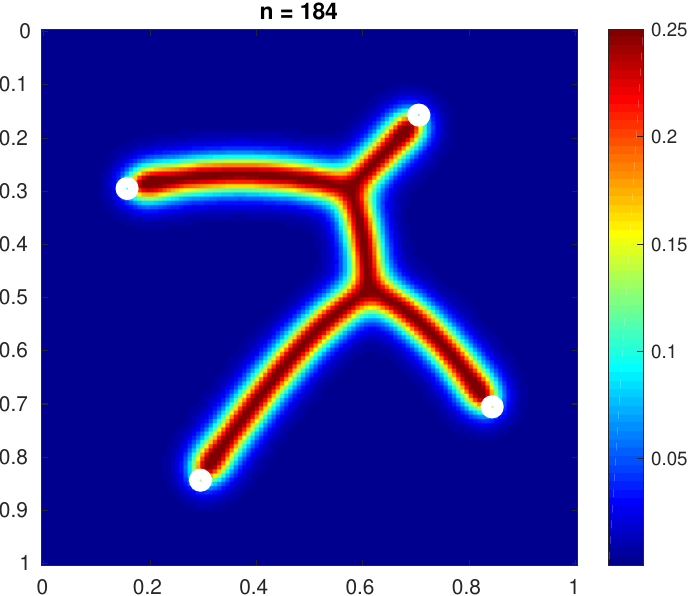}
		\includegraphics[width=.21\textwidth]{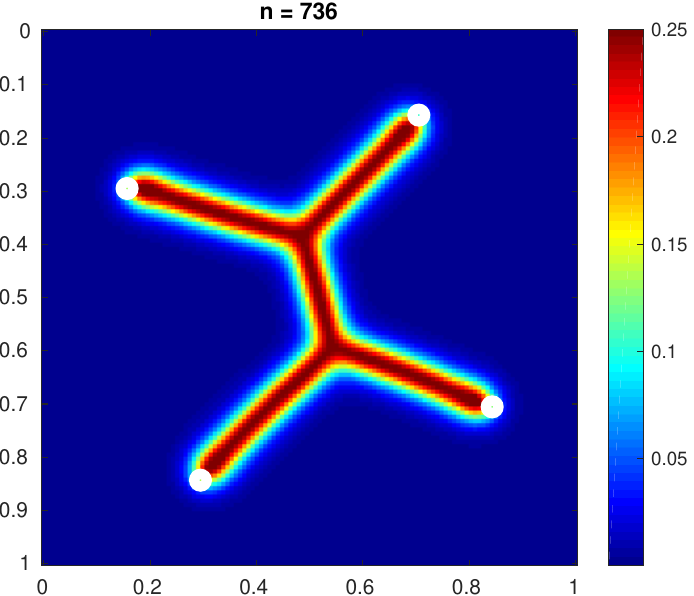}
		\includegraphics[width=.21\textwidth]{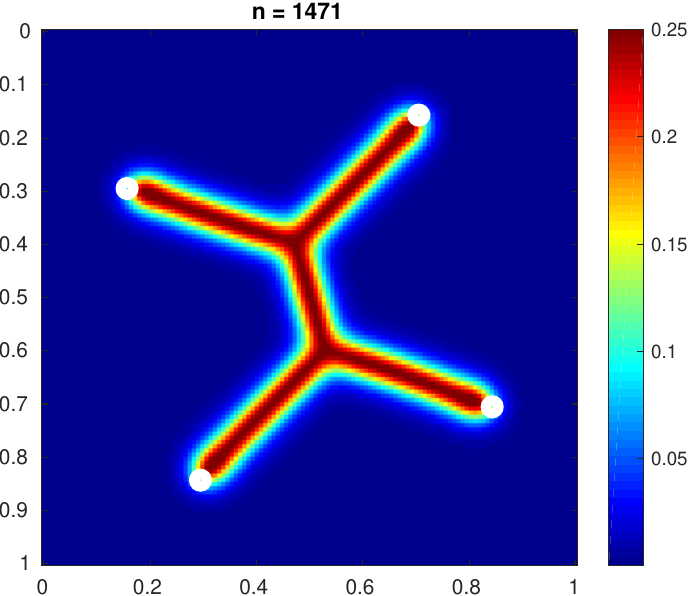}
		
		\caption{ Numerical experiments to approximate solution of Steiner problem using the second order Cahn-Hilliard model;
			First line with $3$ points, second line with $4$ points. In each line, we plot
			the phase field function $u$ at different time $t_n$ along the iterations }
		\label{fig:Steiner_qprim}
	\end{figure}

	In Figure~\ref{fig:Steiner_qprim}, we present two numerical experiments   analogous to the cases presented in Figure~\ref{fig:steiner_ambrosio tortorelli}.
	For these experiments, we are still using a spatial resolution given by $P = 2^7$. However, with this new phase field model, it is possible to use smaller $\varepsilon$ parameters and, in particular, we have used $\varepsilon = 2/P$. We first point that from a qualitative point of view, the numerical results obtained with this method seem comparable with the ones obtained by our original  model (see Figure~\ref{fig:steiner_ambrosio tortorelli}). However, the interface appears to localize better around the Steiner tree, probably due to the choice of the $\varepsilon$ parameter, which is twice as small. Note also that the solution $u$ seems smoother, which was expected with this new profile. 
	From a convergence rate point of view, the use of this new phase field model makes it possible to obtain approximations to Steiner's solution more quickly. In the case of three points, for example, we manage to obtain a correct solution after $1000$ iterations only, which is clearly not the case using $\AT$ term, even though the $\varepsilon$ parameter is smaller.

	In conclusion, this new model appears to be highly beneficial in all respects, justifying our use of this modified version to approximate the numerical solutions to Plateau's problem.

	\begin{rmq}
		The stabilisation of the profile is due to the presence of the second order terms in the energy $P_\eps$, which improves the conditioning of the optimisation problem. However, this Willmore term only makes sense for profiles with regular phase fields and could not have been applied in the case of the standard potential $V(s) = \frac{1}{4}(1-s)^2$.  It is therefore the combination of the potential $W$ with the Willmore term that makes it possible to significantly improve Steiner's phase field model.
	\end{rmq}

	\subsection{Plateau's problem and  geodesic distance approximation}\label{Sect:PlateauNumerics}
	
	The numerical results for the approximation of Steiner's problem have shown the significant numerical advantage of using Cahn-Hilliard's approximation $P_{\varepsilon}$ instead of its original version $\AT$. Consequently, we will use this variant in the case of Plateau's problem by subsequently considering the following model:
	$$\widetilde{F}_\eps(u) = P_\eps(u) + \frac{1}{c_\eps}\sum_{(i,j) \in I_{\gamma}}d^2_{1-4u}(\bpar^i,\bpar^j).
	$$
	The idea is to adapt also the same regularization strategy by considering a relaxed and regularized version of the geodesic term. We are primarily interested in a relaxed version which consists of introducing the variable ${\bf \ell} = (\ell_{i,j})_{(i,j) \in I_{\gamma}}$ where each $\ell_{i,j} \in \curve(\bpar^i,\bpar^j)$ and considering the following relaxed energy
	$$ \widetilde{E}_{\eps}(u,{\bf \ell} ) = P_\eps(u) + R_{\eps}(1-4u,{\bf \ell}),\quad  \text{ with } R_{\eps}(u,{\bf \ell}) = \frac{1}{c_{\eps}} \sum_{(i,j)\in I_{\gamma}} \int_{\surf_{\ell_{i,j}}} (\delta_{\eps} + u^2) d \Haus.$$
	As in the case of Steiner, with $\mathcal{P}_{\gamma} = \{ \ell = (\ell_{i,j})_{(i,j)\in I_{\gamma}}\ ;\  \ell_{i,j} \in \curve(\bpar^i,\bpar^j) \}$, one has
	$$ \widetilde{F}_\eps(u)  = \inf_{\ell \in \mathcal{P}_{\gamma}  } \left\{ \widetilde{E}_{\eps}(u,{\bf \ell} ) \right\},$$
	and the geodesic term is regularized and replaced by
	$$ R_{\varepsilon,\alpha}(u,\ell) =  \frac{1}{c_{\varepsilon}}  \int_{\Omega}  \omega_{\varepsilon,r}[\ell] (\delta_{\varepsilon} + u^2) dx = \frac{1}{c_{\varepsilon}}  \sum_{(i,j) \in I_{\gamma}} \int_{\surf_{\ell_{i,j}}} \rho_{\varepsilon^{r}}*(\delta_{\varepsilon} + u^2) d\Haus, $$
	with
	$ \omega_{\varepsilon,r}[\ell] =  \sum_{(i,j)\in I_{\gamma}} \left(\rho_{\varepsilon^{r}}* \Haus\mres \surf_{\ell_{i,j}}\right) .
	$
	It is then possible to effectively minimize the relaxed and regularized energy
	$$ \widetilde{E}_{\eps,\alpha}(u,{\bf \ell} ) = P_\eps(u) + R_{\eps,\alpha}(1-4u,{\bf \ell}),$$
	by considering the discretization of its $L^2$-gradient flow
	$$
	\begin{cases}
		u_t &= - \nabla_u \widetilde{E}_{\varepsilon,\alpha}(u,\ell) \\
		\ell &= \argmin_{\ell \in \mathcal{P}_{\ell}} \{  \widetilde{E}_{\varepsilon,\alpha}(u,\ell) \}
	\end{cases}
	$$
	with a time decoupled scheme $(u^{n},\ell^n)$ which alternates
	\begin{itemize}
		\item[-] {\bf a geodesic computation :}
		$$ \ell^{n} = \argmin_{\ell \in \mathcal{P}_{\ell}} \widetilde{E}_{\varepsilon,\alpha}(u^n,\ell) =  \argmin_{\ell \in \mathcal{P}_{\ell}} R_{r,\varepsilon}(1-4u^n,\ell);$$
		\item[-] {\bf a phase field evolution} : $u^{n+1} \simeq v(\delta)$ where  $v$ is the solution of the PDE:
		\begin{align*}
			v_t&  =  \varepsilon \left( \Delta u - \frac{1}{\varepsilon^2} W'(u) + \sigma_{\varepsilon} \left( - \Delta (\Delta u - \frac{1}{\varepsilon^2} W'(u))  + \frac{W''(u)}{\varepsilon^2} ( \Delta u - \frac{1}{\varepsilon^2} W'(u))) \right) \right)\\
			& \quad\quad   - \frac{2}{\lambda_{\varepsilon} } \omega_{\varepsilon,r}[\ell^n] (4 v - 1)   ,
		\end{align*}
		with $ v(0) = u^{n}$.
	\end{itemize}
	
	The phase field part is very similar to the approaches presented in the case of Steiner's problem and for which we use exactly the same choices of discretization of the PDE to  estimate $v(\delta_t)$. It therefore did not seem appropriate to repeat these algorithms here. On the other hand, the treatment of the geodesic part is quite different here and requires clarification.

	\subsubsection{Computation of the geodesic}
	
	The first point to understand here is that, from a methodological point of view, the step
	$$ \ell^{n} =  \argmin_{\ell \in \mathcal{P}_{\ell}} R_{r,\varepsilon}(1-4u^n,\ell).$$
	requires some approximation and simplification.
	
	Indeed, let us consider the case of  a Plateau problem where the given boundary contains only a single regular curve $\Gamma^1$ which corresponds to taking $I_{\gamma} = \{(1,2)\}$ and $(\gamma^{1},\gamma^{2})$ such that $\gamma^{1}([0,1])=\Gamma^1$ and $\gamma^{2}([0,1])= x_0 \in \Gamma^1$. Using the initial condition $u^0 = 0$,  it is not difficult to see that the first step of our algorithm, i.e. the computation of $\ell^0$ defined by
	$$ \ell^{0} =  \argmin_{\ell \in \mathcal{P}_{\ell}} R_{r,\varepsilon}(1,\ell),$$
	requires determining the geodesic that connects the curves $\gamma^1$ and $\gamma^2$, whose image precisely corresponds to the desired solution of the Plateau problem.
	It therefore seems unrealistic to solve this optimisation problem and, except in very specific cases,  we treat an approximate version that allows us to obtain one curve in the set $\mathcal{P}_{\ell}$. \\


	\paragraph{Case of a circular curve.}
	The first idea is to calculate a geodesic in a small space. For instance, in the case where $\gamma^{1}$ $\gamma^{2}$ correspond to horizontal circles with respective radii $r_1$ and $r_2$ and  positioned at $ X^1 = (0,0,z_1)$ and $X^2 = (0,0,z_2)$, then at least when $u=0$, we can expect  the geodesic $\ell_{1,2}$ to be identified for all $t \in [0,1]$ with horizontal circles of radius $R(t)$ and positioned at   $X(t) = (0,0,z(t))$. In this specific case, the problem can therefore be reduced to a calculation of geodesics in a $2$ dimensional space corresponding to a radius $r$ and a height $z$.
	The next step is to determine  numerically the weight function $\omega : \R^2 \to \R$ defined by
	$$\forall (r,z) \in \R^2\quad   \omega(r,z) = \int_{0}^{2 \pi} \left( 1 - 4 u(r \cos(\theta),r \sin(\theta), z ) \right)^2 r d\theta.$$
	Thus, as in the case of Steiner's problem in dimension $2$, an approximation of the geodesic $\ell$ will be obtained by
	
	\begin{itemize}
		\item[1)]   computing using the Fast Marching Method  the weighted distance function $x \mapsto d_{\overline{\gamma}_1,\omega}(\gamma)$ corresponding to the distance function from $\overline{\gamma}_1 = (r_1,z_1)$ to $\overline{\gamma} = (r,z)$ and
		associated to the weight $\omega$; 
		\item[2)]  computing the geodesic $\overline{\ell}_{1,2}:[0,1] \mapsto \R^2$ satisfying $\overline{\ell}_{1,2}(0)= \overline{\gamma}_1 $ and $\overline{\ell}_{1,2}(1) = \overline{\gamma}_2 $ by considering a discrete  version of the flow
		$$ (\overline{\ell}_{1,2})'(s) = - \nabla d_{\overline{\gamma}_1,\omega}(\overline{\ell}_{1,2}(s)), \text{ with }  \overline{\ell}_{1,2}(0) = \overline{\gamma}_2.$$
	\end{itemize}
	
\begin{figure}[!htbp]
		\centering
		\includegraphics[width=.3\textwidth]{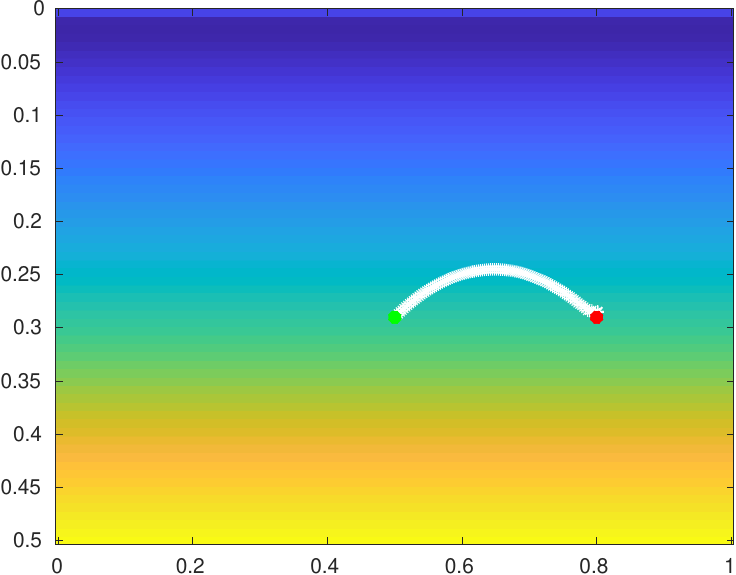}
		\includegraphics[width=.35\textwidth]{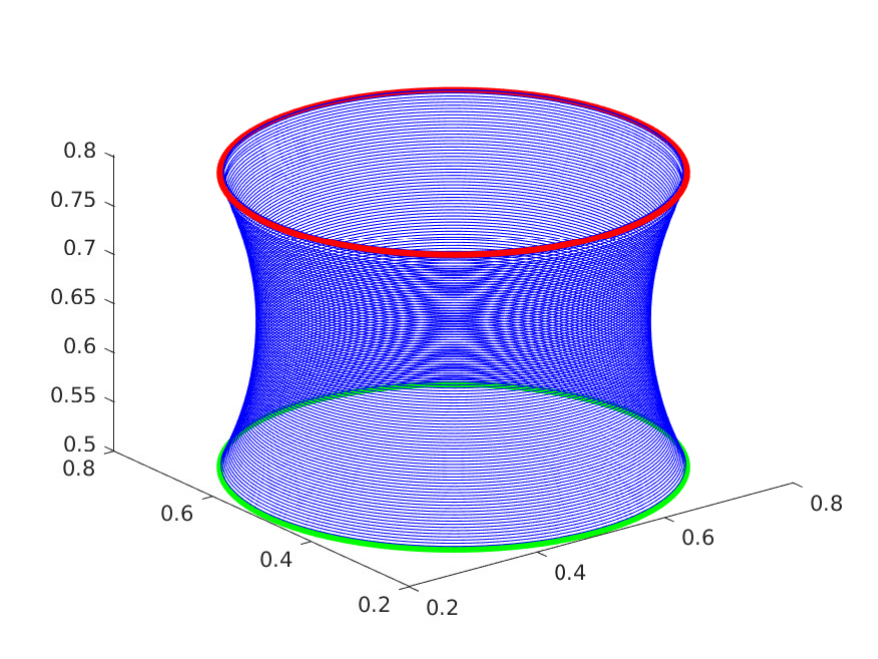}
		\caption{ Example of calculation of a geodesic between two centered and horizontal circles. Left: the weight function $\omega$ and the geodesic $\overline{\ell}$ in the reduced $2$-dimensional space. Right: the geodesic represented in real space in dimension $3$. }
		\label{fig1:plateau_geodesic_espace_reduit}
	\end{figure}

\begin{figure}[!htbp]
		\centering
		\includegraphics[width=.3\textwidth]{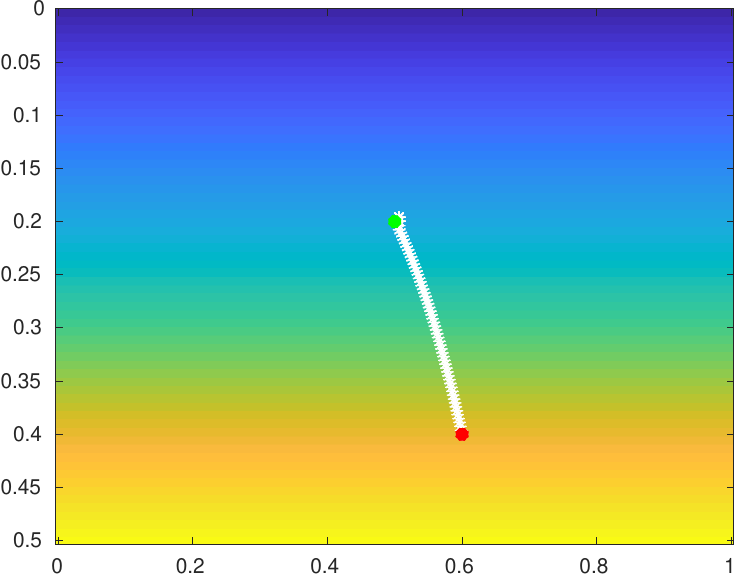}
		\includegraphics[width=.35\textwidth]{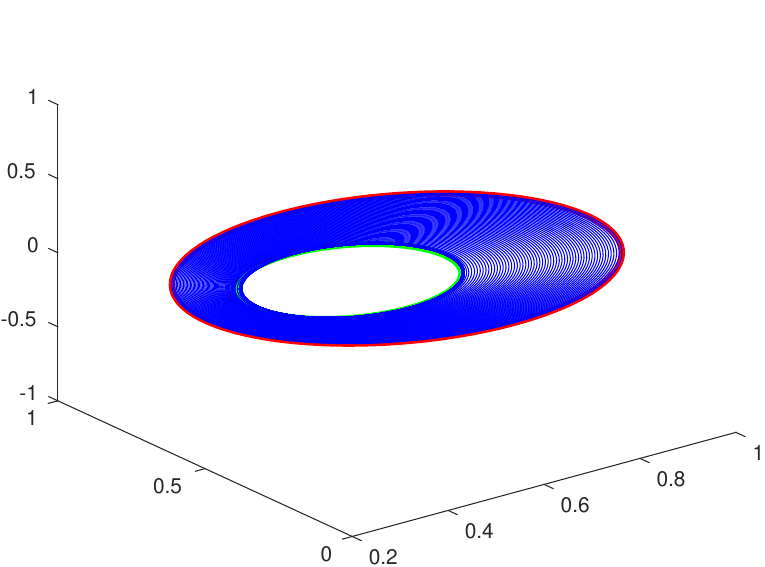}
		\caption{ Example of calculation of a geodesic between two horizontal circles at $z=0$. Left: the weight function $\omega$ and the geodesic $\overline{\ell}$ in the reduced $2$-dimensional space. Right: the geodesic represented in real space in dimension $3$. }
		\label{fig2:plateau_geodesic_espace_reduit}
	\end{figure}

	Figure~\ref{fig1:plateau_geodesic_espace_reduit} shows an example of a geodesic obtained numerically using this reduced space approach. 	In Figure~\ref{fig2:plateau_geodesic_espace_reduit} we present a second example of a geodesic computation. In this case, the reduced space corresponds to the radii of circles of same center $X = (x_1, 0, 0)$, connecting a small circle to a large one in the same horizontal plane.

	This approach produces consistent results, but it is difficult to generalize to more complicated Plateau problems since Fast Marching algorithms are challenging to use in spaces greater than three dimensions. Therefore, we will not use it in the complete model with the phase field approach. \\
	
	\paragraph{General case and approximation.}
	
	The idea is now to examine each curve using the following parametric representation:
	$$ \ell_{1,2}(t) =  \gamma(t) =  \{ \overline{\ell}_{1,2}(t,\theta)  \in \Omega\ ;\  \theta \in [0,2\pi]\}.$$
	once again, the computation of the geodesic will be carried out in two steps:
	\begin{itemize}
		\item[1)]  computation in the real space and using the Fast Marching Method of the weighted distance function $x \mapsto d_{\overline{\gamma}_1,\omega}(\gamma)$ corresponding to the distance function to the curve $\overline{\gamma}_1$ and associated to the weight $\omega = (1-4 u)^2$;
		\item[2)] computation of a geodesic $\ell$ by considering a discrete version of the flow
		$$ \partial_t (\overline{\ell}_{1,2}(t,\theta)) = - \nabla d_{\overline{\gamma}_1,\omega}(\overline{\ell}_{1,2}(t,\theta)), \text{ with }  \overline{\ell}_{1,2}(0,\theta) = \overline{\gamma}_2(\theta), \quad \forall \theta \in [0,2\pi].$$
	\end{itemize}

	\begin{rmq}
		To prevent the discrete approximation of the parametric representation of the curve from revealing discrepancies, which would correspond to a discontinuity, we propose to re-parameterize the curve with respect to the variable $\theta$ when it is no longer sufficiently uniform.
	\end{rmq}
	
\begin{figure}[!h]
		\centering
		\includegraphics[width=.45\textwidth]{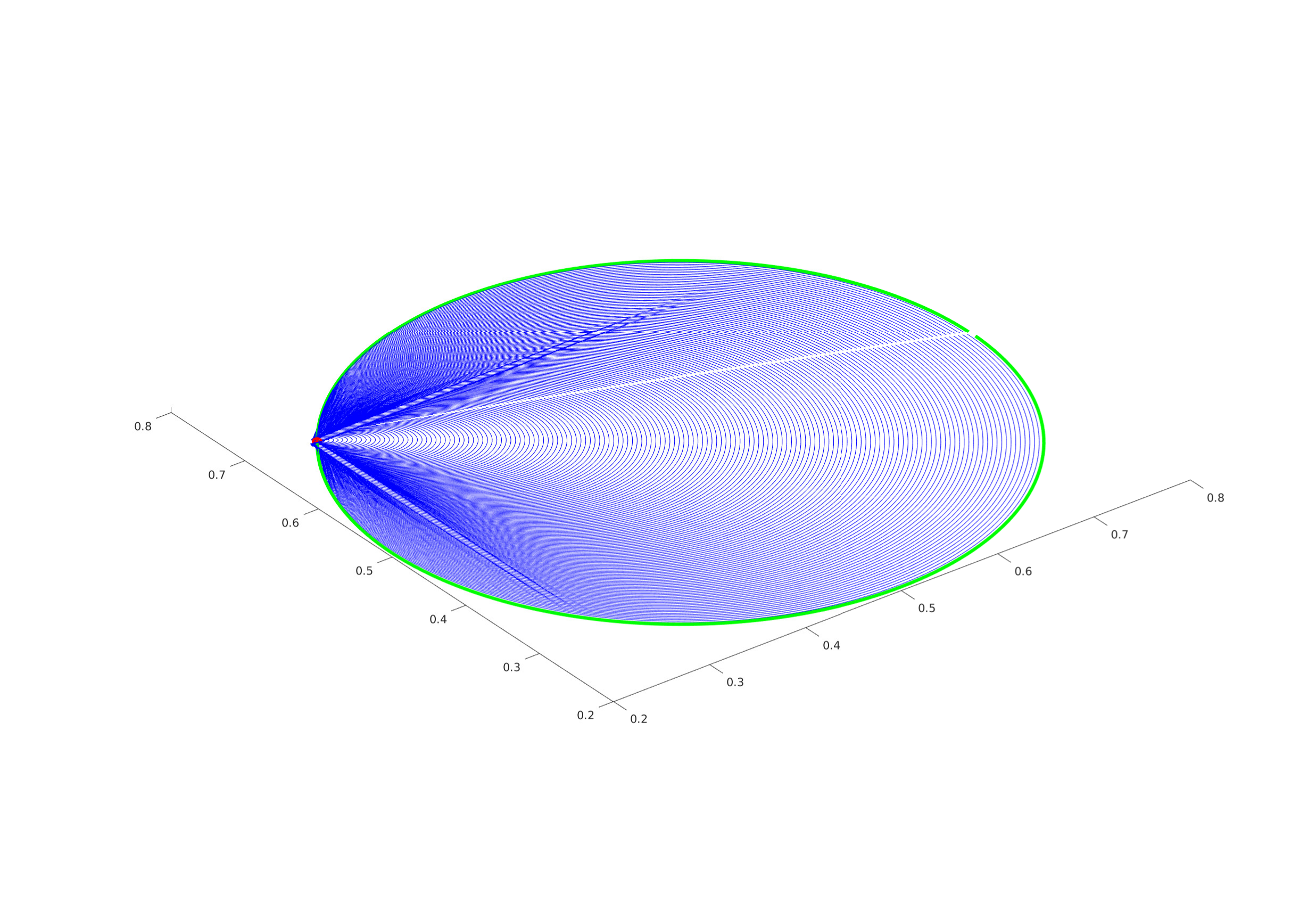}
		\includegraphics[width=.45\textwidth]{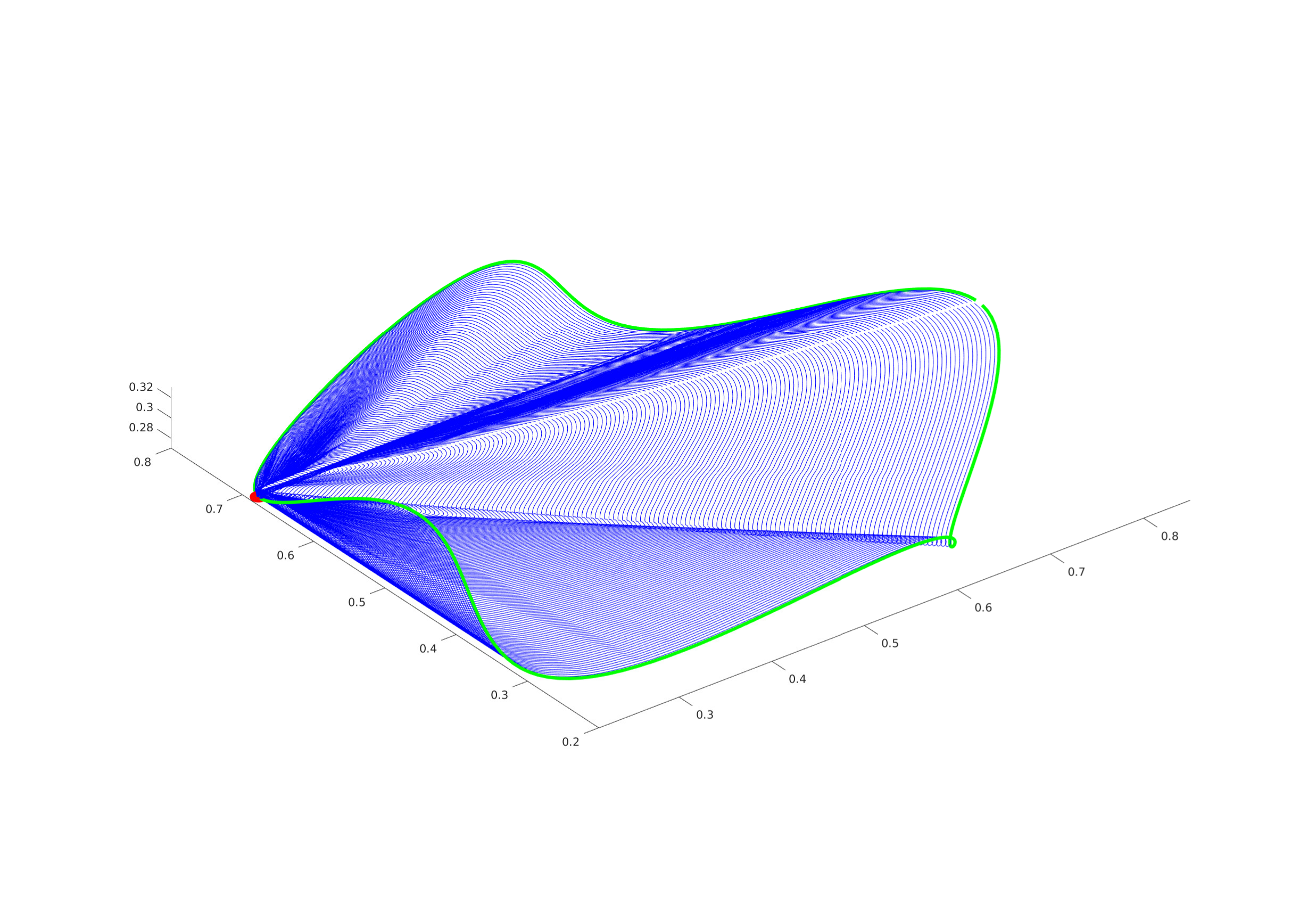}
		\caption{Two numerical examples of geodesic computations
			in the case of a single curve $\Gamma$  with $\gamma^1([0,2\pi]) = \Gamma$ and $\gamma^2([0,2\pi]) = X_2 \in \Gamma$, setting $u=0$.}
		\label{fig:plateau_geodesic_espace_real_onecurve}
	\end{figure}

	Figure~\ref{fig:plateau_geodesic_espace_real_onecurve} shows  examples of geodesic approximations obtained by applying this strategy to the case of a single curve $\Gamma^1$. In the first example, where $\Gamma^1$ is a circle, the computed geodesic is a disk, which effectively corresponds to the solution of Plateau's problem. On the other hand, the second example clearly shows a geodesic that does not correspond to the solution of Plateau's problem. However the computed path $\ell$ connects the green curve $\gamma^1$ to the red point $\gamma^2,$ and therefore properly spans the set $\Gamma^1$. 
	
\begin{figure}[!htbp]
		\centering
		\includegraphics[width=.45\textwidth]{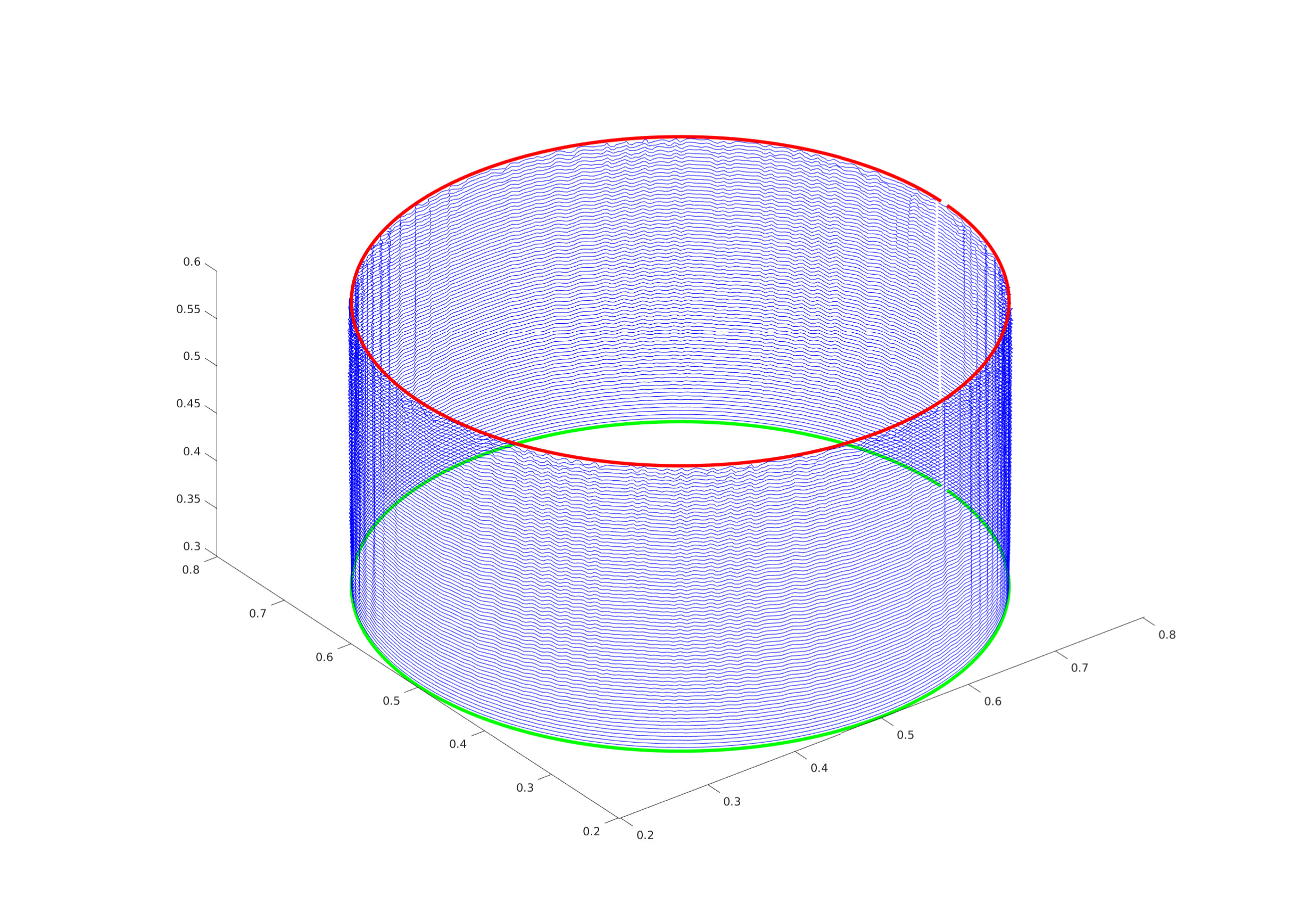}
		\includegraphics[width=.45\textwidth]{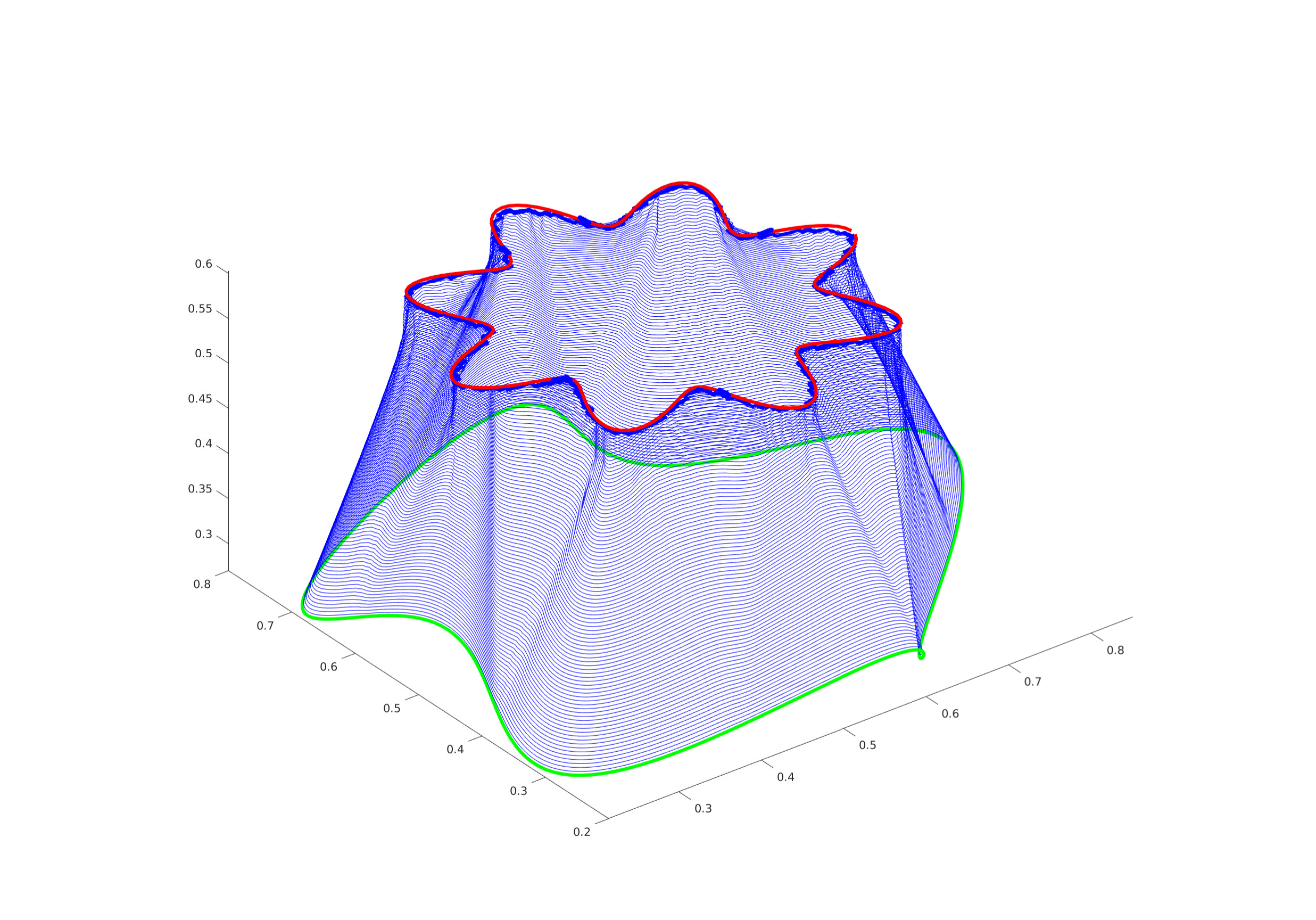}
		\caption{ Two numerical examples of geodesic computations
			in the case of a two curves $\Gamma_1$ and $\Gamma_2$  with $\gamma^1([0,2\pi]) = \Gamma_1$ and $\gamma^2([0,2\pi]) = \Gamma_2 $, setting $u=0$.}
		\label{fig:plateau_geodesic_real_onecurve}
	\end{figure}

	Figure~\ref{fig:plateau_geodesic_real_onecurve} shows examples of geodesic approximations in the case of a Plateau problem with two edges $\Gamma_1$ and $\Gamma_2$. Although this case does is not covered by the Gamma-convergence analysis carried out in Section~\ref{Section:Analysis}, this situation can be dealt with numerically. As before, the geodesic obtained is clearly not the optimal geodesic but does allow us to connect the two green and red curves $\gamma^1$ and $\gamma^2$.

	\begin{rmq}
		As we will see in the next subsection, even if the geodesic computation does not provide the optimal geodesic $\ell$ in the sense of the optimisation problem, this approach allows us to find good approximations to the solution of Plateau's problem by considering the complete phase field model. In a sense, the combined contributions of the phase field and the Willmore Cahn-Hilliard $P_{\eps}$ energy to the minimisation of  the area of the geodesics, yield the convergence to minimal surfaces. 
	\end{rmq}
	
	
	\subsubsection{Numerical experiments}
	
	All the numerical experiments were done  with the same set of numerical parameters, namely  $P = 2^7$,
	$\varepsilon = 2/P$,  $\delta_t = 10\, \varepsilon^2$, $\sigma_{\epsilon} = 1/\varepsilon^2$ and $r = 0.1\,  \varepsilon^2$.
	The scheme applied to the phase field part is exactly the same as that used for the Steiner problem. In particular, we consider all PDEs in the calculation box $\Omega = [0,1]^3$ with additional periodic boundary conditions, using a semi-implicit convex-concave Fourier approach for numerical integration.
	Numerical approximations of geodesic curves are obtained by applying the approximate method, which only enforces the connection between curves $\gamma^i$ and $\gamma^j$. Note that in this case, we cannot guarantee that the optimal surface is calculated, as demonstrated by the example of the two circles, where the cylinder replaces the expected catenoid. Nevertheless, our approach to numerically minimizing the function $F_\varepsilon$ enables us to find good approximations of solutions to Plateau's problem when there are one, two, three or six curves. In particular, we will see that the choice of the $\gamma^i$ curve and how it is connected has a significant impact on the type of Plateau solutions that can be approximated. 
	
	In each figure, we plot the curves $\gamma^i$ in green or red.  The approximation of Plateau's problem is plotted in blue and obtained as a level set of the phase field function $u$:
	$$ S_{\varepsilon} = \{ x \in Q\ ;\ u_{\eps}(x) = 3/16 \}.$$ 
	
	\paragraph{Simple case: one curve, cylindrical case.}
	
\begin{figure}[!htbp]
		\centering
		\includegraphics[width=.32\textwidth]{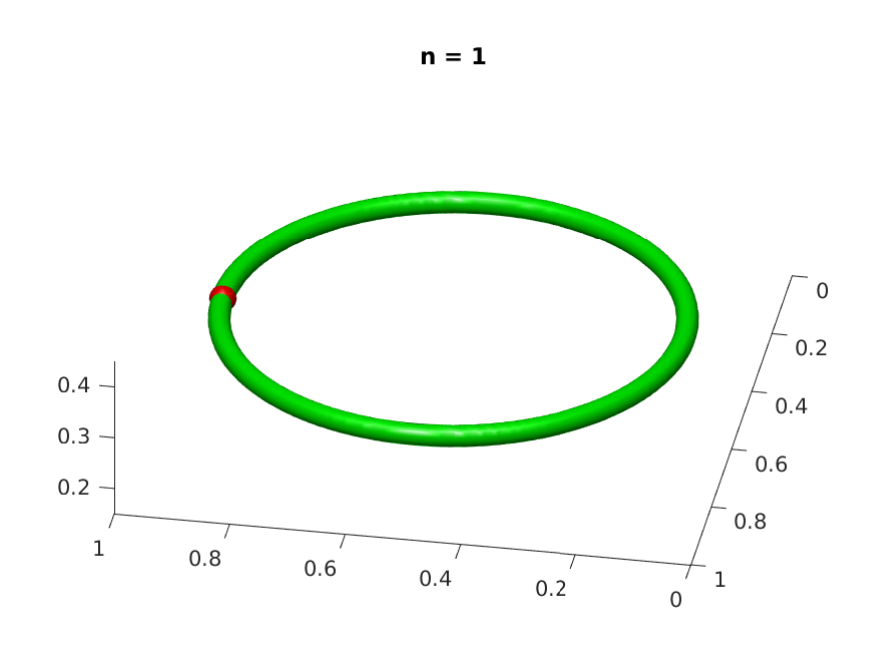}
		\includegraphics[width=.32\textwidth]{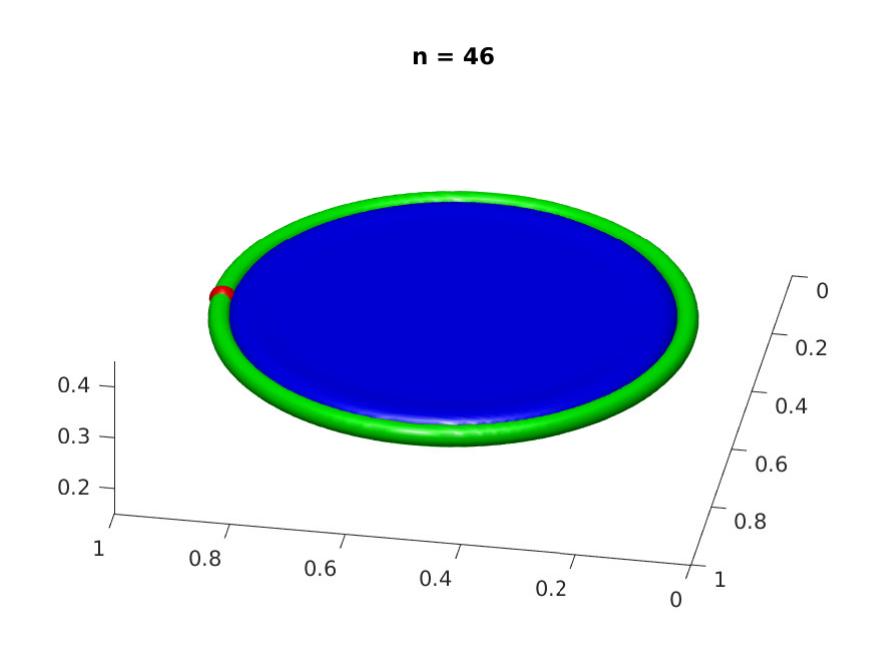}
		\includegraphics[width=.32\textwidth]{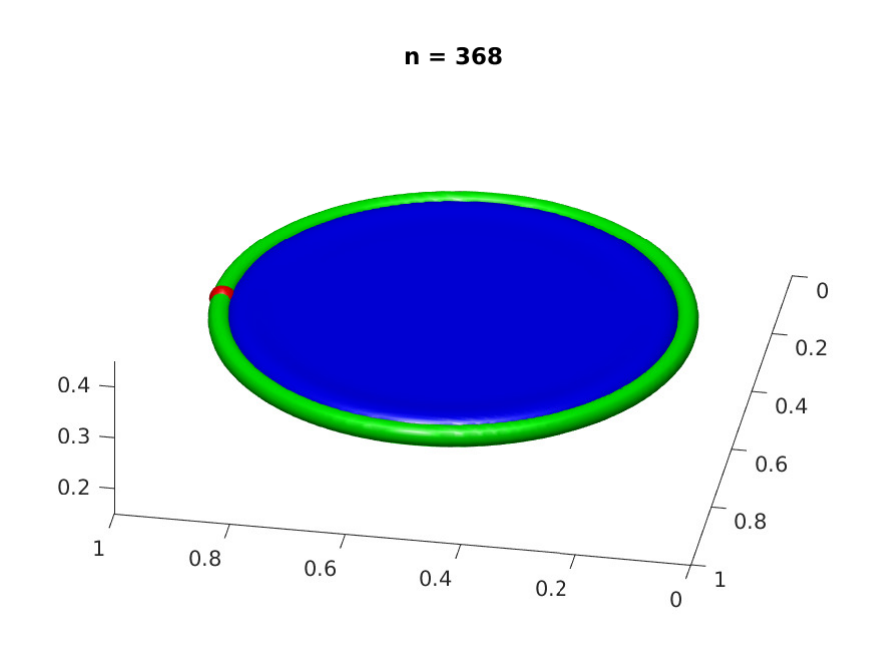} \\
		\includegraphics[width=.32\textwidth]{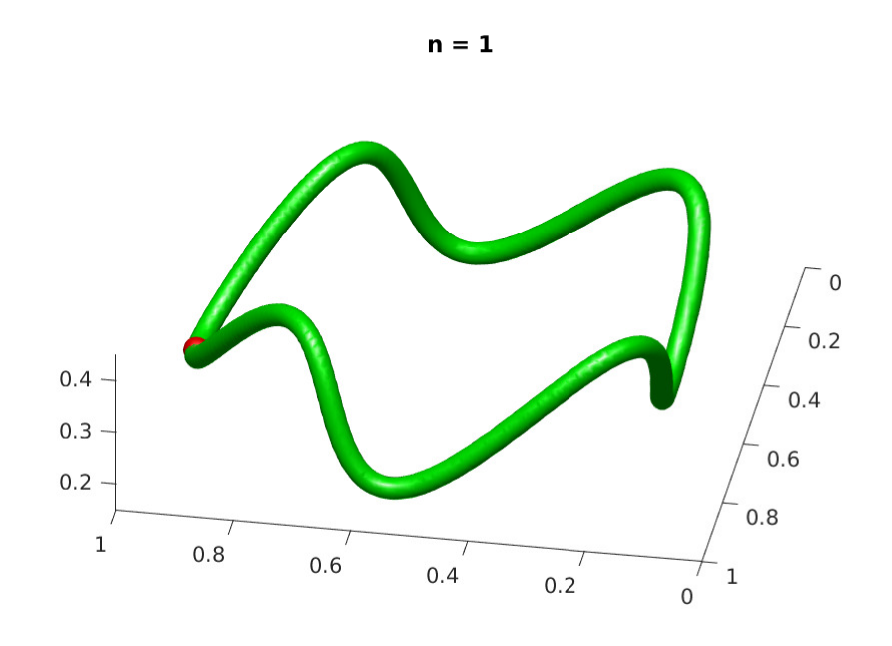}
		\includegraphics[width=.32\textwidth]{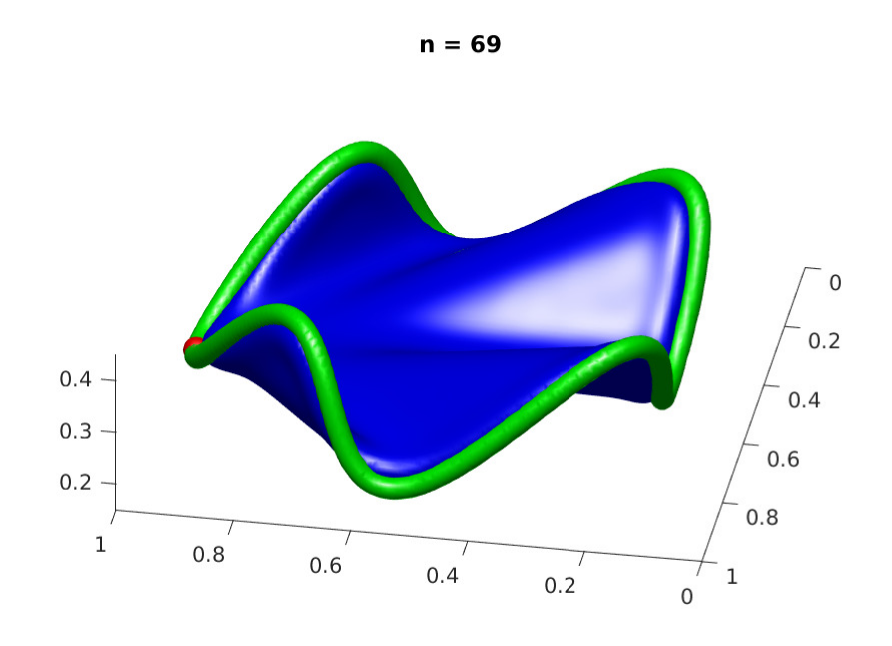}
		\includegraphics[width=.32\textwidth]{./Images/Plateau/Test_plateau_1curve_1_9-eps-converted-to}
		\caption{ Numerical approximations of Plateau's problem ; case $d=1$. In each line, we plot the solution  at different time $t_n$ along the iterations. In green, the curve $\Gamma^1$,
			in red, the point $\gamma^2 = \gamma^2(0)$ and in blue the surface $S_{\varepsilon}$. }
		\label{fig:plateau__onecurve1}
	\end{figure}

	The first two examples are shown in Figure~\ref{fig:plateau__onecurve1} and correspond to the case $d=1$ with a single $\Gamma^1$ curve positioned on a cylinder.

	We plot the approximation $S_{\varepsilon}$ of the minimal surface solution of Plateau's problem on each line at different stages of the optimisation process.
	Initially, the solution resembles the geodesic obtained when $u=0$. Then, under the influence of the Willmore Cahn-Hilliard term, it evolves until it reaches a stationary solution that is a good approximation of the Plateau solution. These first two examples correspond to the case studied in the Gamma-convergence analysis of our phase field model.
	
	\paragraph{Case of two curves.}
	
\begin{figure}[!htbp]
		\centering
		\includegraphics[width=.32\textwidth]{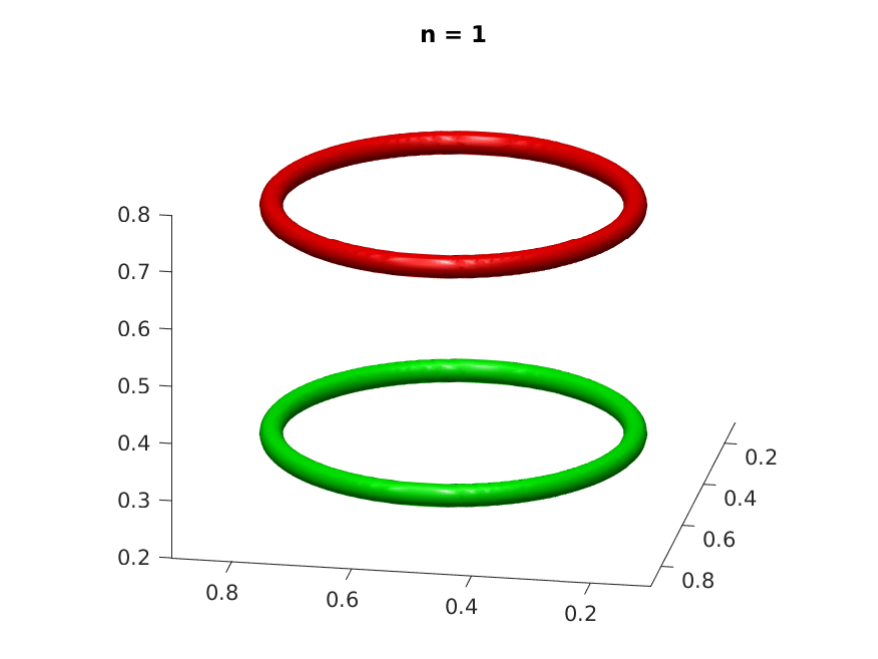}
		\includegraphics[width=.32\textwidth]{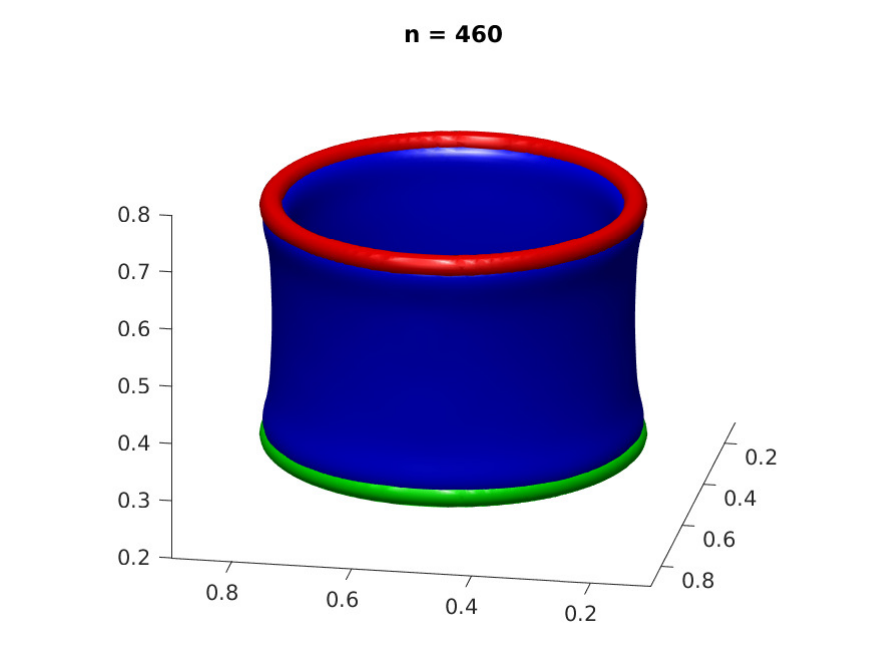}
		\includegraphics[width=.32\textwidth]{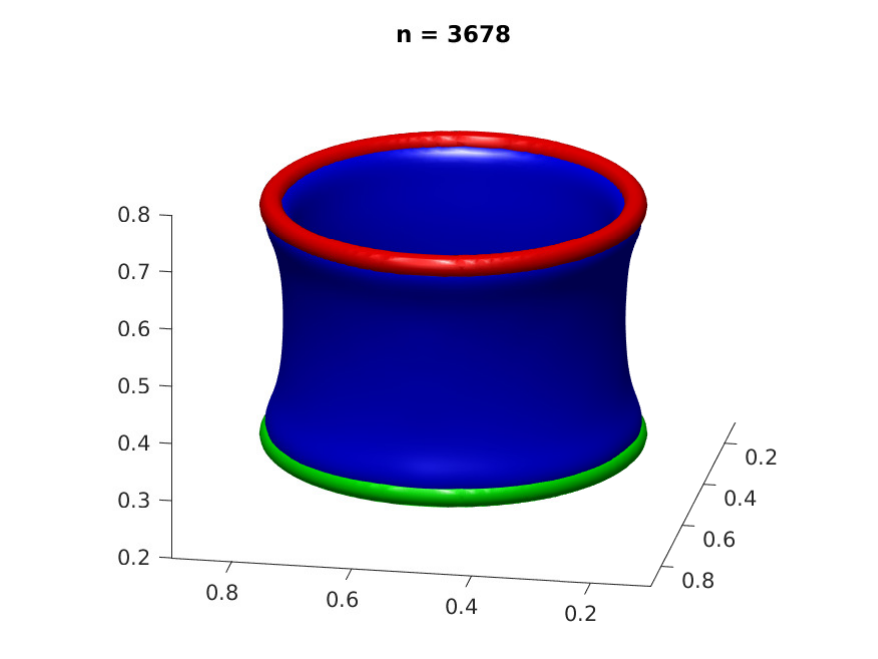} \\
		\includegraphics[width=.32\textwidth]{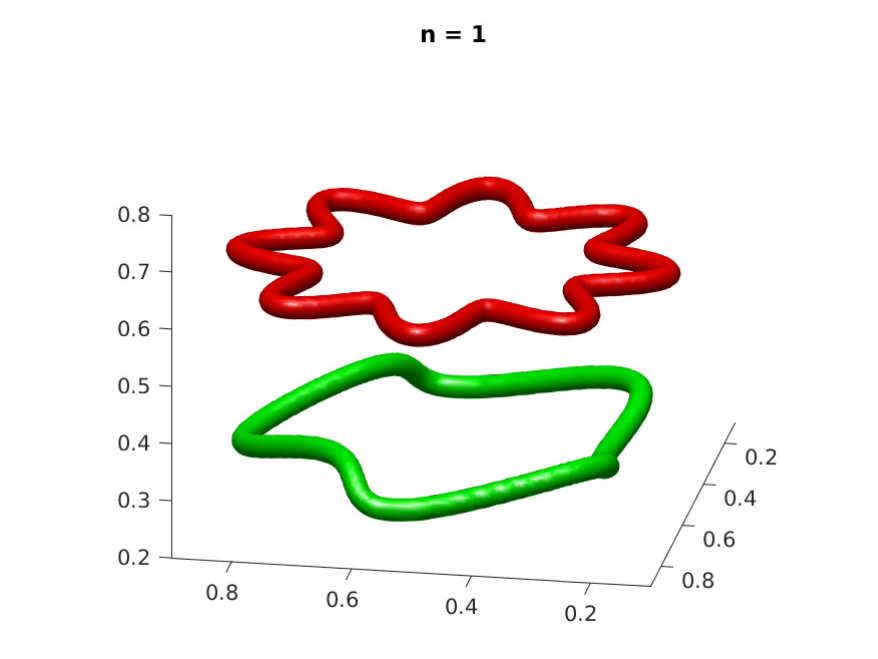}
		\includegraphics[width=.32\textwidth]{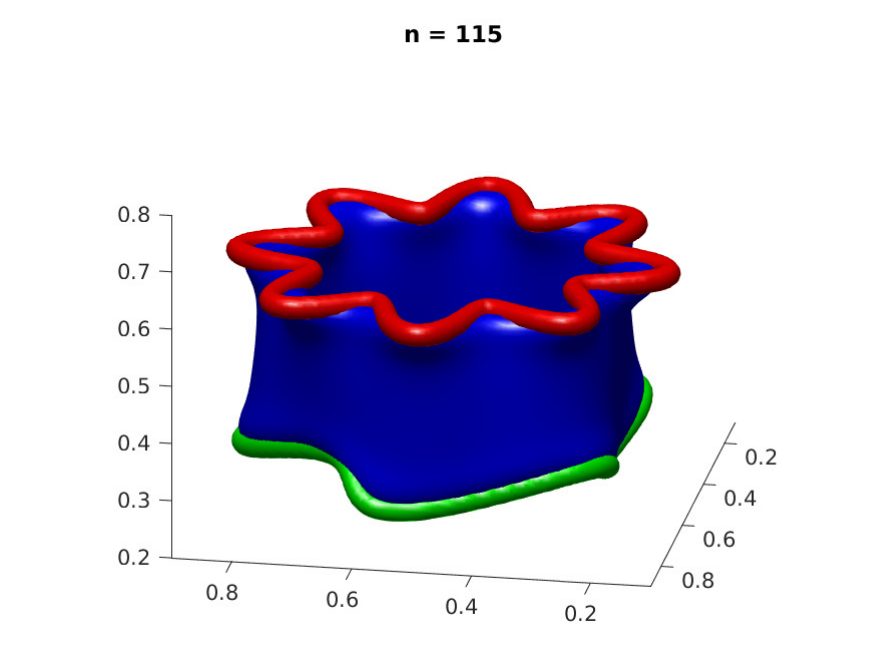}
		\includegraphics[width=.32\textwidth]{./Images/Plateau/Test_plateau_2circle_2_9-eps-converted-to} \\
		\includegraphics[width=.32\textwidth]{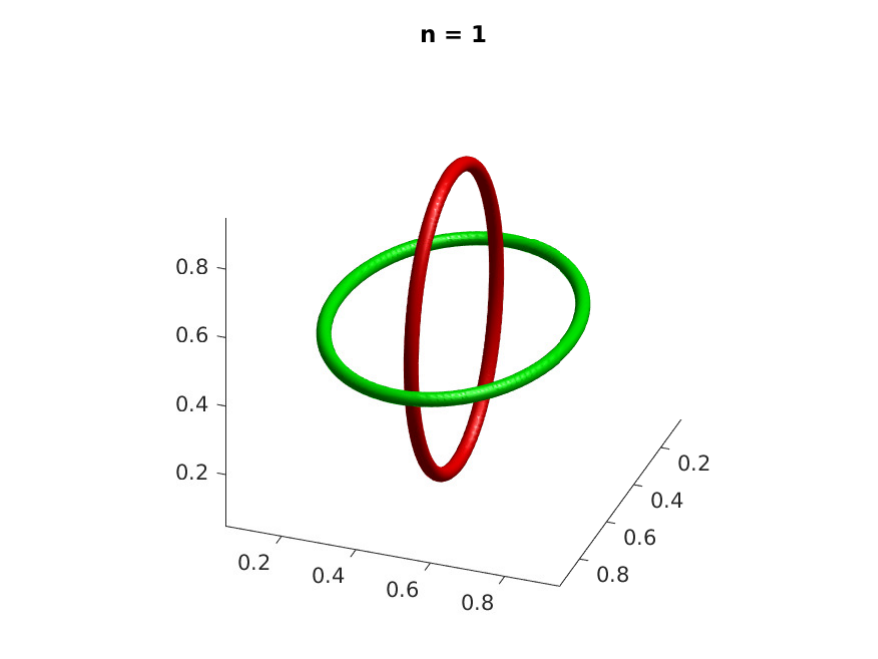}
		\includegraphics[width=.32\textwidth]{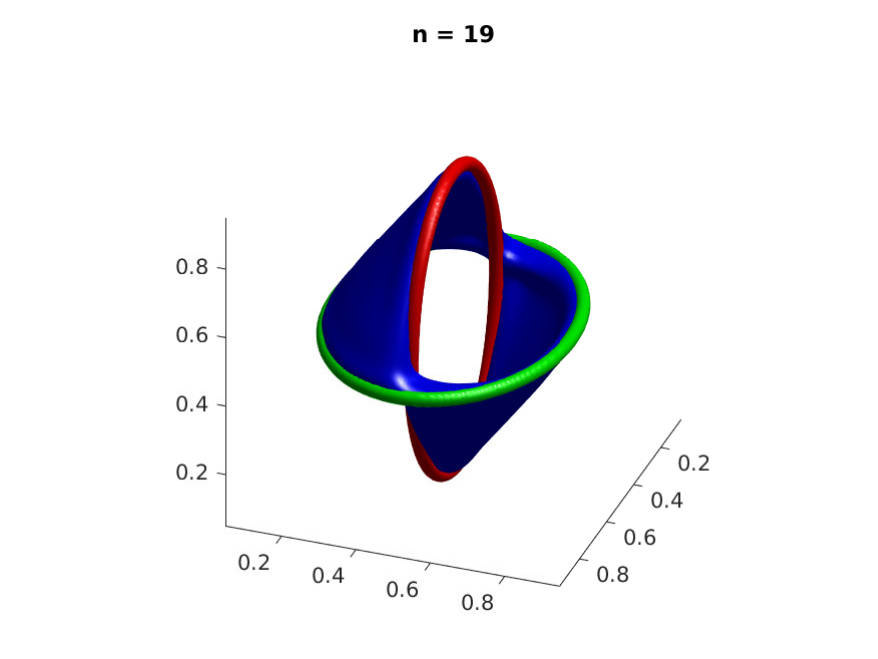}
		\includegraphics[width=.32\textwidth]{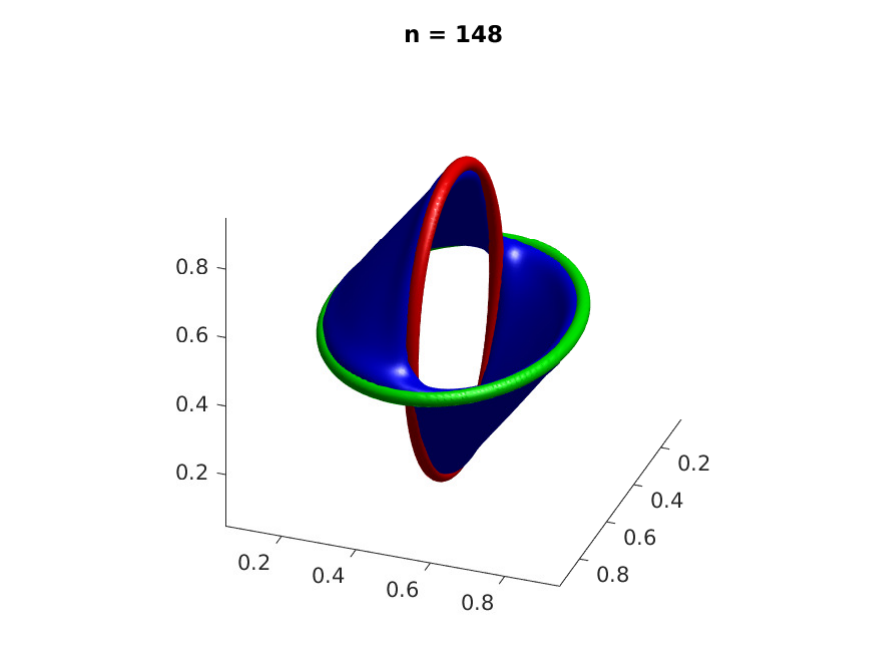} \\
		\includegraphics[width=.32\textwidth]{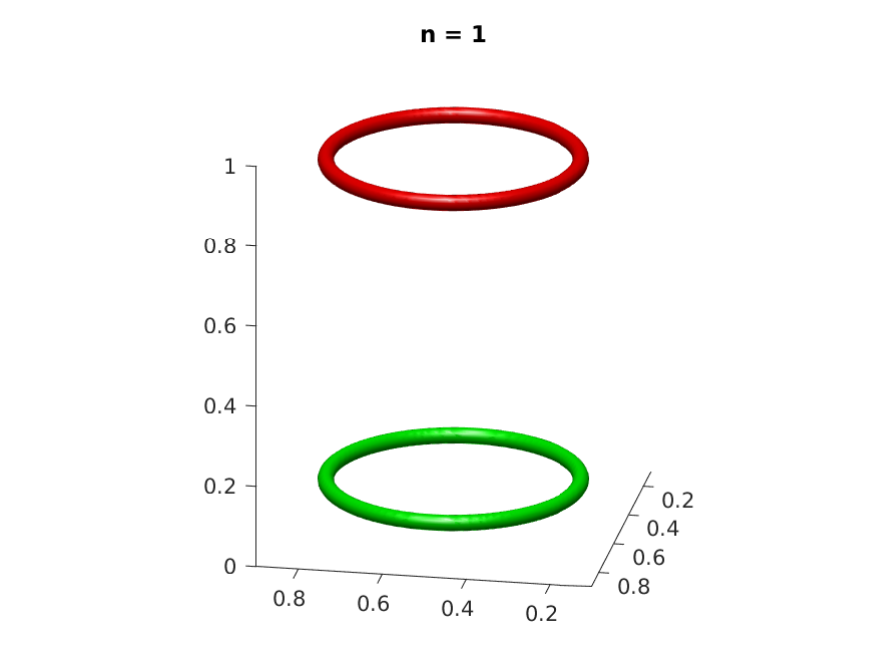}
		\includegraphics[width=.32\textwidth]{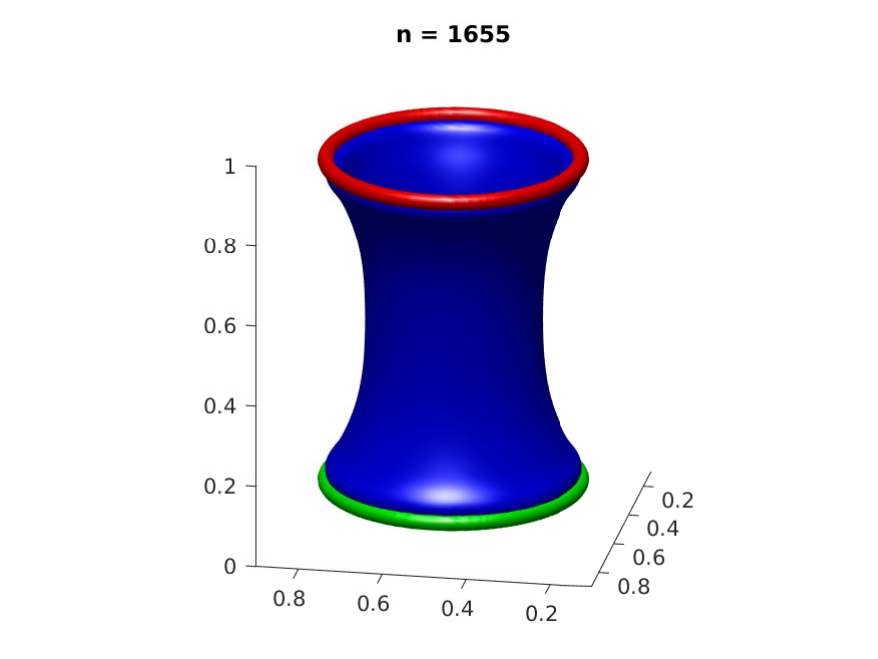}
		\includegraphics[width=.32\textwidth]{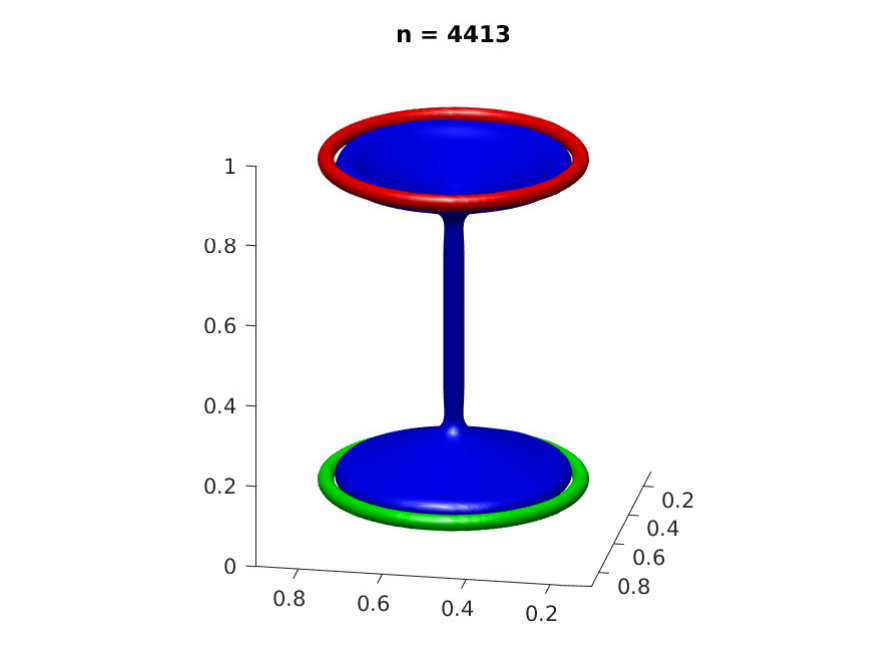}
		
		\caption{
			Numerical approximations of Plateau's problem ; $d=2$ with configuration $1$; $\Gamma^1= \gamma^1([0,2\pi])$ and  $ \Gamma^2 = \gamma^2[0,2\pi]$. In each line, we plot the solution  at different time $t_n$ along the iterations. In green, the curve $\Gamma^1$,
			in red, the curve $\gamma^2$ and in blue the surface $S_{\varepsilon}$.}
		\label{fig:plateau__twocurve1}
	\end{figure}

	Figure~\ref{fig:plateau__twocurve1} shows four new experiments in the case of two curves $\Gamma_1$ and $\Gamma_2$. Here we are in the first configuration of the case $d=2$ (described at the beginning of Section~\ref{Section:numerics}) where $\gamma^1$ and $\gamma^2$ represent a parametrization of the two curves $\Gamma_1$ and $\Gamma_2$ and where the geodesic term aims to  connect $\gamma^1$ and $\gamma ^2$ only. Hence, the energy $F_\eps$ is given by~\eqref{EnergyFepsConfig1}.

	The first line shows the case of two circles that are close enough together (in terms of their radii) that the Plateau solution is a catenoid.	Initially, the solutions resemble a cylinder corresponding to the geodesic approximation computed with $u = 0$, and then converge to the geodesic after $3000$ iterations of the algorithm. 
	
	The second example is very similar, and as previously, it can be seen that in the first iterations the solution resembles the geodesic approximation with $u=0$, before evolving towards a perturbation of the catenoid.
	
	The third example is of particular interest since the obtained solution is not oriented, and consequently does not fall within the framework of the analysis of Plateau's problem in a cylinder.  
	
	The fourth example considers circles that are far enough apart to prevent a catenoid from connecting them. In this case, the solution converges to two disks connected by a narrow tube. Note that the thickness of the tube depends on the parameter $\sigma_{\varepsilon}$. We can expect this tube to disappear when $\eps$ approaches $0$, which corresponds to two disks — the solution to Plateau's problem in such case.

\begin{figure}[!htbp]
		\centering
		\includegraphics[width=.32\textwidth]{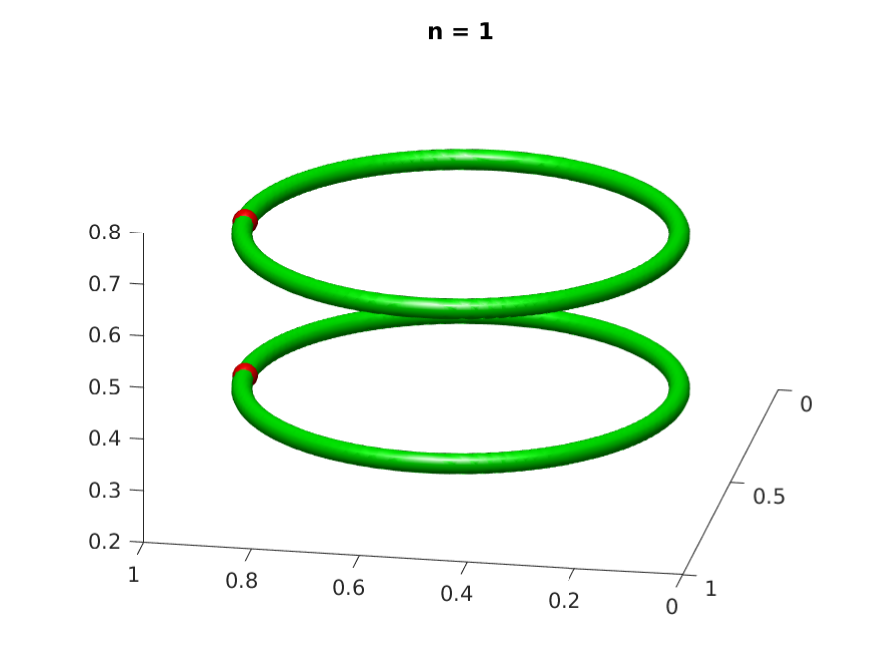}
		\includegraphics[width=.32\textwidth]{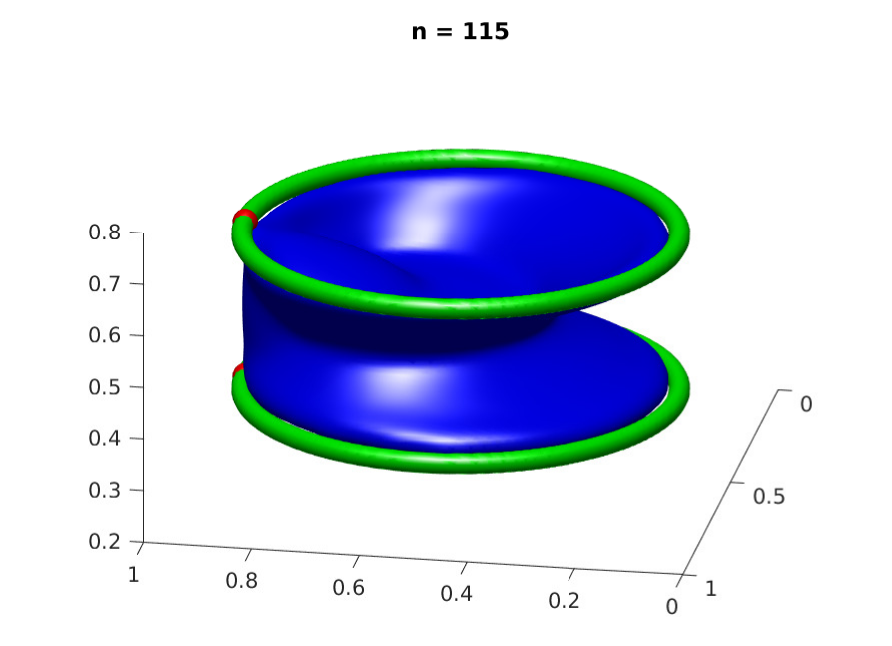}
		\includegraphics[width=.32\textwidth]{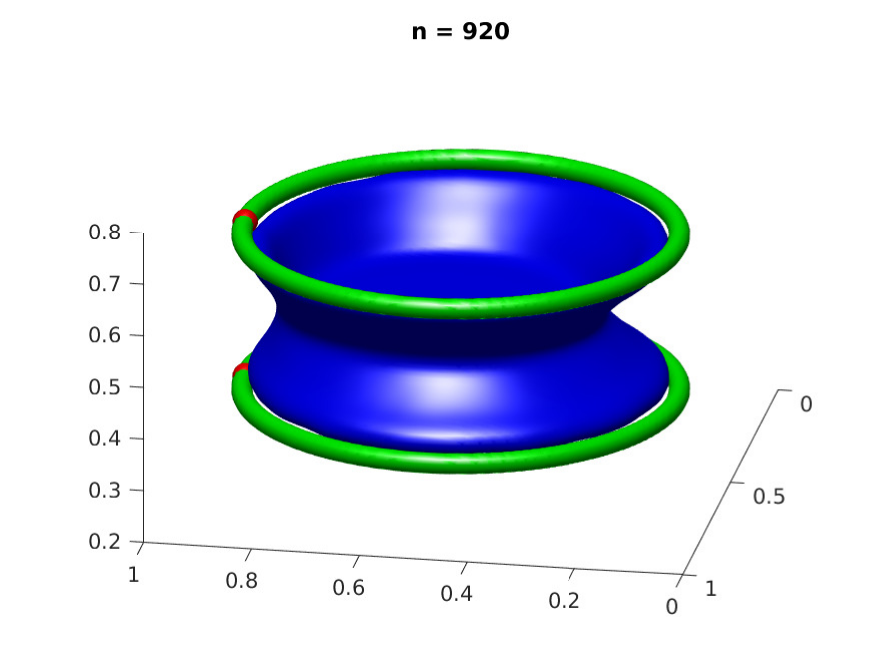} \\
		
		\caption{
			Numerical approximations of Plateau's problem ; $d=2$ with configuration $2$ ;  $\Gamma_1= \gamma^1([0,2\pi])$, $ \Gamma_2 = \gamma^2([0,2\pi])$, $\gamma^3 = \gamma^1(0)$ and $\gamma^4 = \gamma^2(0)$ . In each line, we plot the solution  at different time $t_n$ along the iterations. In green, the curves $\Gamma^1$ and $\Gamma^2$,
			in red, the points $\gamma^3$ and  $\gamma^4$ and in blue the surface $S_{\varepsilon}$.}
		\label{fig:plateau__twocurve2}
	\end{figure}

	\medskip
	Figure~\ref{fig:plateau__twocurve2} presents a numerical experiment
	of configuration $2$ with $d=2$: we consider  two circles 
	$\Gamma^1= \gamma^1([0,2\pi])$, $ \Gamma^2 = \gamma^2[0,2\pi]$ and two points
	$\gamma^3 = \gamma^1(0)$ and $\gamma^4 = \gamma^2(0)$. The geodesic penalization term connects $(\gamma^1,\gamma^3)$
	and $(\gamma^2,\gamma^4)$, and the energy $F_\eps$ is defined by~\eqref{EnergyFepsConfig2}.  
	
	In this situation, the catenoid is indeed a minimal surface. The point here is that this catenoid cannot be approached with the chosen penalty term alone. To do so, it must be supplemented by at least one disk, which is indeed the numerical solution found by our algorithm. This example therefore illustrates the influence of the choice of the geodesic term on the solution of Plateau's problem.

	\paragraph{Extension to multiple curves.}
	
\begin{figure}[!htbp]
		\centering
		\includegraphics[width=.32\textwidth]{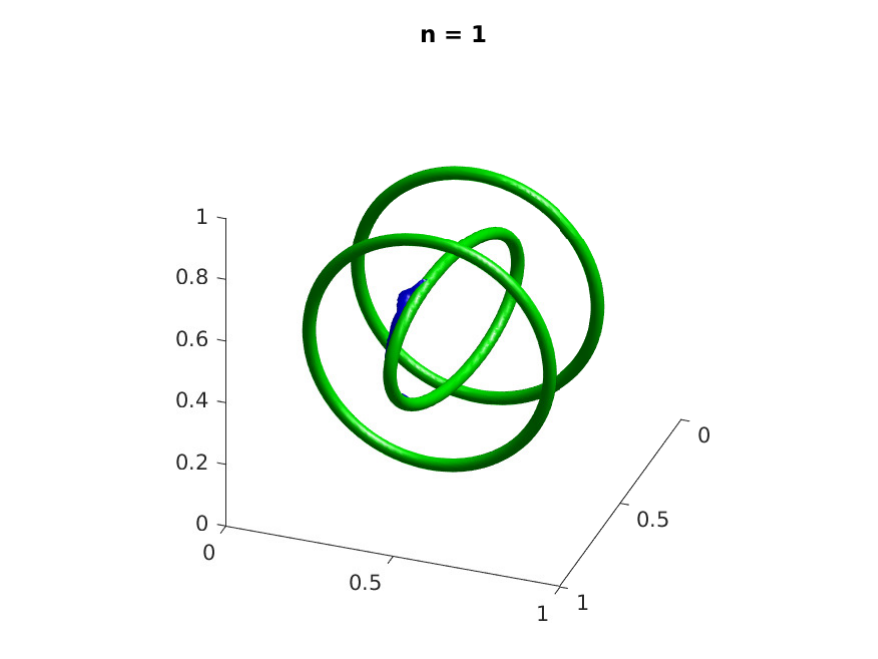}
		\includegraphics[width=.32\textwidth]{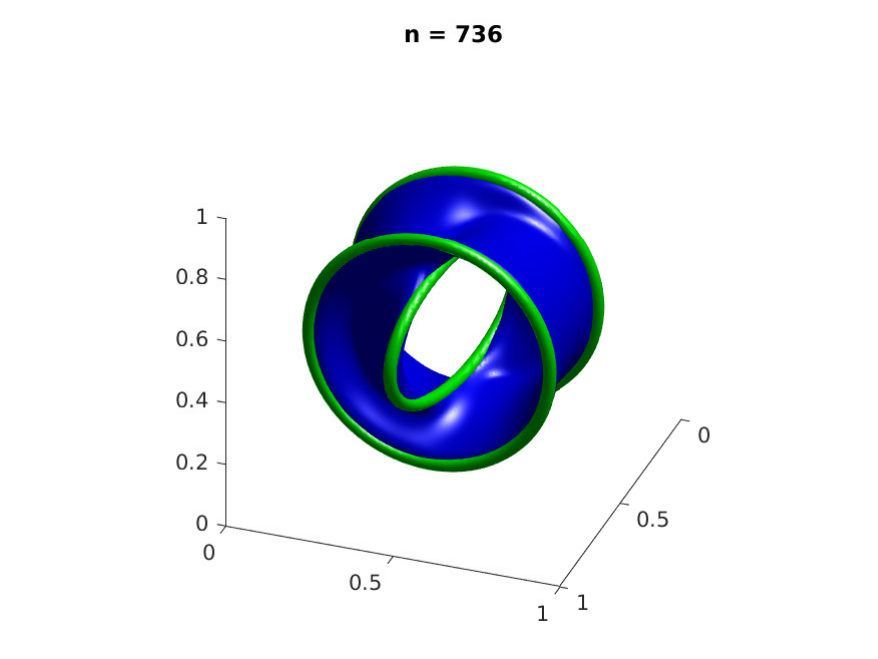}
		\includegraphics[width=.32\textwidth]{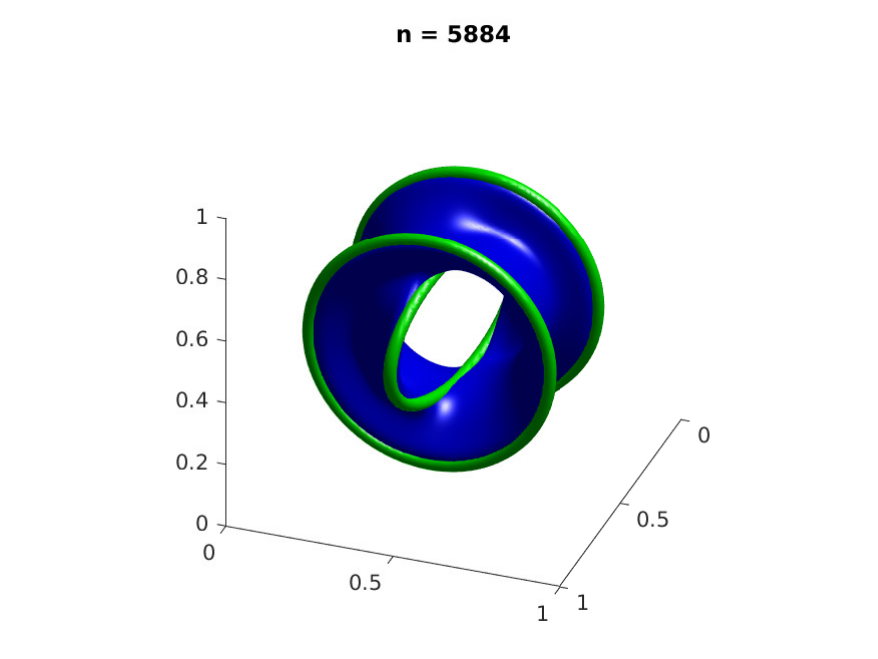} \\
		\caption{  Numerical approximations of Plateau's problem ; $d=3$ with configuration $1$ ;  $\Gamma^1= \gamma^1([0,2\pi])$, $ \Gamma^2 = \gamma^2([0,2\pi])$, $\Gamma^3 = \gamma^3([0,2\pi])$. In each line, we plot the solution  at different time $t_n$ along the iterations. In green, the curves $\Gamma^1$, $\Gamma^2$ and  $\Gamma^3$. In blue the surface $S_{\varepsilon}$.}
		\label{fig:plateau__3curve1}
	\end{figure}

\begin{figure}[!htbp]
		\includegraphics[width=.32\textwidth]{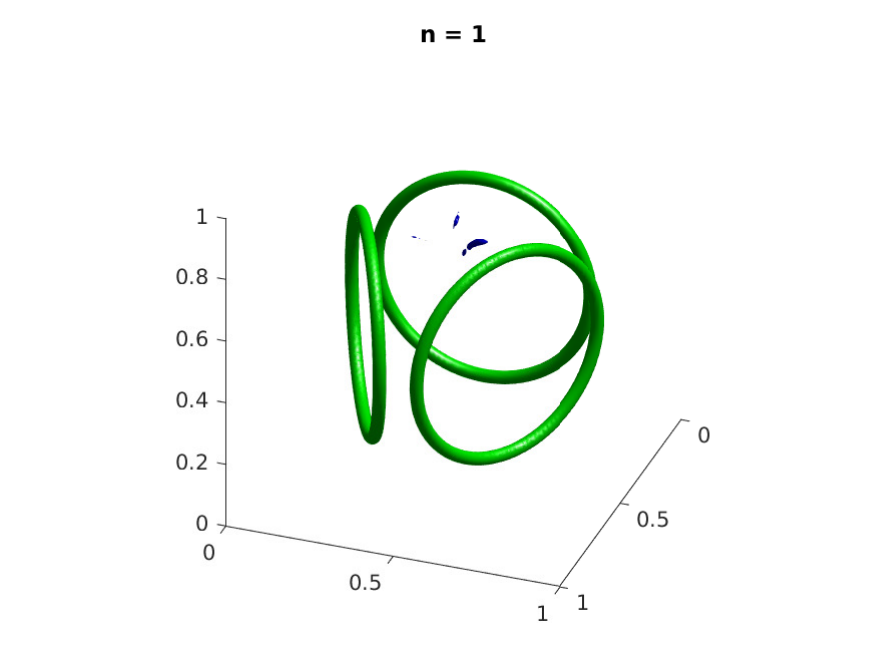}
		\includegraphics[width=.32\textwidth]{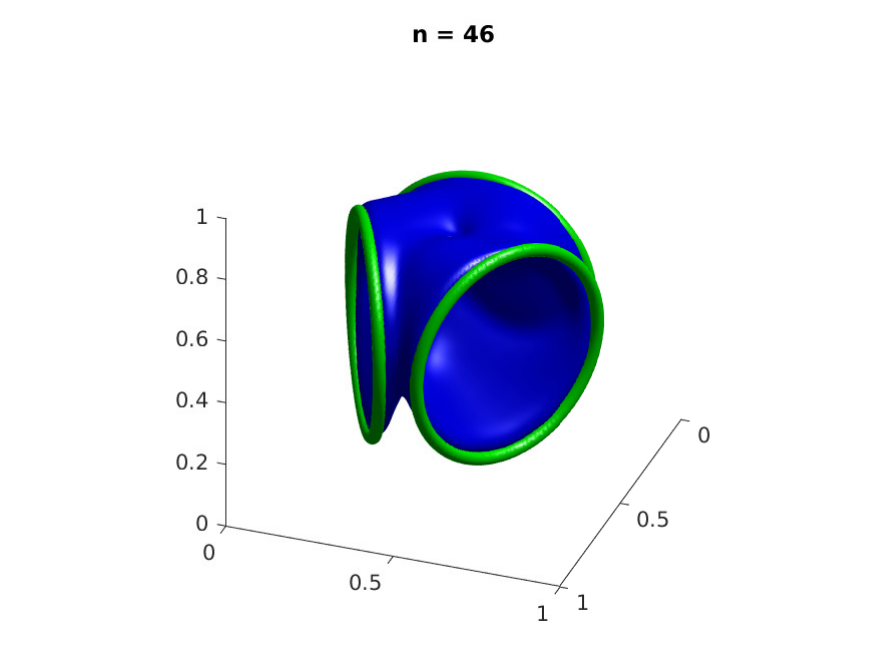}
		\includegraphics[width=.32\textwidth]{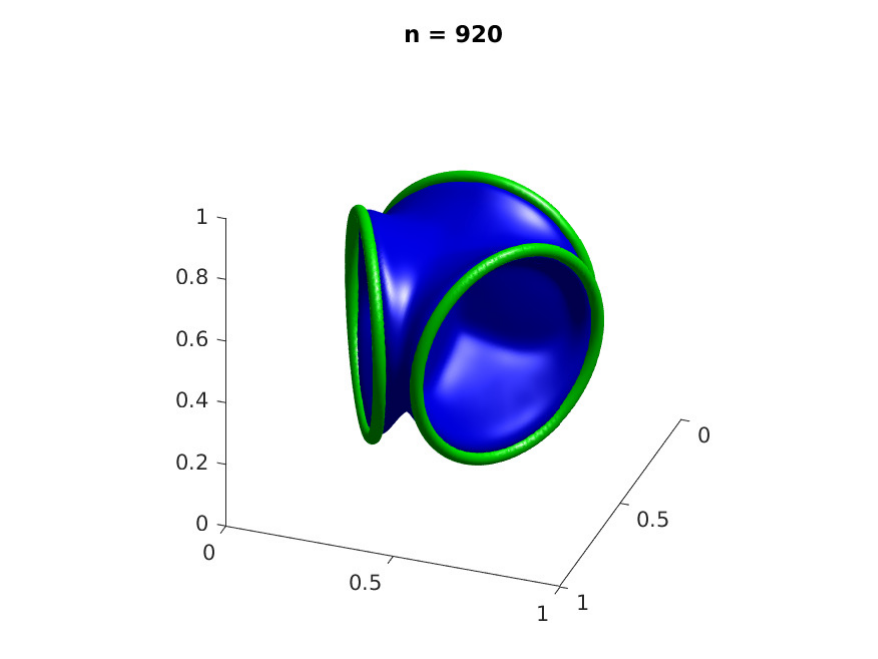}
		
		\caption{ Numerical approximations of Plateau's problem ; $d=3$ with configuration $2$ ;  $\Gamma^1= \gamma^1([0,2\pi])$, $ \Gamma^2 = \gamma^2([0,2\pi])$, $\Gamma^3 = \gamma^3([0,2\pi])$,  $\gamma^4 = \gamma^1(0)$,$\gamma^5 = \gamma^2(0)$ and $\gamma^6 = \gamma^3(0)$ . In each line, we plot the solution  at different time $t_n$ along the iterations. In green, the curves $\Gamma^1$, $\Gamma^2$ and  $\Gamma^3$. In blue the surface $S_{\varepsilon}$.}
		\label{fig:plateau__3curve2}
	\end{figure}

\begin{figure}[!htbp]
		\includegraphics[width=.32\textwidth]{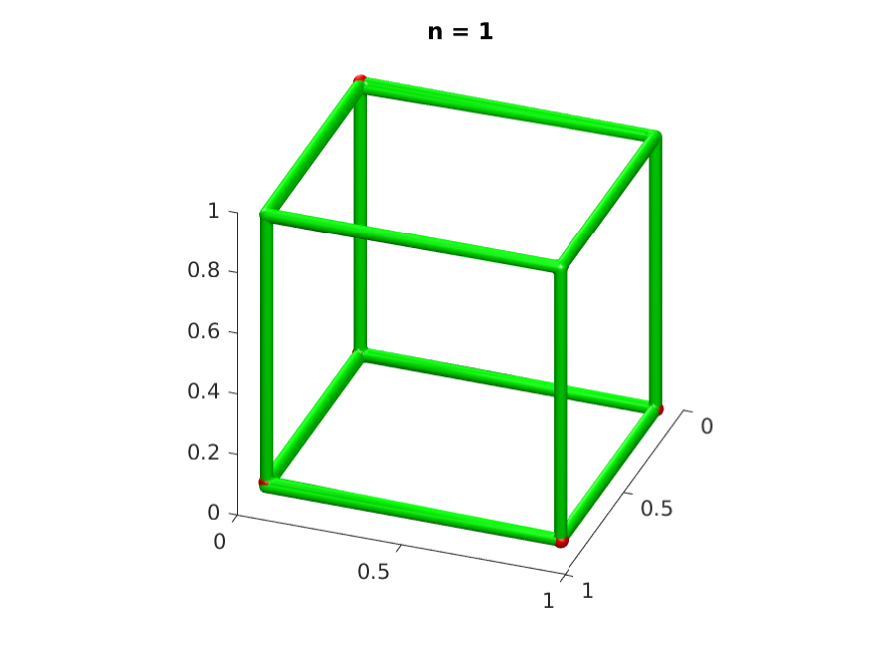}
		\includegraphics[width=.32\textwidth]{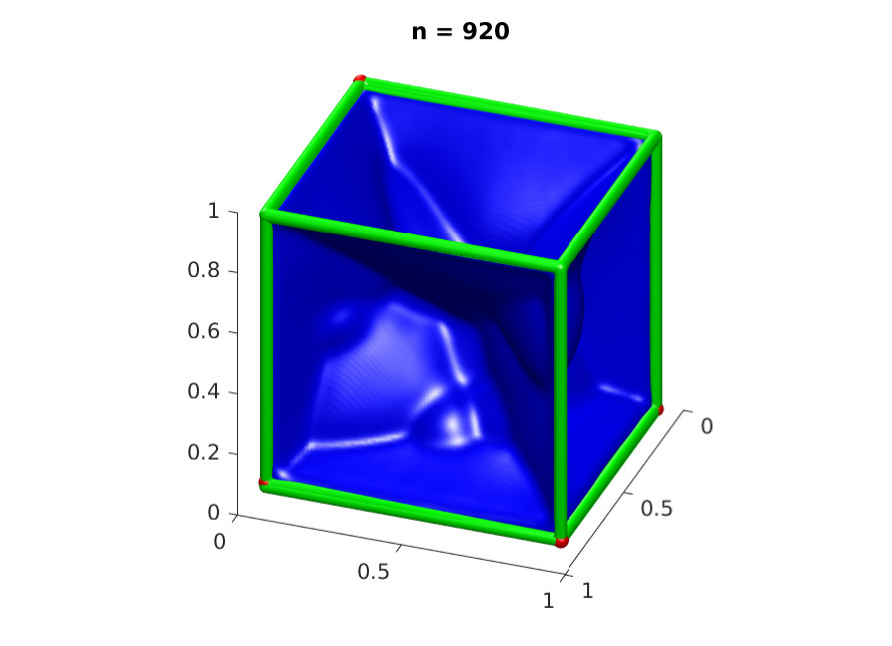}
		\includegraphics[width=.32\textwidth]{./Images/Plateau/Test_plateau_simple_cube_29-eps-converted-to}
		
		\caption{ Numerical approximations of Plateau's problem ; $d=6$ with configuration $2$. In each line, we plot the solution  at different time $t_n$ along the iterations. In green, the curves $\cup_{i=1}^{d}\Gamma^i$. In blue the surface $S_{\varepsilon}$.}
		\label{fig:plateau__Cube}
	\end{figure}

	The last three examples presented in Figures~\ref{fig:plateau__3curve1}, \ref{fig:plateau__3curve2} and \ref{fig:plateau__Cube} illustrate how our approach can easily be extended to more than two curves. In the first case, Figure~\ref{fig:plateau__3curve1}, we  use a geodesic term composed of the three terms associated with the pairs $(\gamma^1,\gamma^2)$,$(\gamma^1,\gamma^3)$ and $(\gamma^2,\gamma^3)$ where each curve $\gamma^i$ is associated  to $\Gamma^i$. This situation corresponds to configuration $1$ with  $d=3$. In the end, the algorithm does find a surface whose boundary corresponds to the set $\{ \Gamma^i \}_{i \in {1,2,3 }}$  and appears to be a minimal surface.
	The next example in Figure~\ref{fig:plateau__3curve2} also corresponds to the case of three circles but using a type $2$ configuration: we take $\gamma^1([0,2\pi])=\Gamma^1$, $\gamma^2([0,2\pi])=\Gamma^2$, $\gamma^3([0,2\pi])=\Gamma^3$,
	$\gamma^4 = \gamma^1(0)$, $\gamma^5 = \gamma^2(0)$ and $\gamma^6 = \gamma^3(0)$ and propose to penalize the distance between $(\gamma^1,\gamma^4)$,
	$(\gamma^1,\gamma^5)$, $(\gamma^2,\gamma^4)$,  $(\gamma^2,\gamma^6)$,
	$(\gamma^3,\gamma^4)$ and $(\gamma^3,\gamma^5)$. As before, the algorithm seems to converge to a stationary solution, which is in fact a Plateau solution to our problem. The last example shown in Figure~\ref{fig:plateau__Cube} corresponds to the classic cube example, for which we break down the given boundary $\Gamma$  as a union of six squares connected with a configuration of type $2$.  In particular, this latest numerical experiment highlights the ability of our phase field approach to find  complex and not smooth solutions to Plateau's problem.

	\appendix 
	
	\bibliographystyle{alpha}
	\bibliography{biblioarticle}

\end{document}